\documentclass[journal]{new-aiaa} 
\usepackage[utf8]{inputenc}
\usepackage{subcaption}
\usepackage{bm}
\usepackage{graphicx}
\usepackage{amsmath}
\usepackage[version=4]{mhchem}
\usepackage{siunitx}
\usepackage{longtable,tabularx}
\usepackage{textcomp}
\title{Classification and Feasibility Assessment of Infinitely Many Iso-Impulse Three-Dimensional Trajectories}
\author{Keziban Saloglu\footnote{Graduate Student, Department of Aerospace Engineering, Auburn University, AL 36849 US, AIAA Student Member.} and Ehsan Taheri\footnote{Assistant Professor, Department of Aerospace Engineering, Auburn University, Auburn, AL 36849, USA, AIAA Senior Member.}}
\affil{Auburn University, Auburn, AL, 36849}

\begin{document}

\maketitle

\begin{abstract}
In two-body dynamics, it is proven that for a sufficiently long flight time, generating infinitely many iso-impulse solutions is possible by solving a number of $\Delta v$-allocation problems analytically. A distinct feature of these solutions is the existence of two impulse anchor positions (APs) that correspond to the locations of the impulses on time-free, phase-free, base solutions. In this paper, the existence and utility of three-impulse base solutions are investigated and their complete solution spaces are characterized and analyzed. Since two- and three-impulse base solutions exist, a question arises: How many APs should base solutions have? A strategy is developed for choosing base solutions, which offers a certificate for $\Delta v$ optimality of general three-dimensional time-fixed rendezvous solutions. Simultaneous allocation of $\Delta v$ at two and three APs is formulated, which allows for generating $\Delta v$-optimal solutions while satisfying a constraint on individual impulses such that $\Delta v \leq \Delta v_\text{max}$. All iso-impulse solutions are classified in four layers: 1) base solutions, 2) feasible solution spaces, 3) solution families, and 4) solution envelopes. The method enables us to characterize the complete solution space of minimum-$\Delta v$, iso-impulse, three-dimensional trajectories under the nonlinear two-body dynamics.
To illustrate the utility of the method, interplanetary and geocentric problems are considered. 


\end{abstract}

\section{Introduction}
 The optimization of spacecraft trajectories in terms of flight time and propellant consumption is an essential step in space mission design ~\cite{strange2002graphical,ellison2018application, shirazi_spacecraft_2018}.
 In addition, many optimal (with respect to fuel consumption and flight time) trajectories are generated for determining the mass budget during the preliminary spacecraft design phases \cite{woolley2019optimized,woolley2019cargo}. Spacecraft trajectories are significantly affected by the type of their propulsion systems, which can be classified broadly as chemical rockets \cite{olympio2009designing}, solar-powered electric propulsion systems \cite{nurre2023duty,petukhov2019joint} and propellant-free systems (e.g., solar sails) \cite{pengyuan2023design}. 
 
Practical spacecraft trajectory optimization problems are solved using numerical methods that are classified broadly as direct and indirect methods \cite{betts_survey_1998}. Indirect methods \cite{bertrand2002new,zhu2017solving,taheri2018generic,taheri2021costate,li2021homotopy} form a (problem-dependent) set of necessary optimality conditions and find multiple solutions to them through numerical continuation and homotopy methods for impulsive \cite{olympio2007global,shen2014indirect} and low-thrust trajectories \cite{kluever1995optimal,taheri2016enhanced,pan2020practical}. For impulsive trajectory optimization, indirect-based methods refer to the methods that use the primer vector theory (PVT) based approaches \cite{lawden_optimal_1963}.
 The direct methods, however, discretize continuous-time optimization problems into parametric optimization ones~\cite{conway_survey_2012} and incorporate various types of state-only and mixed-state-control path constraints \cite{englander2017automated,mannocchi2022homotopic,nurre_comparison_2022,oshima_regularizing_2024}. The quality of the initial guess influences the convergence performance of the nonlinear programming (NLP) problem solvers \cite{conway_survey_2012, ottesen_unconstrained_2021,sowell2024eclipse}. 
 
 Impulsive trajectory optimization problems are important for crewed missions, as they require shorter mission times compared to low-thrust trajectories~\cite{shirazi_spacecraft_2018}. Determining when and where to apply a finite, yet unknown number of impulses, their magnitude and direction is a fundamental task \cite{edelbaum_how_1966}. Early work on impulsive trajectories often categorizes solutions concerning initial and target orbit shape analytically~\cite{gobetz_survey_1969}. When the terminal orbits are sufficiently close, the Hohmann transfer produces minimum-$\Delta v$ solutions for planar circular-to-circular and coaxial elliptical-to-elliptical transfers with two tangential impulses at the apsides~\cite{hohmann_attainability_1960, curtis2013orbital}. If the terminal orbits are sufficiently far apart, the bi-elliptic transfer with three impulses requires less $\Delta v$ than the Hohmann transfer~\cite{edelbaum1959some}. The bi-parabolic transfers characterize the theoretical limit of time-free, three-impulse trajectories with a midcourse, zero-magnitude impulse at infinity~\cite{vallado_fundamentals_2001}. It is shown that for transfers between coaxial planar ellipses, at most, three impulses are required~\cite{ting1960optimum}. The surveys on early work for impulsive trajectory optimization can be found in Refs.~\cite{edelbaum_how_1966, gobetz_survey_1969}. 
 Recent investigations include the trajectories of the planar elliptic transfers. 
 Carter and Humi proposed the closed-form solutions of optimal two- and three-impulse transfers when the terminal ellipses are coaxial for the planar case~\cite{carter_two-impulse_2020}. They also used a transformation to have linear equations of motion to investigate the time-free and time-fixed transfers between planar ellipses with two and three impulses~\cite{carter_generalized_2021}. The non-coaxial planar elliptical transfers are solved by a simplified optimization algorithm based on a homotopic approach, starting from a two-impulse Hohmann solution using the PVT~\cite{li_minimum-fuel_2024}. 
 
 The PVT constitutes a set of first-order necessary conditions for the local optimality of impulsive solutions. It offers a mechanism to determine the number of impulses, impulse locations, and directions. Ref.~\cite{jezewski_efficient_1968} demonstrates using the PVT combined with an NLP-based approach to obtain optimal impulsive solutions. 
 Obtaining an extremal solution based on the PVT starts with solving two- or many-impulse problems with the formulated NLP problem. Once a solution is obtained, the primer vector time history is inspected to ensure that the necessary conditions are satisfied. If there are one or multiple violations of the necessary conditions, an impulse is inserted at that time instant and the problem is resolved~\cite{sandrik_primer-optimized_nodate, bell_primer_nodate, chiu_optimal_nodate, bokelmann_halo_2017}. PVT is used for impulsive rendezvous between close, circular, co-planar orbits for time-fixed problems~\cite{ prussing_optimal_1986,prussing_optimal_2003}. Hughes et al. compared the direct and primer-vector-based indirect optimization methods and concluded that the methods are not superior to each other~\cite{hughes_comparison_2003}. It is also shown that the PVT can be used for efficient impulsive maneuver placement in complex trajectory optimization problems \cite{landau2018efficient}.
 

Impulsive maneuvers can be viewed as the limiting case of low-thrust trajectory optimization problems when the magnitude of the thrust is sufficiently large \cite{gergaud2007orbital,taheri_how_2020}. This fundamental connection between impulsive and finite-thrust solutions is leveraged successfully to impulsive trajectory optimization problems under two-body dynamics \cite{arya_generation_2023}, and recently, within the circular restricted three-body problem models \cite{saloglu2024acceleration}.


The connection between the low-thrust and impulsive solutions is established with the notion of optimal switching surfaces, introduced by Taheri and Junkins~\cite{taheri_how_2020} to answer the ``How Many Impulses?'' question posed by Edelbaum~\cite{edelbaum_how_1966}. More specifically, they reported that for the problem of transferring a spacecraft from the Earth to asteroid Dionysus (Earth-to-Dionysus problem), four impulsive solutions require the exact same total $\Delta v$. Therefore, it is shown that Edelbaum's ``How Many Impulses?'' question does not have a unique answer. Even though they considered time-fixed, rendezvous-type classes of maneuvers, there are some key features in the solutions of the Earth-to-Dionysus problem, including the fact that all impulses are applied at two distinct locations and that there is a fundamental arc common between all reported equal-$\Delta v$ solutions. Later, the common features of the equal-$\Delta v$, multiple-impulse solutions of the Earth-to-Dionysus problem are investigated and explained by introducing a method to generate impulsive solutions by Saloglu et al.~\cite{saloglu_existence_2023, saloglu_existence_conf_2023}. 
It is shown that, for long-time-horizon missions, impulsive trajectory optimization problems can be solved starting from a two-impulse, time-free, phase-free base solution, where impulse anchor positions (APs) are the positions at which the impulses occur on the base solution~\cite{saloglu_existence_2023}. At these impulse APs, impulses can be divided into multiple smaller-magnitude impulses. This way, thruster requirements can be satisfied without trading the $\Delta v$ optimality, and in some cases, only at the expense of requiring a maneuver over a longer flight time. In essence, these smaller-magnitude impulses inject the spacecraft into several phasing orbits. Similarly, the phasing of the spacecraft on the geostationary orbit can be accomplished without considering additional phasing maneuvers~\cite{barea2022dual2}.
In Saloglu et al.~\cite{saloglu_existence_2023}, the problem of dividing the impulse at an impulse AP is formulated as a $\Delta v$-allocation problem, which is solved analytically. This method represents a substantially faster alternative to the continuation-based methods \cite{arya_generation_2023}. 

Our previous study \cite{saloglu_existence_2023} has been the motivation to develop and propose a generalized approach to the design of minimum-$\Delta v$, iso-impulse maneuvers. The time-free, three-impulse transfer maneuvers are relevant in spacecraft trajectories and appear frequently in circular-to-circular or elliptic-to-elliptic transfers. In addition, there exist impulsive trajectories that require less $\Delta v$ for three-impulse solutions, similar to bi-elliptic transfers. To generalize the $\Delta v$-allocation problem, three-impulse base solutions have to be considered in addition to the two-impulse base solutions. Moreover, there can be thruster constraints to divide impulses into smaller ones at multiple APs such that individual impulses satisfy the $\Delta v \leq \Delta v_\text{max}$ constraint with $\Delta v_\text{max}$ denoting a mission-specific maximum $\Delta v$ threshold. An ideal $\Delta v$-allocation formulation has to take this constraint into the problem formulation while retaining $\Delta v$ optimality.  

The contributions of this paper are as follows. 1) We consider maneuver scenarios with three-impulse, time-free, phase-free base solutions for generating infinitely many iso-impulse minimum-$\Delta v$ solutions. Classification of the solutions is presented in terms of the number of phasing orbits at possibly all APs. These solution families are represented with solution envelopes in the form of polygons~\cite{saloglu_existence_2023}.
Similar to bi-elliptic transfers between coplanar circular and coaxial elliptical orbits, one has to consider the possibility of having lower-$\Delta v$ solutions compared to the two-impulse base solutions. By considering two- and three-impulse base solutions, we can cover all minimum-$\Delta v$ solution possibilities in a general, unified manner. This contribution is based on our recent work \cite{saloglu2024new}, but we provide additional details and examples compared to the conference paper. 2) In \cite{saloglu_existence_2023}, we introduced the $\Delta v$-allocation problem and applied it only to a single AP that requires the maximum amount of impulse among the potential impulse APs. In this paper, we extend the $\Delta v$-allocation problem by deriving the relations that encompass two and three APs simultaneously. This generalization allows us to divide the impulses simultaneously at any potential AP under a time-feasibility criterion. We extend the time-feasibility relation to base solutions when multiple APs are considered; we further visualize feasible solution spaces based on the total number of revolutions that can occur on all possible phasing orbits. Consequently, we can extend the solution space and enable the option of dividing the impulse magnitude at all impulse APs. From a practical point of view, one can reduce all individual impulse magnitudes and make the entire transfer realizable by thrusters such that individual impulses satisfy the $\Delta v \leq \Delta v_\text{max}$ constraint. 
However, for problems with smaller values for $\Delta v_\text{max}$, it is inevitable to extend the time of maneuver if $\Delta v$-optimality is to be preserved, especially for fixed-time cases. 3) If we consider transfer maneuvers between non-coaxial, three-dimensional, elliptical orbits, a fundamental question arises: How many APs should the base solution have? To address this question, we propose a strategy for selecting between two- and three-impulse base solutions, which provides a certificate for $\Delta v$ optimality of general three-dimensional impulsive solutions. This certificate and its associated time-feasibility criterion offer an elegant mechanism to modify locally suboptimal time-fixed transfer and time-fixed rendezvous problems into theoretically optimal ones. For the maneuvers between coaxial, inclined circular orbits, we provide the details to determine the optimality of all two- and three-impulse base solutions. The optimality of two- and three-impulse solutions between non-coaxial, planar, elliptical orbits exist in the literature \cite{pontani_simple_2009,li_minimum-fuel_2024}. However, under the hypothesis that transfer between general three-dimensional, non-coaxial, and elliptical orbits can only have at most three impulses (excluding the limiting bi-parabolic ones), our proposed method can be used for turning suboptimal solutions into theoretically optimal ones for the complete set of two- and three-impulse bases solutions. 4) We show that all iso-impulse solutions can be characterized and classified in four layers: a) two- or three-impulse base solutions, b) feasible solution spaces (in the form of polytopes and polygons), c) solution families (characterizing possible orbital revolutions), and d) solution envelopes (characterizing orbital periods ranges). Characterization and significance of each layer are explained for two- and three-impulse base solutions, which extends our ability to characterize the complete solution space of minimum-$\Delta v$, iso-impulse, three-dimensional solutions under the nonlinear two-body dynamics.

The paper is organized as follows. In Section~\ref{sec:genmultimp}, the $\Delta v$-allocation problem for a single impulse AP is briefly reviewed. 
The generation of three-impulse base solutions is explained in Section~\ref{sec:threeimpbase} and details for selecting two- or three-impulse base solutions are explained. Feasibility equations for multiple impulse APs for having infinitely many solutions are given in Section~\ref{sec:multianchor}.  In Section~\ref{sec:results}, numerical examples are given for generating interplanetary and planet-centric minimum-$\Delta v$, iso-impulse solutions. Discussions on the selection of the base solutions by numerical examples are given in Section~\ref{sec:discussion_base}. Finally, conclusions are presented in Section \ref{sec:conclusion}. Detailed analytical relations to produce solution envelopes are given in the Appendix. 

\section{Generation of Infinitely Many Solutions from Two-Impulse Base Solutions}\label{sec:genmultimp}
\begin{figure}[htbp!]
    \centering
    \vspace{-4mm}
    \includegraphics[width=0.6\columnwidth]{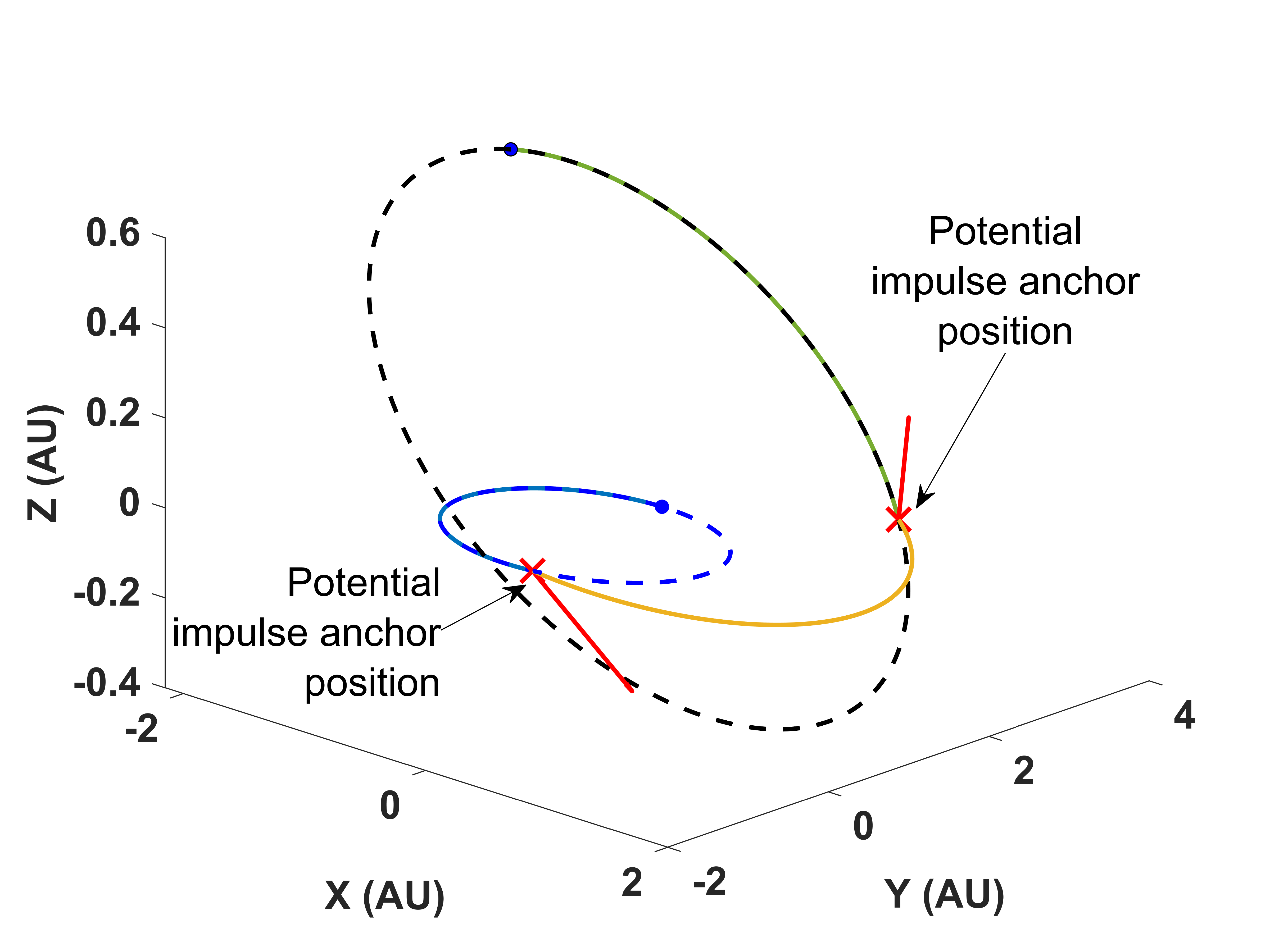}
    \caption{Earth-to-Dionysus phase-free two-impulse base solution and potential impulse APs ($t_{c_1} \neq 0, t_{c_2} \neq 0$).}
    \label{fig:twoimp_base}
\end{figure}

We review two-impulse base solutions but limit our analysis to only one of the APs~\cite{saloglu_existence_2023}. The definitions are common among the two- and three-impulse base solutions, which facilitates the discussion of the more complex cases. 
Let $\Delta \bm{v}_i$ and $\Delta \bm{v}_f$ denote the initial and final impulse vectors on a time-free, phase-free base solution applied at two impulse locations, respectively. Minimization of the total impulse, $\Delta v_\text{total}$, subject to the two-body dynamics, can be stated as,
\begin{align}  \underset{\theta_i,~\theta_f,~t_\text{pf}}{\text{minimize}}\quad \Delta v_\text{total} = ||\Delta \bm{v}_i|| + ||\Delta \bm{v}_f||,  
    \label{eq:biimpopt}
\end{align}
where $\theta_i \in [0,2\pi]$ and $\theta_f\in [0,2\pi]$ are the initial and target orbit true anomalies, respectively, and $t_\text{pf}$ denotes the time spent on the phase-free (subscript `pf') connecting arc between the two impulses and $||. ||$ denotes the Euclidean norm. The locations of the impulses denote the potential impulse APs. One can choose any or both of the impulse APs to divide the total impulse amount by introducing phasing orbits without sacrificing $\Delta v$-optimality under a time-feasibility criterion, which is discussed in Section~\ref{sec:feasTime}. An example two-impulse base solution is shown in Fig.~\ref{fig:twoimp_base} for the benchmark Earth-to-Dionysus problem \cite{taheri_how_2020}. The connecting arc is colored yellow. 
The blue and green solid arcs denote two coast arcs from the Earth to the AP on the Earth orbit and from the AP on the orbit of Dionysus to the location of Dionysus. The cross markers are the impulse locations and the impulses are denoted as red vectors. Let $t_{c_1}$ and $t_{c_2}$  denote the duration of the terminal coast arcs, which is determined by calculating the time from departure or arrival points to the impulse locations on the connecting arc. AU stands for the astronomical unit. 

If the total mission time is fixed (which is the case for time-fixed rendezvous maneuvers), there is a constraint on the available time of flight, $TOF$. We highlight the difference between the total mission time, denoted as $t_f - t_0$, and the available time of flight, $TOF$, which can be written as,
\begin{equation} \label{eq:Tava2imp}
    TOF = \underbrace{t_f-t_0}_{\text{total mission time}} - \underbrace{(t_{c_1}+t_\text{pf}+t_{c_2})}_{\text{time spent on base coast arcs}},
\end{equation}
where $t_f$ and $t_0$ correspond to the dates used for determining the states of the Earth and Dionysus on their orbits using the ephemerides, respectively. Let $T_0$ and $T_f$ denote the orbital periods of the initial and target orbits (Earth and Dionysus in this case), respectively, with $T_0 < T_f$. We have $t_{c_1} \in [0,T_0]$ and $t_{c_2} \in [0, T_f]$ since the coast arc values on the initial and final orbits can not be larger than the orbital periods. For the Earth-to-Dionysus problem with $t_f - t_0 = 3534$ days, the following values are obtained \cite{saloglu_existence_2023}: $\theta_i = 179.27^\circ$, $\theta_f = 149.20^\circ$, $t_{c_1} = 193.24$ days, $t_\text{pf} = 348.46$ days, and $t_{c_2} = 501.81$ days. $\Delta v_\text{total} = 9.907425$ km/s. Substituting the values for all coast times in Eq.~\eqref{eq:Tava2imp}, we have $TOF = 2490.48$ days. We emphasize that $t_f - t_0 = 3534$ days is based on the benchmark problem statement \cite{taheri2016enhanced}.

The available time of flight, $TOF$, can be used to introduce additional intermediate phasing orbits and multiple revolutions on the other segments of the trajectory to generate infinitely many iso-impulse solutions. We provide clarification for intermediate phasing orbits. In the remainder of the paper and without loss of generality, we assume a maneuver from a low-energy initial orbit to a higher-energy final orbit ($T_0 < T_f$). Upon introducing $n_p$ number of intermediate phasing orbits at the impulse AP on the initial orbit, a constraint on the available time of flight can be written as, 
\begin{equation}
\sum_{k=1}^{n_p} N_k T_k(\alpha_k) + N_0T_0+N_fT_f+N_{\mathrm{pf}} T_{\mathrm{pf}}=TOF,
\label{eq:equality_cond}
\end{equation}
where the value for $n_p$ is not known in advance, but it has an upper bound that can be determined analytically. Orbital periods are denoted as $T_s$ for $s\in\{k,0,f,\text{pf}\}$. $N_k$ is the integer number of revolutions on the $k$-th intermediate phasing orbit, $T_k(\alpha_k)$ is the orbital period of the $k$-th intermediate phasing orbit. Note that $T_k$ is a nonlinear function of $\alpha_k$. Subscripts ``$0$'' and ``$f$'' correspond to the initial and final orbits. In essence, Eq.~\eqref{eq:equality_cond} states that a spacecraft has the option of spending the available time of flight (which can be interpreted as extra time) on, at most, four main segments of an impulsive solution. The spacecraft can make 1) $N_0$ revolutions on the initial orbit, 2) $N_f$ revolutions on the final orbit, 3) $N_\text{pf}$ revolutions on the orbit corresponding to the phase-free connecting arc. Note that the connecting arc, as plotted in Fig.~\ref{fig:twoimp_base}, is part of an ellipse, and the spacecraft can make multiple revolutions to make up for the extra time, and 4) $N_k$ revolutions on $n_p$ number of distinct intermediate phasing orbits corresponding to the sum term. 
\subsection{Definition of solution families and solution envelopes}
Impulsive solutions can be classified based on the number of revolutions that spacecraft make on each segment of the trajectory. Each feasible combination of $N_s$ for $s\in\{k,0,f,\text{pf}\}$ is referred as a solution family. If the combination of all $N_s$ values is such that $k\geq 2$, i.e., at least two intermediate phasing orbits are possible, then, the orbital periods of those phasing orbits can take an infinite number of values \cite{saloglu_existence_2023}. The number of feasible families of impulsive solutions depends on the value of $TOF$, which depends on the base solution parameters and the mission time for time-fixed rendezvous maneuvers. The feasibility of a solution family is defined as the existence of infinitely many solutions. We prove that the lower bound on $TOF$ can be used as a time-feasibility criterion in Section~\ref{sec:feasTime}. Suppose more than two phasing orbits are possible, i.e., $k\geq 3$. In that case, the orbital periods of the infinitely many solutions can be characterized and visualized as polygons, due to Eq.~\eqref{eq:equality_cond} being a linear function of the orbital period values. They are defined as solution envelopes \cite{saloglu_existence_2023}, see Sec.~\ref{sec:sol_env}. 

\subsection{Formulation of impulse-allocation problem at one anchor position}
Let $\alpha_k \in (0,1)$ denote the $\Delta v$ ratio parameter of the $k$-th phasing orbit. The direction of $\Delta \bm{v}$ vectors at an impulse AP is determined by the PVT and we do not change the direction of the impulse; but, $\Delta \bm{v}$ magnitude is modulated upon introducing phasing orbits. Therefore, $\alpha_k$ determines the ratio between the total impulse magnitude at an AP and the $k$-th impulse magnitude. For the $k$-the impulse, we can write $\bm{v}_{p,k} = \bm{v}_0 + \alpha_k \Delta \bm{v}$.

\begin{figure}[htbp!]
    \centering
    \vspace{-4mm}
\includegraphics[width=0.3\linewidth]{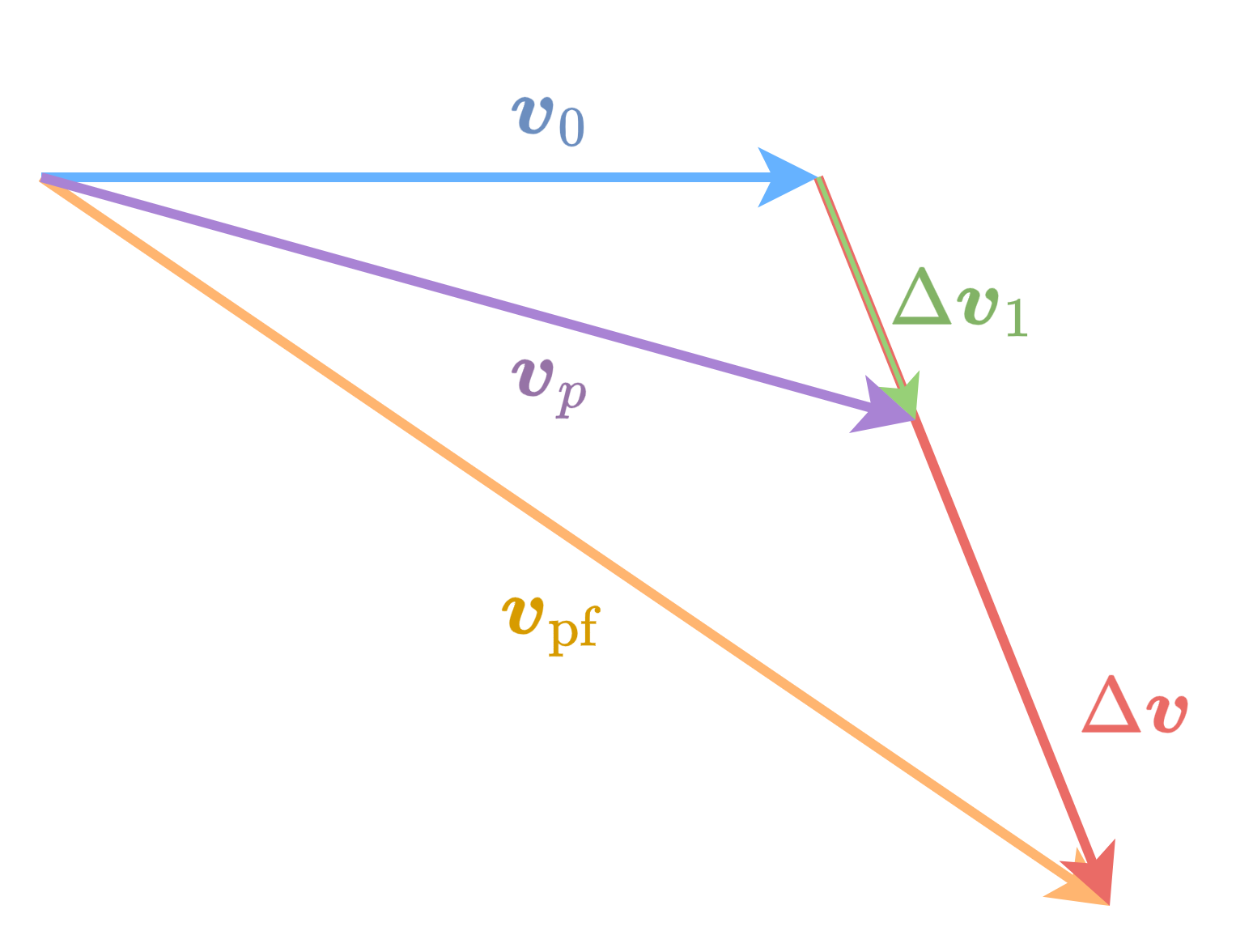}
    \caption{Schematics for the addition of one phasing orbit at the first AP on the Earth orbit in Fig.~\ref{fig:twoimp_base}.}
    \label{fig:alpha_param}
    \vspace{-4mm}
\end{figure}

In Fig.~\ref{fig:alpha_param} schematics of a $\Delta v$-allocation step is drawn in which $\bm{v}_0$ is the initial orbit velocity vector, $\bm{v}_p$ is the phasing orbit velocity vector, and $\bm v_\text{pf}$ is the phase-free arc velocity vector. The $
\Delta \bm{v}_1$ is the required impulse to inject the spacecraft into a phasing orbit and it is shown with a green vector. It is aligned with the total $\Delta \bm{v}$, which is denoted with the red vector. Then, $\Delta \bm{v}_1 = \alpha_1 \Delta \bm{v}$ with $\alpha_1 \in (0,1)$. Note that $\alpha_1 = 0$ and $\alpha_1 = 1$ correspond to the initial orbit and the phase-free arc and are excluded. 
Assuming there are $n_p$ phasing orbits, a set of inequality constraints on the value of $\alpha$ parameters and the orbital period can be formed as, 
\begin{equation}
0<\alpha_1<\cdots<\alpha_{n_p}<1, \quad T_k(\alpha_k)>0, \quad \text { for } \quad  k=1, \ldots, n_p.
\end{equation}

To generate infinitely many impulsive solutions, one can discretize $\alpha_1,\alpha_2,\cdots,\alpha_{n_p-1}$ variables in the range $(0,1)$. Once $\alpha$ values are chosen, the corresponding $T_1, T_2, \cdots T_{n_p-1}$ values can be obtained by using the energy equation as,
\begin{align}
\frac{v_p^2}{2}-\frac{\mu}{r} & = -\frac{\mu}{2 a},&  \xrightarrow[]{T_p = 2 \pi \mu^{-1/2} a^{3/2}} & & T_p&=2 \pi \mu \sqrt{-\frac{1}{8\left(\left(v_p^2 / 2\right)-(\mu / r)\right)^3}},
\label{eq:TOF_alpha}
\end{align}
where $v_p^2 = ||\bm{v}_0 + \alpha_k\Delta \bm{v}||$ with $\bm{v}_0$ denoting the initial orbit velocity (at the AP), $\Delta \bm{v}$ is the total impulse at the impulse AP, $\mu$ is the gravitational parameter and $r$ denotes the magnitude of the position vector at the impulse AP. Then, using the equality condition given in Eq.~\eqref{eq:equality_cond}, $T_{n_p}$ is determined. If $T_{n_p}$ is a feasible value, i.e., $T_0 < T_{n_p} < T_\text{pf}$, then $\alpha_{n_p}$ can be determined by solving the quadratic equation given in Eq.~\eqref{eq:alpha_T} as, 
\begin{align}
v_p^2 = ||\bm{v}_0 + \alpha_k\Delta \bm{v}||,& & \xrightarrow[]{\text{using Eq.}~\eqref{eq:TOF_alpha}} & &\boldsymbol{v}_0^{\top} \boldsymbol{v}_0+2 \alpha_{n_p} \boldsymbol{v}_0^{\top} \Delta \boldsymbol{v}+\alpha_{n_p}^2 \Delta \boldsymbol{v}^{\top} \Delta \boldsymbol{v}+2\left(\frac{\mu^2 \pi^2}{2 T_{n_p}^2}\right)^{1 / 3}-\frac{2 \mu}{r}=0,
\label{eq:alpha_T}
\end{align}
and by choosing the feasible solution of the two roots. Thus, the $\Delta v$-allocation step is a completely analytic process \cite{saloglu_existence_2023}. 

\section{Generation of Multi-Impulse Solutions Using Three-Impulse Base Solutions}\label{sec:threeimpbase}
\begin{figure}[htbp!]
    \centering
    \vspace{-4mm}
    \includegraphics[width=0.7\columnwidth]{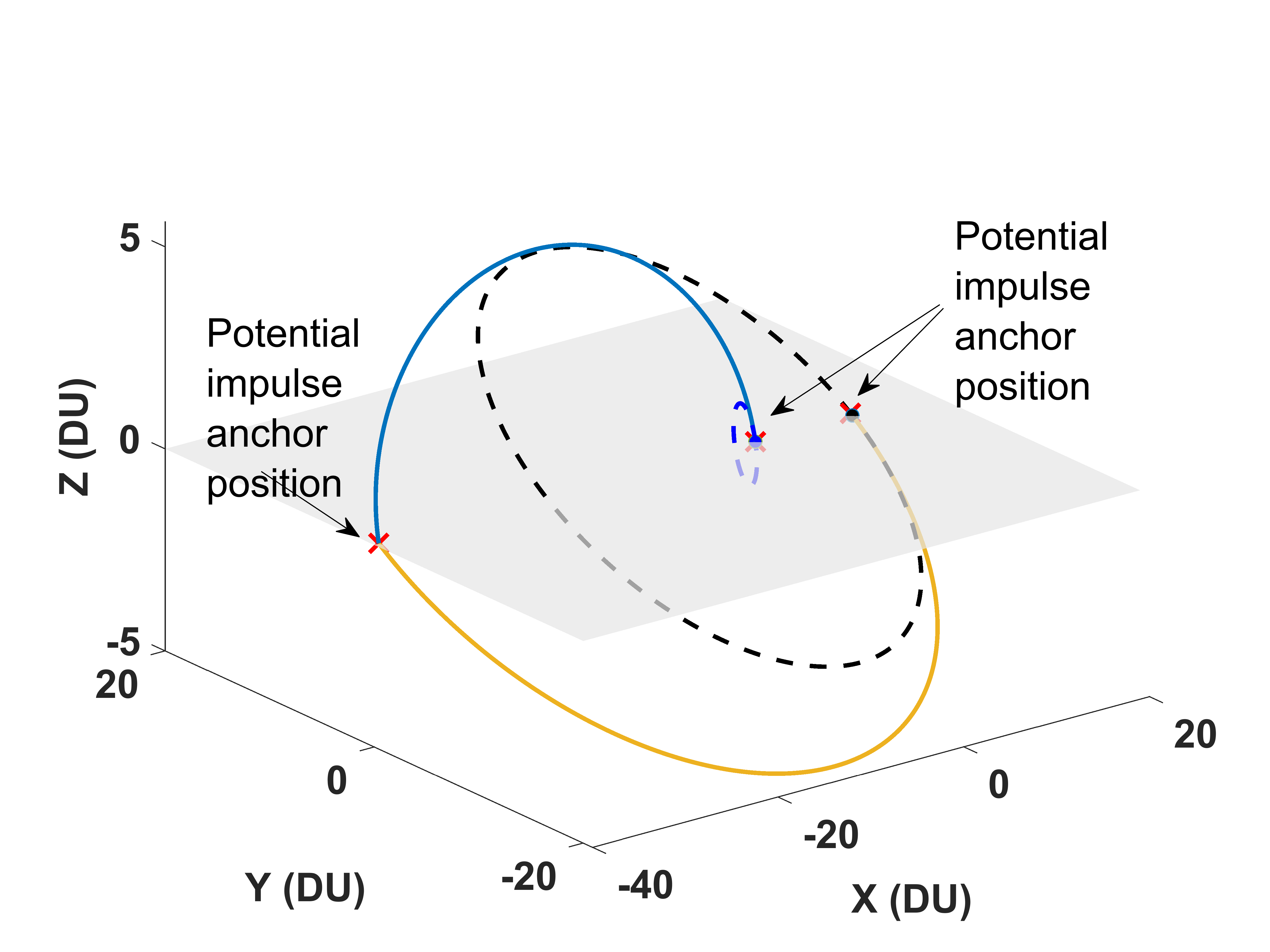}
    \caption{A representative geocentric phase-free base solution with three potential impulse APs ($t_{c_1} = t_{c_2} = 0$).}
    \label{fig:threeimp_base}
   \vspace{-4mm}
\end{figure}

Under the assumption that there are three impulses, the base solution optimization problem can be stated as,
\begin{align}   
\underset{\theta_i, \theta_f, t_{\text{pf},1}, t_{\text{pf},2}, r_\text{mid}, \phi, \lambda}{\text{minimize}}\quad \Delta v_\text{total} = \sum_{i=1}^{3}||\Delta \bm{v}_i|| , 
    \label{eq:threeimpopt}
\end{align}
where similar to the two-impulse phase-free problem, the true anomalies on the initial and final orbits are design variables ($\theta_i$ and $\theta_f$). Here, $t_{\text{pf},1}$ and $t_{\text{pf},2}$ denote the times associated with the two phase-free coast arcs between two consecutive impulses. The midcourse impulse position vector magnitude, ($r_\text{mid}$), azimuth ($\phi$), and elevation ($\lambda$) angles are decision variables to determine the middle impulse position vector. 
Once a three-impulse base solution is obtained, the locations of the three impulse APs are known. The generation of multiple-impulse extremal solutions follows the same steps summarized in Section~\ref{sec:genmultimp} for a single impulse AP. However, any potential impulse APs can be considered for introducing additional intermediate phasing orbits. An example three-impulse base solution is shown in Fig.~\ref{fig:threeimp_base} with three potential impulse APs for a geocentric case (DU stands for distance unit and is equal to the equatorial radius of the Earth). The blue and yellow arcs are the two connecting arcs between three impulses. Typically, one can divide the impulse at an AP with the largest value of impulse. This is to comply with thruster constraints at each single impulse. If more than one impulse AP exists with large magnitude impulses, all potential APs can be considered as impulse APs for solving the $\Delta v$ allocation problem. 
For a three-impulse base solution, we can express the available time of flight as,
\begin{equation} \label{eq:Tava3imp0}
TOF = \underbrace{t_f - t_0}_{\text{mission time}} - \underbrace{(t_{c_1} + t_\text{pf,1} + t_\text{pf,2} + t_{c_2})}_{\text{time spent on base coast arcs}}.
\end{equation}
\subsection{Selection Between Two- or Three-Impulse Base Solutions}\label{sec:basesol}
Given that two- or three-impulse base solutions are possible, a procedure has to be developed for selecting two- or three-impulse base solutions. One of the main distinguishing criteria is $\Delta v_\text{total}$. Therefore, the base solution that requires a lower $\Delta v_\text{total}$ value is the best option with respect to Edelbaum's question \cite{edelbaum1967many}. While computationally more demanding compared to a two-impulse solution, one can first solve for the three-impulse base solution and if the solver tends to make one of the impulses zero, it is typically an indication that the two-impulse base solution has a lower $\Delta v_\text{total}$ value. Maneuvers similar to bi-elliptic ones are more likely to have a lower $\Delta v_\text{total}$ value with a three-impulse base solution than the two-impulse base solution. Thus, transfer maneuvers between highly inclined orbits or those with large changes in the energy are, in general, candidate scenarios for three-impulse base solutions.

\begin{figure}[htbp!]
    \centering
    \includegraphics[width=0.85\textwidth]{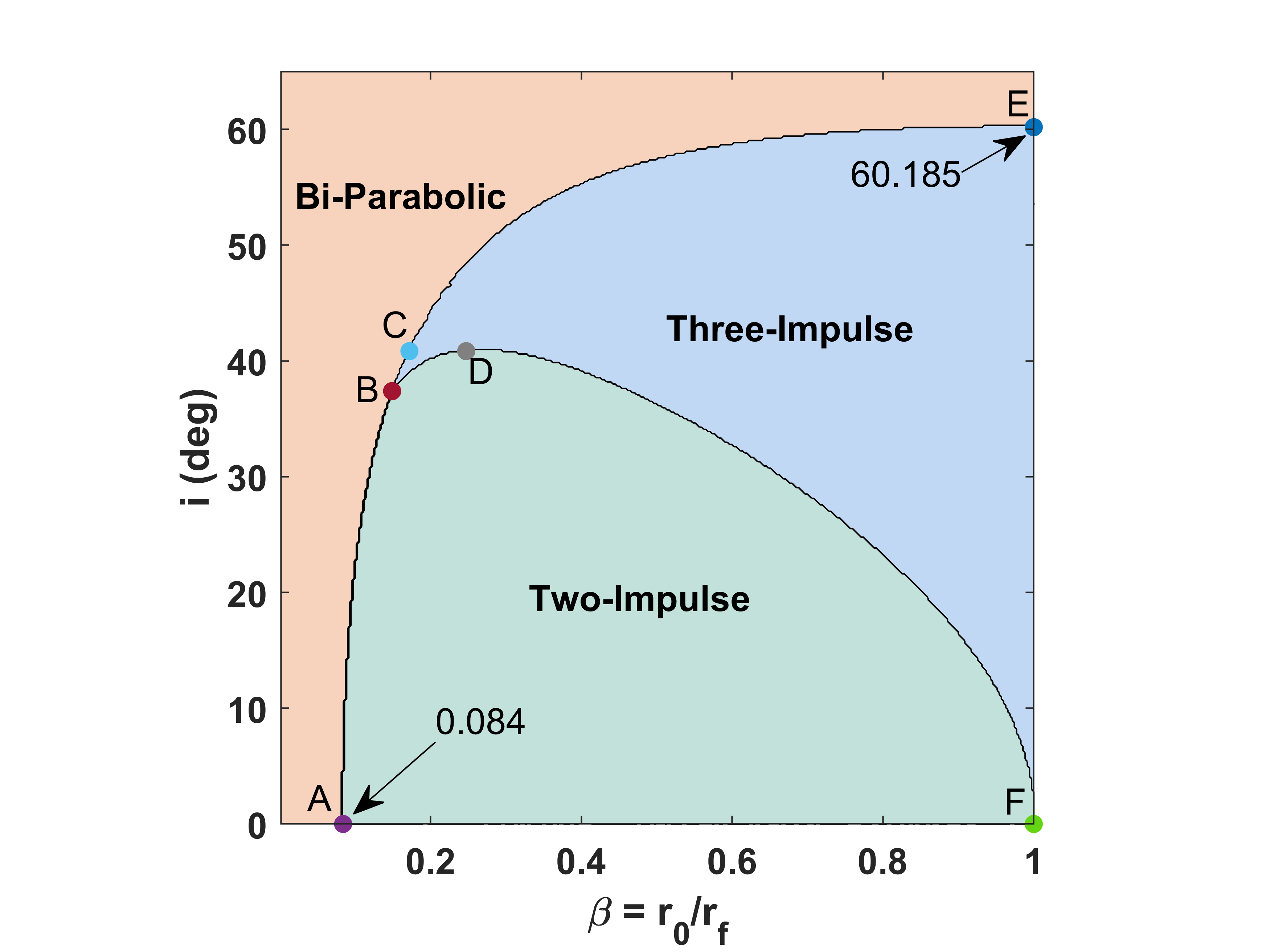}
    \caption{Circle-to-circle maneuvers:  $\Delta v$-optimality for different types of maneuvers and number of impulses.}
    \label{fig:two_three_imp}
\end{figure}

It is theoretically proven that two-impulse (i.e., Hohmann) and three-impulse (i.e., bi-elliptic) transfer trajectories can have the same total $\Delta v$ value (see the separatrix in Fig 6.8 in \cite{curtis2013orbital}) for coplanar, coaxial circular orbits. Let $r_0$ and $r_f$ denote the radii of the initial and final orbits, respectively, and let $\beta = r_0/r_f$. The separatrix is when $11.94 \leq r_f/r_0 \leq  15.58$. For time-free transfer maneuvers between three-dimensional, coaxial, inclined, circular orbits, a categorization on the number of impulses is given by Marec depending on the orbit radii ratio and the inclination change~\cite{marec_1979}. 
We reproduced that plot in Fig.~\ref{fig:two_three_imp}, but we included additional details and data points that facilitate the partitioning of the $(i,\beta)$-space. These points are A($i = 0^{\circ}, \beta = 0.084$), B($i = 37.389^{\circ}, \beta = 0.149$), C($i = 40.841^{\circ}, \beta = 0.1719$), D($i = 40.841^{\circ}, \beta = 0.2473$), E($i = 60.185^{\circ}, \beta = 1.0$), and F($i = 0^{\circ}, \beta = 1.0$). Point B is a bifurcation at which all $\Delta v_\text{total}$ values are equal. Point D is the maximum inclination at which the two-impulse and three-impulse solutions have an equal $\Delta v_\text{total}$ value. At point C, three-impulse and bi-parabolic $\Delta v_\text{total}$ values are equal at the maximum inclination value of point D. 
Each region denotes the required number of impulses for $\Delta v_\text{total}$ optimality. The (hypothetical) bi-parabolic transfer is the case when the inclination change is achieved at an infinity radius with no $\Delta v$ cost. Bi-parabolic transfers, however, are still classified under the three-impulse base solutions. 
When the inclination is zero, the two-impulse solution becomes $\Delta v_\text{total}$-optimal for $\beta \geq 0.084$, or $1/\beta = 11.94$, which is the same value given in \cite{curtis2013orbital}. 
For all practical purposes that exclude the bi-parabolic region, the following categorization of the $\Delta v$ optimality based on the number of impulses can be considered:  
\begin{itemize}
    \item Considering points A and B, we can have both two- and three-impulse base solutions. For each value of $0.084\leq \beta \leq 0.149$, there is a unique inclination that defines the boundary line   
    at which two- and three-impulse base solutions have the same $\Delta v_\text{total}$ value. For inclination values below the boundary line, two-impulse base solutions are optimal. 
    \item Considering point C, for $\beta > 0.1719$ and $i> 40.841^\circ$, three-impulse base solutions are optimal.
    
    \item Considering point D, if $i< 40.841^\circ$ and $\beta > 0.2473$, then, depending on the inclination values, it is possible for two- or three-impulse base solutions to be optimal. It is also possible for two- and three-impulse solutions to have the exact same $\Delta v_{\text{total}}$ corresponding to the separatrix from D to F between the two regions. In such cases, the time of flight is a factor that determines the base solution. 
    \item Considering points B and D, if $i< 40.841^\circ$ and $0.149<\beta <0.2473$, it is possible for either two- or three-impulse base solutions to be optimal. It is also possible for both two- and three-impulse solutions to have the exact same $\Delta v_{\text{total}}$ value, corresponding to the separatrix between the two regions. In such cases, the time of flight is a factor that determines the base solution.
\end{itemize}

We highlight that there can be identical solutions in terms of $\Delta v_\text{total}$ for circular-to-circular transfer cases with different impulse locations. 
In such cases, if a fixed-time rendezvous case is considered, the impulse locations of the base solution can be selected such that the terminal coast times are the lowest (see discussion in Sec.~\ref{sec:geocircir}). The selection of base solutions is affected by the different maneuver types. The goal is to consider all types of maneuvers including, time-free transfer, time-free rendezvous, and time-fixed rendezvous. 
We can consider the following cases:
\begin{itemize}
    \item Time-free transfer: If the goal is a phase-free maneuver without any time constraint, the obtained base solution is the minimum-$\Delta v$ transfer and the $\Delta v$-allocation step can be performed. 
If both two- and three-impulse base solutions have the same $\Delta v$ value, which is a possibility, the base solution that has a shorter flight time is chosen.
\item Time-free rendezvous: These maneuvers can be built upon the phase-free solutions by introducing initial or terminal coast arcs as well as a to-be-determined number of phasing orbits. Since time is free, we can extend the mission time to satisfy the correct phasing corresponding to the rendezvous scenario and then introduce intermediate phasing orbits and terminal coast arcs.
\item Time-fixed rendezvous: Regardless of the number of impulses on the base solution, the time of flight associated with the base solution has to be greater than the given constrained mission time. Therefore, the base solution has to be a time-feasible one. If the base solution has the exact locations with the given initial and final points, the fixed mission time can be equal to the time of flight of the base solution. Additionally, if the mission time of flight is greater than the base solution time of flight, the mission time must be greater than the time of flight of the base solution by, at least, one orbital period of the orbit with the least orbital period. This is required so that the correct phasing is maintained by adding phasing orbits, as we stated in our previous study~\cite{saloglu_existence_2023}. 
\end{itemize}


Time-fixed rendezvous maneuvers require a careful treatment of times of flight. When a time-fixed rendezvous-type mission is considered, the departure and arrival points and the total mission time are all known. We explain the constraint for two-impulse base solutions, but the same constraint applies to three-impulse base solutions. To apply the $\Delta v$-allocation step, 
we have a constraint for time-feasibility of minimum-$\Delta v$ maneuvers, which has to be satisfied by the available time of flight, $TOF$, and can be written as, $TOF \geq T_0$.
This constraint indicates that, at least, one phasing orbit (corresponding to the smallest orbital period) should be possible to be added to the base solution. This is to match the total mission time while maintaining the correct phasing between departure and arrival points. If the time-feasibility constraint is not satisfied, we cannot recover the theoretically minimum-$\Delta v$ base solution under a time-fixed rendezvous scenario. Therefore, the $\Delta v$-optimality has to be sacrificed to find a solution with the given fixed time and (rendezvous) boundary conditions. However, if it is possible to relax the problem to a time-free rendezvous problem, the mission time can be increased, for instance, with multiples of the orbital period of the target orbit, $T_f$, as
 \begin{equation} \label{eq:addkappaTf}
     TOF = (t_f - t_0) - (t_{c_1} + t_\text{pf} + t_{c_2}) + \kappa T_f, \quad \text{for} \quad \kappa \in \mathbb{N}.
 \end{equation}

In Eq.~\eqref{eq:addkappaTf}, multiples of the initial orbit period, $\kappa T_0$ can also be added, but it requires careful treatment by introducing the correct number of revolutions on the $T_f$ to keep the correct phasing. Increasing the mission time allows for introducing phasing orbits to the solution such that the correct phasing and $\Delta v$-optimality are preserved. Although $\Delta v$ is optimal, the new total mission time is increased. Therefore, in such scenarios, maintaining $\Delta v$-optimality requires a trade-off between total mission time and $\Delta v$. If the sacrifice in the $\Delta v_\text{total}$ value is significant, mission time can be extended so that the minimum $\Delta v_\text{total}$ value is preserved. In addition, among the base solutions, one of them can have a feasible solution (with respect to satisfying the boundary conditions and the mission time of flight), but has a worse $\Delta v_\text{total}$. In that case, the sacrifice in the $\Delta v_\text{total}$ value can be made such that the base solution is feasible in terms of time. If the mission time is less than the time of flight of both base solutions, then, the minimum-$\Delta v$ cannot be recovered at all. Thus, our proposed method offers a mechanism for determining if theoretical minimum-$\Delta v$ is possible based on the mission type and time of flight. In cases that the theoretical minimum-$\Delta v$ solution can not be obtained, a Lambert problem \cite{gooding_procedure_1990} can be solved for the given mission time or any other impulsive trajectory optimization can be performed with the knowledge that a sacrifice has been made on optimality with respect to the minimum $\Delta v_\text{total}$ value. Thus, the time-feasibility constraint provides an invaluable certificate for $\Delta v$ optimality.

\begin{figure}[!htbp]
\begin{subfigure}{1\textwidth}
  \centering
  \vspace{-8mm}
  \includegraphics[width=1\columnwidth]{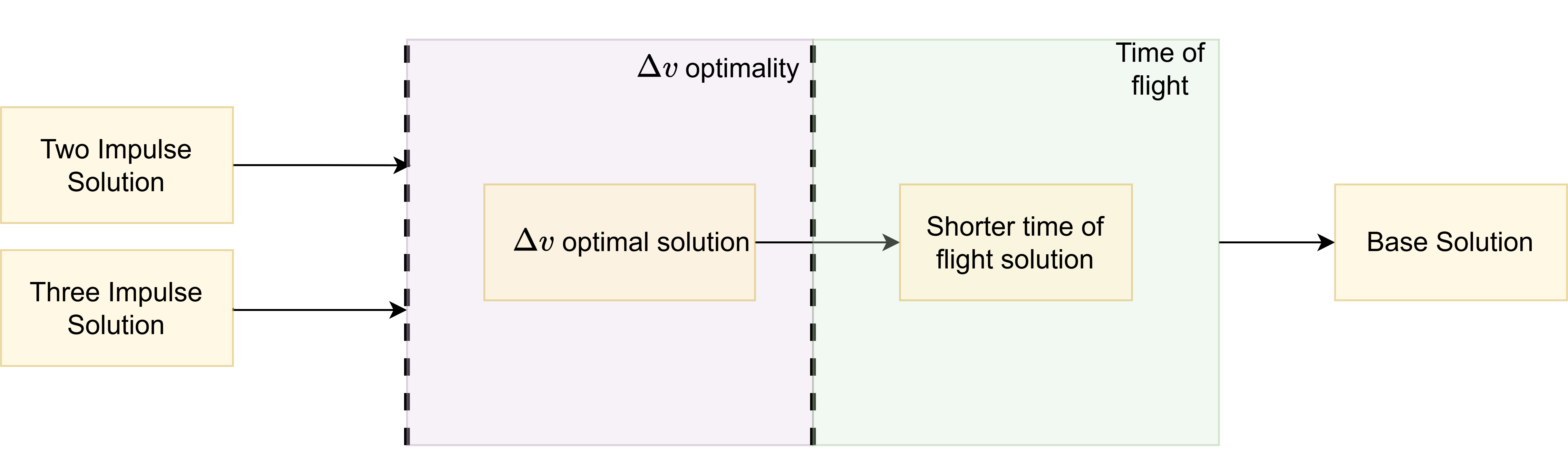}
  \caption{Time-free transfer and time-free rendezvous}
  \label{fig:freetime_base}
\end{subfigure}%

\begin{subfigure}{1\textwidth}
  \centering
  \includegraphics[width=1\columnwidth]{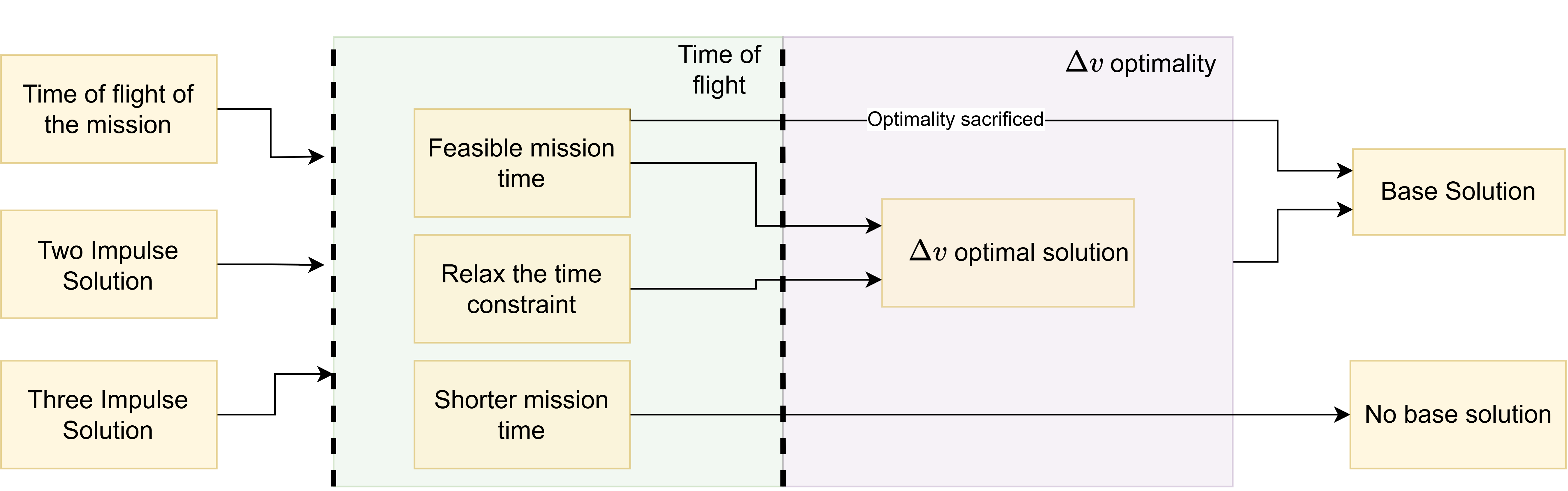}
  \caption{Time-fixed rendezvous}
  \label{fig:fixedtime_base}
\end{subfigure}
\caption{Flowchart of selection of base solutions depending on the maneuver types.}
\label{fig:flowchart_base}
\end{figure}

In principle, we seek solutions to a nonlinear dynamical system corresponding to the inverse-square gravity model. Potential solutions to nonlinear dynamical systems may (and frequently do) exhibit bifurcations with respect to time and other parameters of the problem (e.g., boundary conditions or orbital parameters). It is crucial to take both of these components into consideration when our goal is to classify and analyze potentially all solutions to complex nonlinear dynamical systems. In summary, there are two main cases for the selection of a base solution. The first case is for the time-free maneuvers and the second is for the time-fixed rendezvous. The latter represents the most constrained case. We start with solving the optimization problems for base solutions (per Eqs.~\eqref{eq:biimpopt} and \eqref{eq:threeimpopt}). 
 Then, the first case is the selection of a base solution for the time-free maneuvers. 
 The optimality among the base solutions for the circular-to-circular transfers can be obtained in advance using Fig.~\ref{fig:two_three_imp} with orbit ratio and inclination as the input parameters. For the most general case of inclined elliptical orbits, both problems have to be solved. These optimization problems in Eqs.~\eqref{eq:biimpopt} and \eqref{eq:threeimpopt} can be initialized with the line of nodes as the potential initial location of the impulses. It is known that the best solutions do not correspond to nodal points, but nodal points serve as good initial points \cite{vinh1988optimal}. Once the two- and three-impulse solutions are obtained, the minimum-$\Delta v$ solution is determined for the time-free case. Since there might be solutions with the same $\Delta v_\text{total}$, the one that requires less time is the base solution and the $\Delta v$-allocation problem can be solved. Two layers of sorting can be performed on the two- and thee-impulse solutions, which are shown in Fig.~\ref{fig:freetime_base}. 
 
 The first layer of sorting is performed in terms of $\Delta v_\text{total}$, which is depicted as a pink-colored block in Fig.~\ref{fig:freetime_base}. If the $\Delta v_\text{total}$ values are equal, the next layer of ``filtering" is applied in terms of the required time of flight among the two- and three-impulse solutions, where the solution with a shorter time of flight is selected (green block in Fig.~\ref{fig:freetime_base}). Here, filtering is referred to as a sorting/selection mechanism. Finally, the output is the selected base solution. 
 
 For the second case, the time-fixed rendezvous maneuvers are considered. Another constraint, in this case, is the mission time. Therefore, we have to consider this constraint as one of the inputs for the selection process, as shown in Fig.~\ref{fig:fixedtime_base}. For the second case, two layers of filtering are considered, in terms of the time of flight and the $\Delta v$-optimality criterion, but in reverse order compared to the first case. The flowchart is shown in Fig.~\ref{fig:fixedtime_base} with an additional input for the mission time. Once the two- and three-impulse phase-free base solutions are obtained, the total times of these solutions are analyzed. The solutions with the feasible time of flight (i.e., $TOF\geq T_0$) pass to the second layer. At the second layer, $\Delta v$-optimality is considered. Then, the solution that requires the least amount of $\Delta v$ is selected as the base one. In the case when the time of flight constraint of the mission can be relaxed, the minimum-$\Delta v$ solution is the base one, even if its time of flight is not originally feasible. This is denoted by the ``relax the time constraint'' option in Fig.~\ref{fig:fixedtime_base}. Finally, if the mission time is shorter than both of the phase-free solutions, that means the minimum-$\Delta v$ cannot be recovered and none of the base solutions can be used, which is the ``shorter mission time'' case in Fig.~\ref{fig:fixedtime_base}. In that case, a Lambert problem can be solved, which increases the $\Delta v$ value, but the solution satisfies the mission time constraint. Also, the obtained phase-free solution can be used as an initial guess or even the first step of a homotopic approach for the time-fixed impulsive optimization process.

\section{Generation of Multi-Impulse Solutions with Three Impulse Anchor Positions}\label{sec:multianchor}
We present the steps for generating infinitely many iso-impulse solutions (corresponding to the minimum-$\Delta v$ base solutions) if two or three impulse APs are considered simultaneously for time-fixed rendezvous maneuvers. In addition, we explain the feasibility of having infinitely many solutions for a candidate solution family. 
Finally, the generation of solution envelopes is discussed. 

\begin{figure}[htbp!]
    \centering
    \includegraphics[width=1.0\columnwidth]{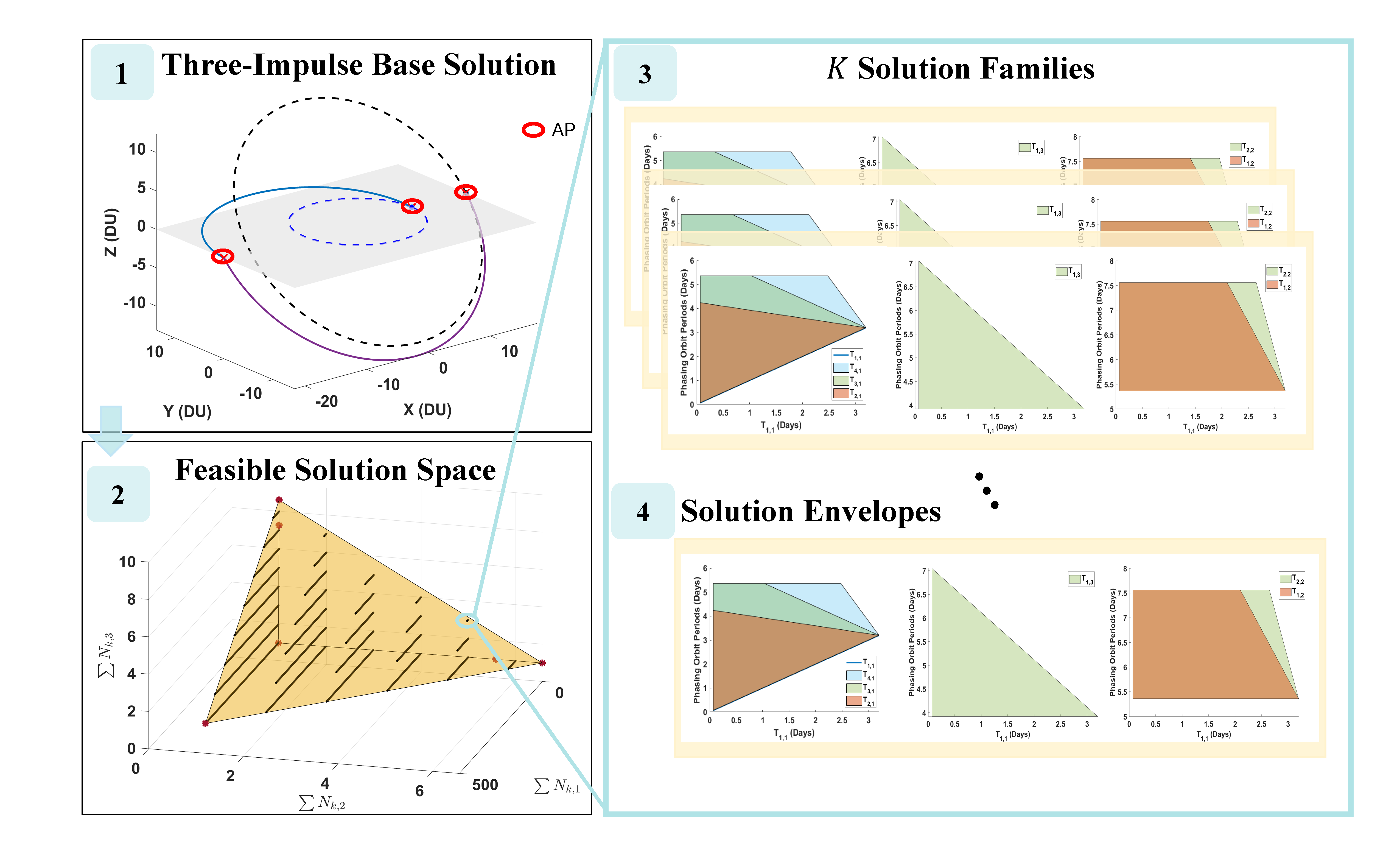}
    \caption{Summary of connection between the base solution, feasible solution space, and the solution envelopes.}
    \label{fig:vis_abs}
\end{figure}

At a high level, we can analyze and categorize all possible iso-impulse solutions in terms of four layers of a superstructure illustrated in Fig.~\ref{fig:vis_abs} for a three-impulse base solution and labeled as 1) the base solution, 2) feasible solution spaces, 3) solution families, and 4) solution envelopes. In the first layer, a two- or three-impulse base solution is selected according to the procedure outlined in Sec.~\ref{sec:basesol} and based on $\Delta v$-optimality or time-feasibility criterion. Intermediate phasing orbits are added to the trajectory at the APs of the chosen base solution. In the second layer, the feasible solution space is determined in terms of the total number of revolutions at each AP. That is, at each AP the lower and upper bounds of the total number of revolutions are determined, which results in a polytope for three-impulse base solutions. The details on how to obtain the feasible solution spaces are explained in Sec.~\ref{sec:feasTime}, but we note that for two-impulse base solutions, the solution space looks like a polygon (see Figs.~\ref{fig:2AP_solspace_E2D} and \ref{fig:2AP_solspace}). 

In the third layer, each feasible solution creates $K$ distinct solution families, where $K$ is an integer value defined as the total number of solution families belonging to a feasible solution. For instance, if we can add 30 revolutions on intermediate phasing orbits at the first AP, we can have the 30 revolutions occur either on a unique phasing orbit or split between multiple distinct phasing orbits. The combination of all possibilities at all APs constitutes the $K$ solution families. The generation of these iso-impulse solution families is discussed in Sec.~\ref{sec:isoimp_3AP}. 
In the fourth layer, each solution family can be depicted in terms of multiple solution envelopes corresponding to each phasing orbit to visualize the ranges of the phasing orbits' orbital period values. For each solution family, there are distinct sets of solution envelopes. In this fourth layer, each solution family consists of infinitely many solutions. The details on the solution envelopes are explained in Sec.~\ref{sec:sol_env}. We emphasize that all solutions have the same $\Delta v_\text{total}$ value, which is equal to the $\Delta v_\text{total}$ value of the chosen base solution.  

\subsection{Generation of Iso-Impulse Solutions}\label{sec:isoimp_3AP}

If multiple impulse APs exist, using Eq.~\eqref{eq:Tava3imp0}, the time constraint has to be written as, 
\begin{equation}
\sum_{i=1}^{n_i}\sum_{k=1}^{n_{p,i}} N_{k,i}T_{k,i}(\alpha_{k,i}) + N_{0}T_{0} + N_f~T_f + N_\text{pf,1}~T_\text{pf,1} + N_\text{pf,2}~T_\text{pf,2} = T O F,
\label{eq:TOF_gen}
\end{equation}
where $n_i$ denotes the number of impulse APs and $n_{p,i}$ is the number of phasing orbits at the $i$-th impulse AP. 
Intermediate phasing orbits are denoted by subscripts ``$k, i$'', and subscripts ``$\text{pf,1}$'' and ``$\text{pf,2}$'' denote  the phase-free arcs. 
Solving the $\Delta v$-allocation problem is similar to the one summarized in Section~\ref{sec:genmultimp}. The inequality constraints on the impulse ratio values are true for each AP, 
\begin{equation}
0<\alpha_{1,i}<\cdots<\alpha_{n_p,i}<1, \quad T_{k,i}>0, \quad \text { for } k=1, \ldots, n_{p,i}, \ \text{and} \ i = 1,2,3.
\label{eq:ineq_new}
\end{equation}

To generate infinitely many iso-impulse solutions, we discretize the $\alpha_{k,1}, \cdots, \alpha_{(n_{p,1}-1),1}$ for the first AP. For the second and third APs, we discretize all the associated $\alpha$ parameters, such as $\alpha_{k,i}, \cdots, \alpha_{n_{p,i},i}$ for $i \in \{2,3\}$. Then, the corresponding orbital periods, $T_{1,1}, \cdots, T_{(n_{p,1}-1),1}$ and $T_{k,i}$ are determined by using Eq.~\eqref{eq:TOF_alpha} for each discretized value of $\alpha$ at each AP. Then, the remaining $T_{n_{p,1},1}$ is determined from Eq.~\eqref{eq:TOF_gen}. The remaining $\alpha_{n_{p,1},1}$ parameters are determined by solving Eq.~\eqref{eq:alpha_T}. Finally, the inequality conditions given in Eq.~\eqref{eq:ineq_new} are used to filter out the feasible solutions generated for each solution family with different combinations of $N$ parameters for each segment of the solution.

\subsection{Time-Feasibility Criterion of Solution Families}\label{sec:feasTime}
The feasibility of having infinitely many iso-impulse solutions for a solution family can be determined analytically. The relation for the feasibility of infinitely many solutions has to be updated when there are multiple impulse APs. For each impulse AP, we first determine the maximum number of intermediate phasing orbits that can be added. Therefore, we define a new available time that can only be distributed among the intermediate phasing orbits,
\begin{equation}
    TOF_p = TOF - N_fT_f - N_\text{pf,1}T_\text{pf,1} - N_\text{pf,2}T_\text{pf,2} - N_0T_0,
\end{equation}
which is the time that can be spent on the intermediate phasing orbits after subtracting the times spent on the initial and final orbits, as well as the phase-free arcs from $TOF$. The feasibility equation with two or three APs is
\begin{equation}
    T_0 \leq \frac{TOF_p - \sum_{i=2}^{3}\sum_{k=1}^{n_{p,i}}N_{k,i}T_{k,i}(\alpha_{k,i} = 0)}{\sum_{k=1}^{n_{p,1}}N_{k,1}} \leq T_{k,1}(\alpha_{k,1} = 1),
    \label{eq:feas2}
\end{equation}
where $\sum_{k=1}^{n_{p,i}}N_{k,i}$ denotes the total number of revolutions on the phasing orbits at the $i$-th AP. $T_{k,i}(\alpha_{k,i} = 0)$ is the value of the orbital period when no impulse is applied at the $i$-th AP. It should be noted that Eq.~\eqref{eq:feas2} still valid either $\sum_{k=1}^{n_{p,2}}N_{k,2} = 0$ or $\sum_{k=1}^{n_{p,3}}N_{k,3} = 0$. However, we assume at least one phasing orbit at each AP to consider the most general case. To determine the upper bound of the $\sum_{k=1}^{n_{p,i}}N_{k,i}$, we can derive the following relation from Eq.~\eqref{eq:feas2} with the assumption $\sum_{k=1}^{n_{p,1}}N_{k,1} = 1$, using the left-hand side of the inequality,
\begin{equation}
    TOF_p - \underbrace{\sum_{k=1}^{n_{p,2}}N_{k,2}T_{k,2}(\alpha_{k,2} = 0)}_{\text{time spent at the 2nd AP}}  - \underbrace{\sum_{k=1}^{n_{p,3}}N_{k,3}T_{k,3}(\alpha_{k,3} = 0)}_{\text{time spent at the 3rd AP}} \geq T_0.
    \label{eq:feas1}
\end{equation}

Using Eq.~\eqref{eq:feas1}, we can determine the maximum number of phasing orbits at the second or the third AP, denoted as $n_{p,2}$ and $n_{p,3}$, respectively, assuming one revolution on each phasing orbit. Similarly, the upper bound of the total number of revolutions on phasing orbits, $\sum_{k=1}^{n_{p,2}}N_{k,2}$ or $\sum_{k=1}^{n_{p,3}}N_{k,3}$, can be determined. In Eq.~\eqref{eq:feas1}, when $\sum_{k=1}^{n_{p,2}}N_{k,2} = 1$, the upper bound on $\sum_{k=1}^{n_{p,3}}N_{k,3}$ can be found and vice versa. When $\sum_{k=1}^{n_{p,2}}N_{k,2} = \sum_{k=1}^{n_{p,3}}N_{k,3} = 1$, the upper bound of the $\sum_{k=1}^{n_{p,1}}N_{k,1}$ is determined from Eq.~\eqref{eq:feas2}. The upper bounds are helpful when the feasibility of every possible combination of phasing orbits is checked using Eq.~\eqref{eq:feas2}. 

To fully represent the feasible solution space as a polytope (see Sec.\ref{sec:GEO3APs}), the range (upper and lower bounds) of $\sum_{k=1}^{n_{p,i}}N_{k,i}$ have to be determined. The right-hand side of the two-sided inequality in Eq.~\eqref{eq:feas2} is used to calculate the lower bounds with the assumption $\sum_{k=1}^{n_{p,1}}N_{k,1} = 1$, 
\begin{equation}
    TOF_p - \sum_{k=1}^{n_{p,2}}N_{k,2}T_{k,2}(\alpha_{k,2} = 0) - \sum_{k=1}^{n_{p,3}}N_{k,3}T_{k,3}(\alpha_{k,3} = 0) \leq T_{k,1}(\alpha_{k,1} = 1).
    \label{eq:feas3}
\end{equation}

The lower bound of $\sum_{k=1}^{n_{p,1}}N_{k,1}$ is determined from Eq.~\eqref{eq:feas2} when $\sum_{k=1}^{n_{p,2}}N_{k,2} = \sum_{k=1}^{n_{p,3}}N_{k,3} = 1$. Equations~\eqref{eq:feas1} and \eqref{eq:feas3} are both derived from Eq.~\eqref{eq:feas2}. Therefore, Eq.~\eqref{eq:feas2} can be directly used to calculate the ranges of each $\sum_{k=1}^{n_{p,i}}N_{k,i}$.

It should be noted that $T_{k,i}(\alpha_{k,i} = 0)$ has to provide the minimum value of the orbital period that can be added at the impulse AP. If the total impulse at the AP decreases the energy of the orbit, then $T_{k,i}(\alpha_{k,i} = 1)$ can give the minimum orbital period and that value has to be used. 

\subsection{Generation of Solution Envelopes}\label{sec:sol_env}
Solution envelopes introduced by Saloglu et al. \cite{saloglu_existence_2023} represent the complete space of orbital period values of the phasing orbits of all infinity of iso-impulse solutions. 
 We have shown that solution envelopes take polygonal shapes, which means that the vertices of the solution envelopes contain sufficient data for characterizing and visualizing them. Analytic relations for determining the corner points of these solution envelopes have to be modified depending on the number of impulse APs, which extends the work in \cite{saloglu_existence_2023}. This modification is required since $TOP_p$ can be distributed simultaneously between two or three impulse APs. The procedure for introducing phasing orbits at two or three distinct impulse APs is not coupled through the $\Delta v_\text{total}$ requirement. This decoupling requires a special treatment compared to when phasing orbits are added at the same impulse AP since we want to plot the solution envelopes in terms of $T_{1,1}$. 
 
\begin{figure}[!htbp]
\vspace{-3mm}
\begin{subfigure}{.5\textwidth}
  \centering
\includegraphics[width=0.8\columnwidth]{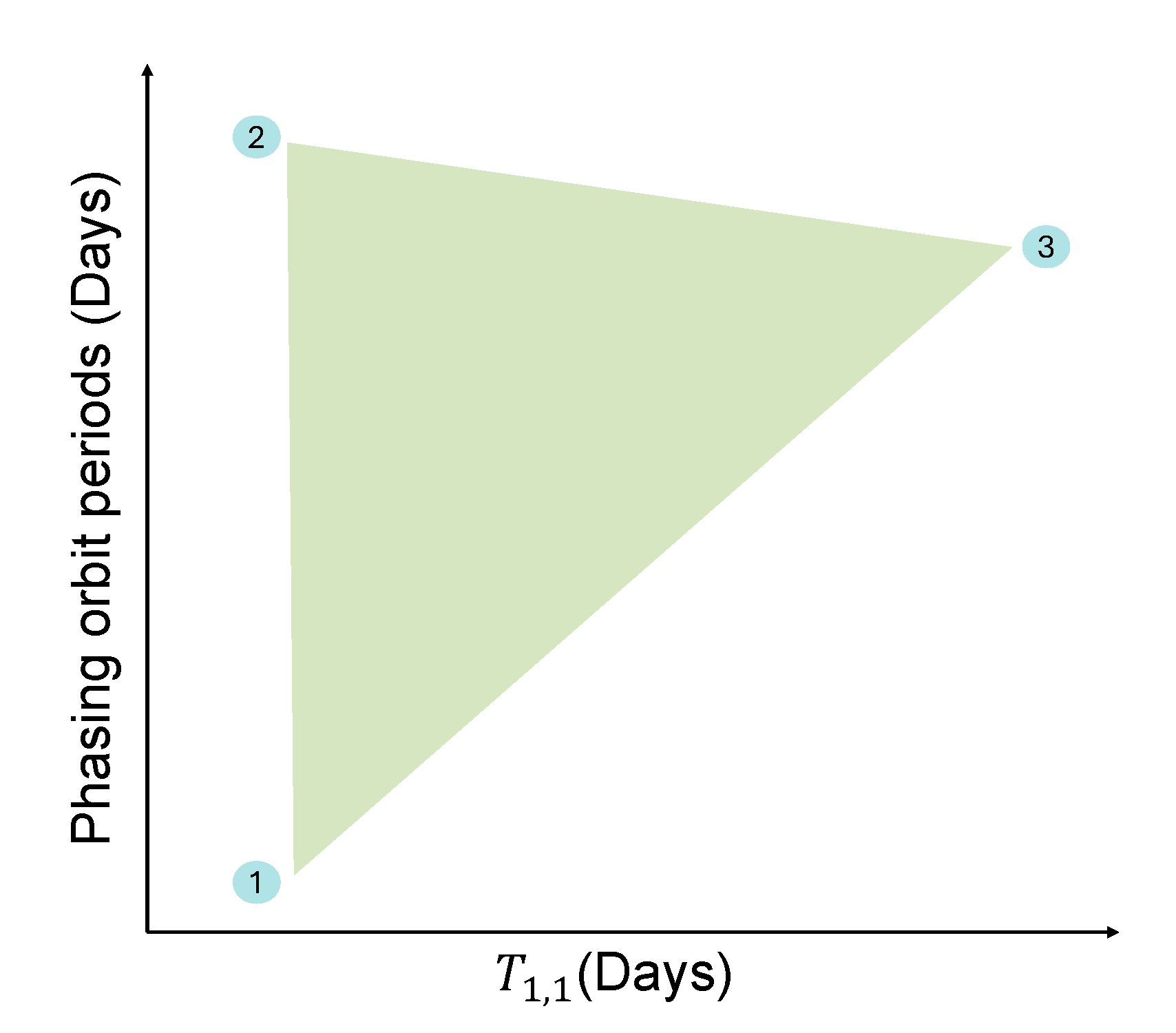}
  \caption{Three main corners}
  \label{fig:solenv1}
\end{subfigure}%
\begin{subfigure}{.5\textwidth}
  \centering
  \includegraphics[width=0.8\columnwidth]{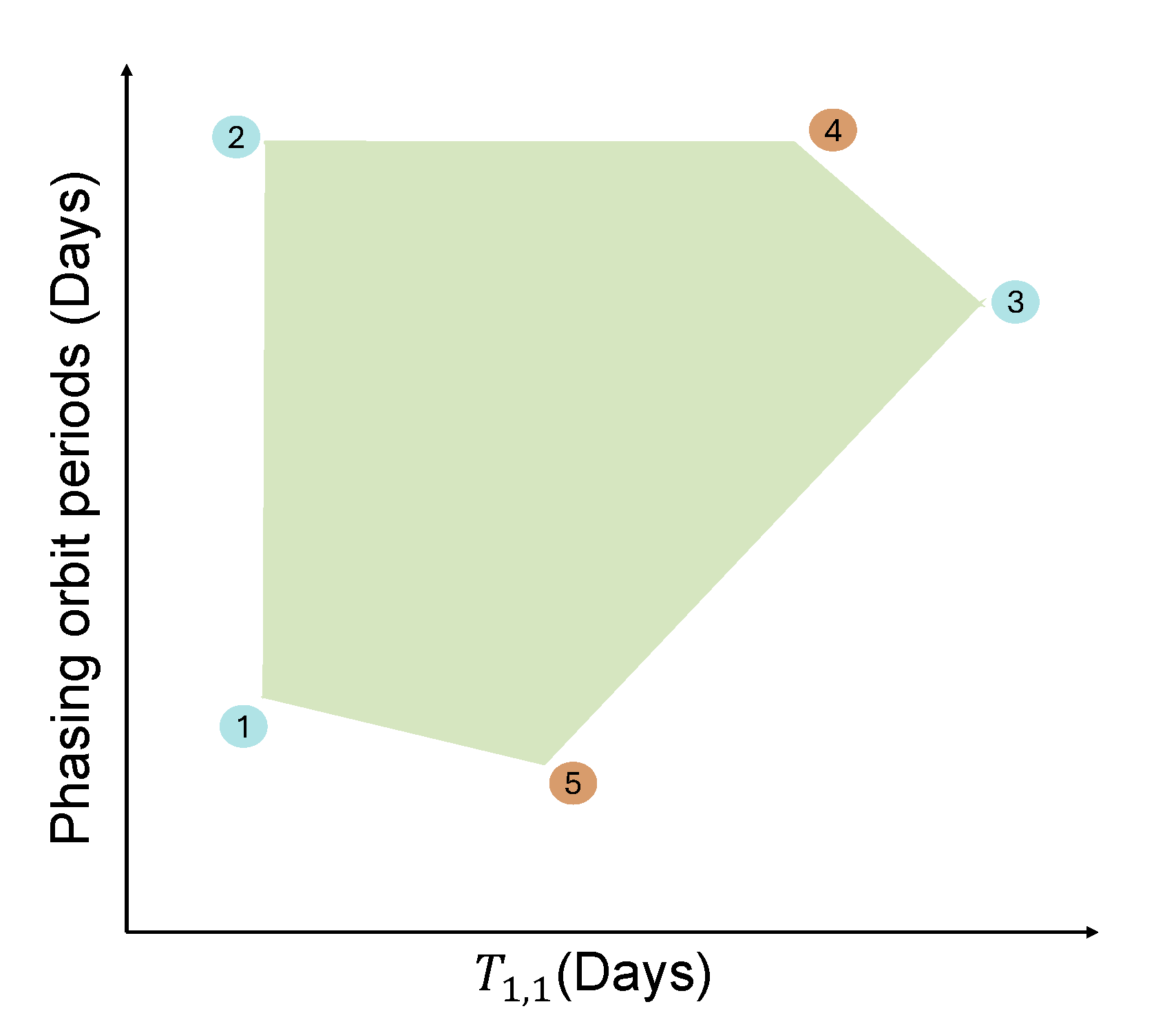}
  \caption{Three main corners and additional corners}
  \label{fig:solenv2}
\end{subfigure}
\caption{Polygonal shapes of representative solution envelopes with numbered corner points}
\vspace{-3mm}
\label{fig:solenv_sketch}
\end{figure}

 Two general polygonal regions are possible, and representative solution envelopes are shown in Fig.~\ref{fig:solenv_sketch}. To precisely determine the corner points of the solution envelopes, we have to answer two questions: 
 \begin{enumerate}
     \item What are the maximum and minimum limits of the orbital period values, when $T_{1,1}=T_0$? (see points ``1'' and ``2'' in Figs.~\ref{fig:solenv1} and \ref{fig:solenv2}), and
     \item What are the values of orbital periods, when $T_{1,1}$ takes its maximum value? (see corner ``3'' in Figs.~\ref{fig:solenv1} and \ref{fig:solenv2}).
 \end{enumerate}

Considering these questions, three main corners must be determined for each phasing orbit solution envelope. These corners are numbered ``1'', ``2'', and ``3'' in Fig.~\ref{fig:solenv1}. However, suppose any of these values are beyond the maximum or minimum boundary that can be taken at the corresponding impulse AP. In those cases, there is a possibility for one additional corner point for each case. These additional corner points are numbered ``4`'' and ``5'' in Fig.~\ref{fig:solenv2}, but a solution envelope may have only one additional corner point (e.g., 4 or 5). Fig.~\ref{fig:solenv2} shows a general case with two additional corner points. When we consider the second and third APs, the relations to determine the corner points are similar. When generating solutions, different AP phasing orbits are not related in terms of $\Delta v_\text{total}$. 
Therefore, for the first AP, we treat the second, and third AP phasing orbit periods as constant parameters and subtract them from the available time for phasing orbits, $TOF_p$. The detailed relations to generate the solution envelopes are given in the Appendix. 

All the corner points can also be determined by solving a set of linear programming (LP) problems by changing the objective of the LP problem for each corner point. For instance, the minimum and maximum values of $T_{2,2}$ when $T_{1,1}=T_0$ are determined by minimizing and maximizing the corresponding element of the design variable vector. The design vector consists of all the phasing orbit orbital period parameters at all the APs, $\boldsymbol x = [T_{1,1}, \cdots, T_{n_{p,1},1}, T_{1,2}, \cdots, T_{n_{p,2},2}, T_{1,3}, \cdots, T_{n_{p,3},3}]$ where their upper and lower bounds are known. The linear equality constraint is the total $TOF$ constraint given in Eq.~\eqref{eq:TOF_gen} and linear inequality constraints are $T_{k,i} \leq T_{k+1,i}$ for each $k = 1,\cdots, n_{p,i}$ where $i \in \{1, 2, 3\}$. We have used an LP-based formulation to validate the analytic relations in the Appendix and to verify the obtained solution envelopes for the considered test problems.





\section{Results}\label{sec:results}
We consider a geocentric transfer problem where a three-impulse base solution generates the minimum-$\Delta v$ solution. First, we show the solutions with a single AP for the geocentric case. In addition, we present solutions with two APs simultaneously for the Earth-to-Dionysus problem and three APs simultaneously for the geocentric problem. We present the solution envelopes for two different families of solutions for the Earth-to-Dionysus and geocentric problems. 

\subsection{Iso-Impulse Solutions Using a Three-Impulse Base Solution}
A geocentric transfer example problem is solved to show the cases where a three-impulse base solution requires a lower $\Delta v_\text{total}$ than a two-impulse base solution. Orbital elements of the initial and target co-axial elliptical orbits are given in Table~\ref{tab:boundary}. 
A large inclination change between the initial and the target orbits is considered. 

\begin{table}[htbp!]
\caption{Geocentric problem: classical orbital elements of the initial and target orbits.}
    \centering
    \begin{tabular}{|c|c|}
    \hline
      Orbit   & $[a~(km), e, i, \omega, \Omega, \theta]$ \\
      \hline
      Initial & $[7000, 0.02, 60^{\circ}, 0^{\circ}, 0^{\circ}, \theta_i]$ \\
      \hline
       Target  & $[105000, 0.3, 12^{\circ}, 0^{\circ}, 0^{\circ}, \theta_f]$ \\
       \hline
    \end{tabular}
    \label{tab:boundary}
\end{table}
\begin{figure}[!htbp]
\begin{subfigure}{.5\textwidth}
  \centering
  \includegraphics[width=1\columnwidth]{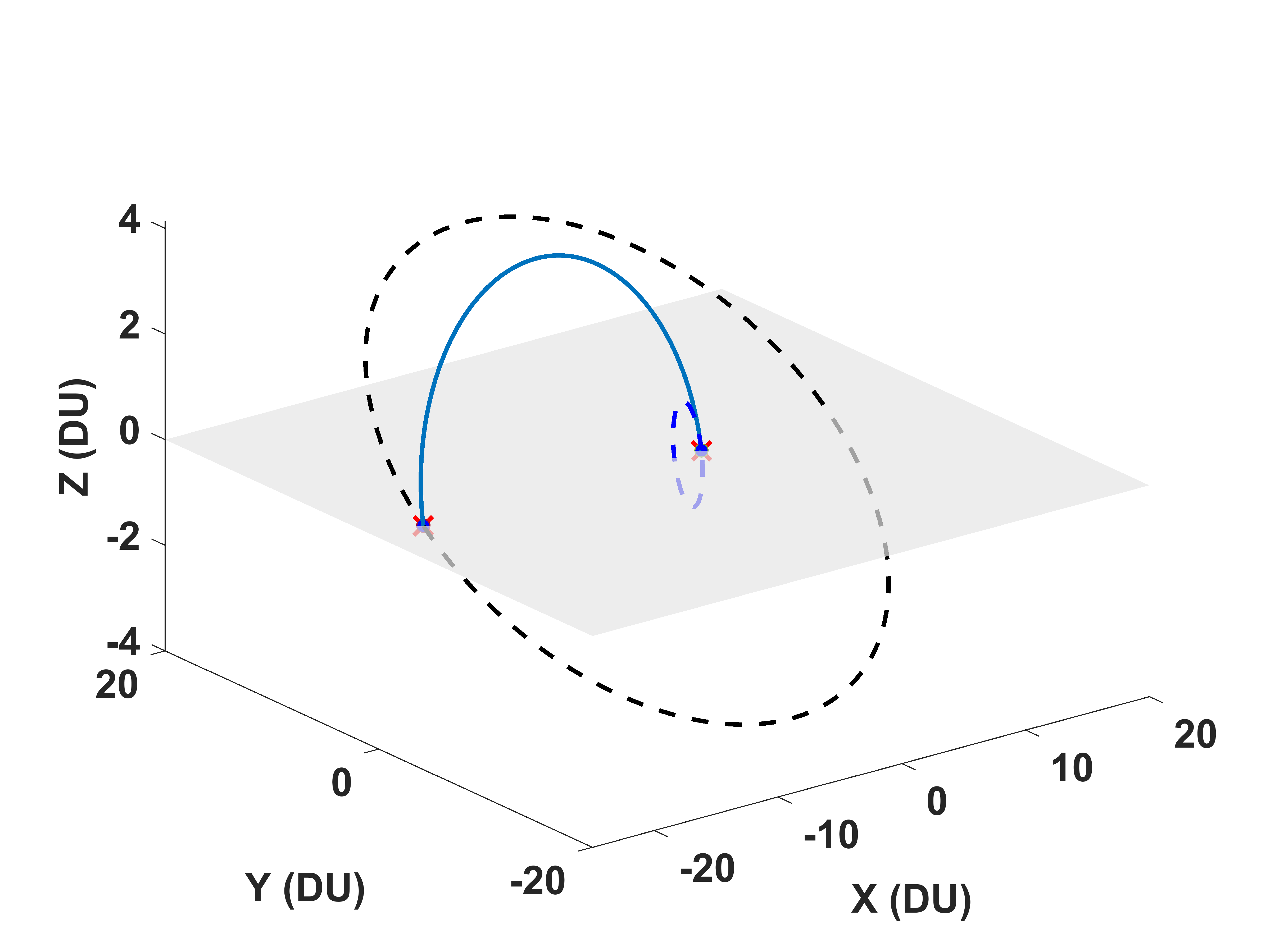}
  \caption{Trajectory (DU = Earth radius)}
  \label{fig:2imptraj}
\end{subfigure}%
\begin{subfigure}{.5\textwidth}
  \centering
  \includegraphics[width=1\columnwidth]{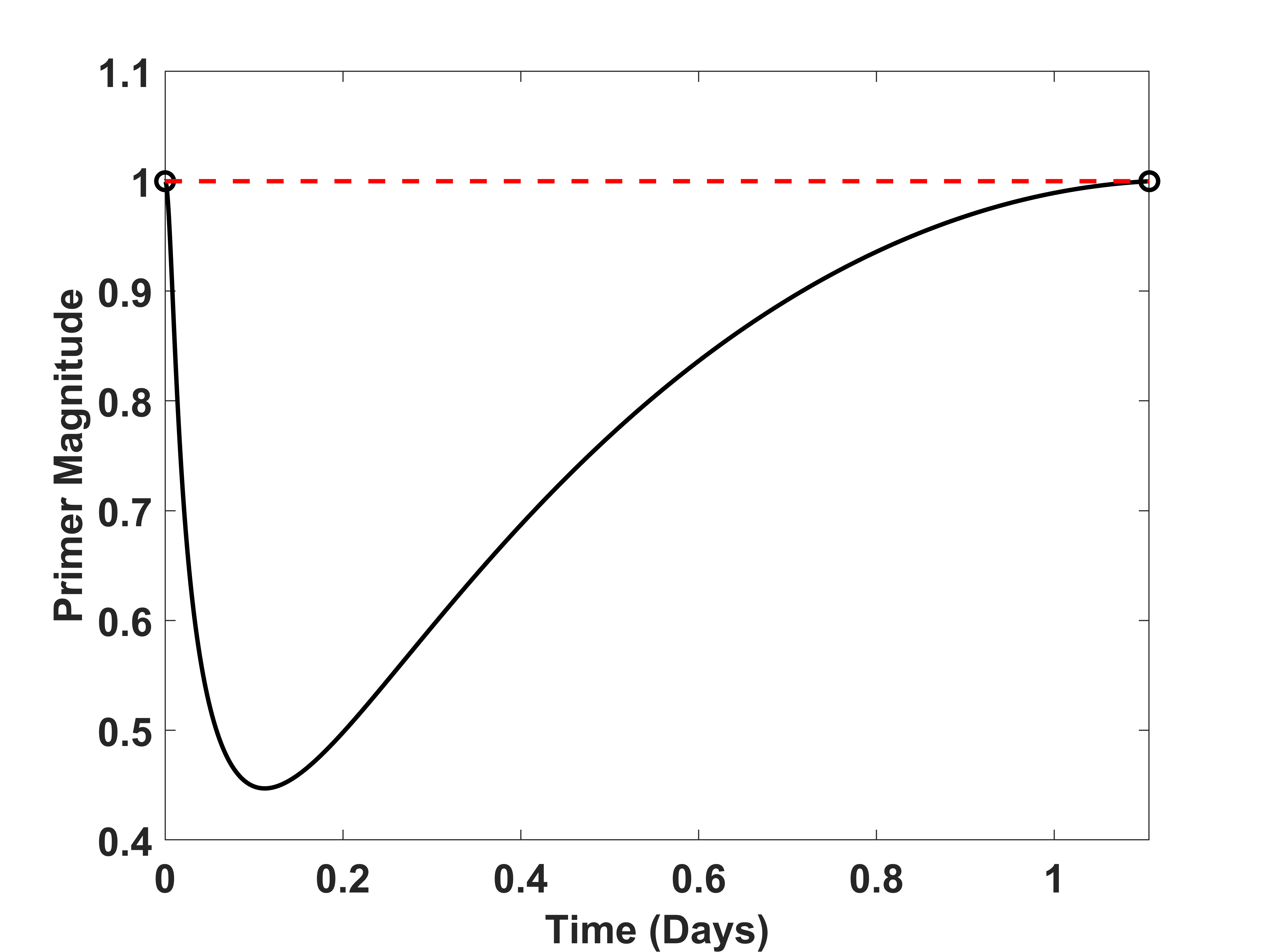}
  \caption{Primer vector magnitude vs. time}
  \label{fig:2imp_prim}
\end{subfigure}
\caption{Geocentric problem: two-impulse base solution ($t_{c_1} = t_{c_2} = 0$). }
\label{fig:twoimp}
\end{figure}
The two-impulse base solution is shown in Fig.~\ref{fig:twoimp}. In Fig.~\ref{fig:2imptraj}, the orbits and the two-impulse base solution are shown. The plane of $z = 0$ is shown to denote the ascending and descending node locations. 
Dashed blue and black orbits represent the initial and target orbits, respectively. Red cross markers denote the impulse locations. Departure and arrival points are shown with blue and black dot markers. The transfer starts from the perigee point of the initial orbit and terminates at the apogee point of the target orbit. Note that the initial and final true anomalies are free and $\Delta v_\text{total}$ is $3.9618011$ km/s with of $\Delta v_1 = 2.8246140$ km/s and $\Delta v_2 =1.1371871$ km/s individual impulses. Most of the inclination change occurs at the second impulse AP, when the spacecraft is at the apogee point of the target orbit. The primer vector magnitude time history is given in Fig.~\ref{fig:2imp_prim}. The primer vector magnitude stays below 1 during the coast arc between the two impulses, characterizing an extremal bi-impulse solution.

\begin{figure}[!htbp]
\begin{subfigure}{.5\textwidth}
  \centering
  \includegraphics[width=1\columnwidth]{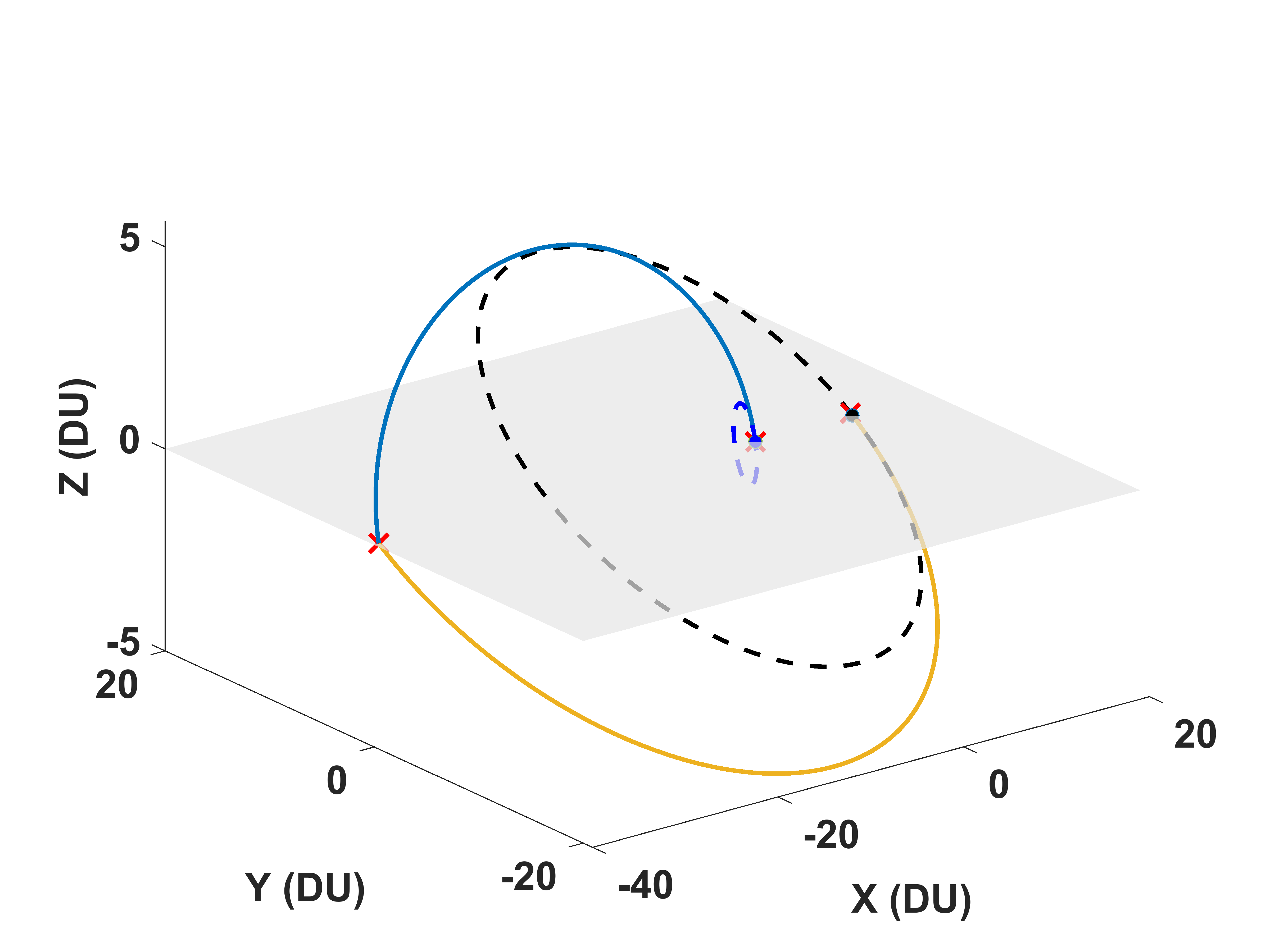}
  \caption{Trajectory}
  \label{fig:3imp_traj}
\end{subfigure}%
\begin{subfigure}{.5\textwidth}
  \centering
  \includegraphics[width=1\columnwidth]{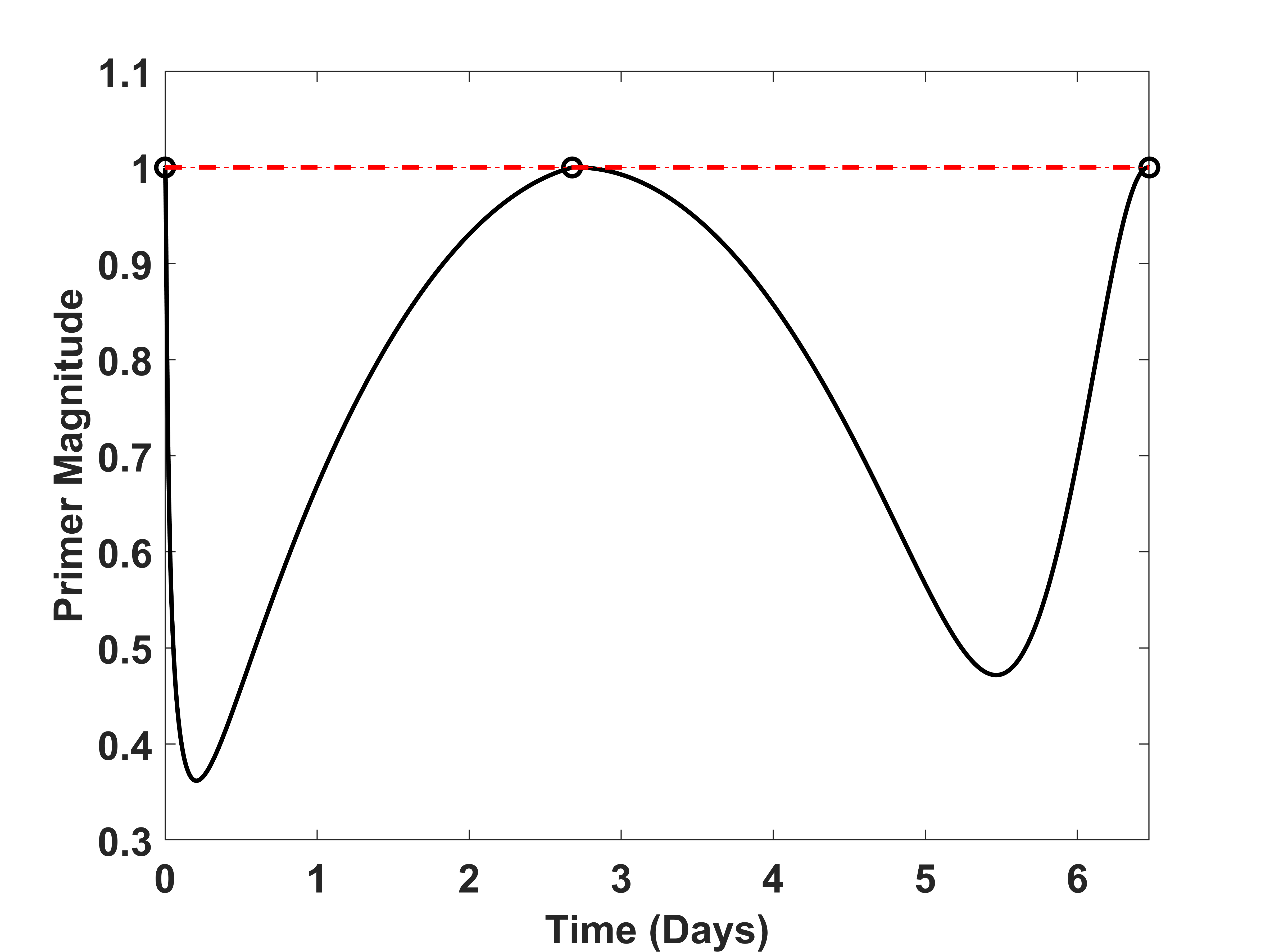}
  \caption{Primer vector magnitude vs. time}
  \label{fig:3imp_prim}
\end{subfigure}
\caption{Geocentric problem: the three-impulse base solution ($t_{c_1} = t_{c_2} = 0$).}
\label{fig:threeimp}
\end{figure}

To compare the two- and three-impulse base solutions, we solve the optimization problem given in Eq.~\eqref{eq:threeimpopt} using MATLAB's \texttt{particleswarm} algorithm. The solution is given in Fig.~\ref{fig:threeimp} with the trajectory depicted in Fig.~\ref{fig:3imp_traj}.
The  $\Delta v_\text{total}$ is $3.8641159$ km/s with $\Delta v_1 = 2.9390$ km/s, $\Delta v_2 = 0.6815$ km/s and $\Delta v_3 = 0.2436$ km/s. The primer vector magnitude time history is shown in Fig.~\ref{fig:3imp_prim}. 
Comparing the $\Delta v_\text{total}$ values of the two- and three-impulse base solutions, the three-impulse solution is optimal. 
We can choose an impulse AP as the next step starting from the three-impulse base solution since it is the minimum-$\Delta v$ solution. The largest impulse occurs at the first AP. Therefore, the best option is to divide the impulse at that location into smaller impulses if that impulse is large to be realized by the rocket engine.

The total coast time associated with the three-impulse base solution is 6.4738 days, calculated by adding up the coast times of the two phase-free arcs. With this total mission time, it is not feasible to have any phasing orbits added to the trajectory simply because this is the minimum mission time to reach the arrival point at the target orbit, and the available $TOF_p$ for phasing orbits is zero. Following Eq.~\eqref{eq:addkappaTf}, we add 1 orbital revolution of the target orbit to the total mission time (i.e., $\kappa = 1$) to demonstrate the possibility of adding phasing orbits. This way, the correct phasing of the departure and arrival remains the same. The orbital period of the target orbit is $T_f = 3.9191$ days; therefore, when we consider no revolutions on any segments of the trajectory (i.e., $N_0 = N_f = N_{\text{pf},1} = N_{\text{pf},2} = 0$), we have $TOF_p = 3.9191$ days. Using the feasibility relation given in Eq.~\eqref{eq:feas2}, the lower and upper bounds of the $\sum_{k=1}^{n_{p,1}}N_{k,1}$ are determined as $0.73 \leq \sum_{k=1}^{n_{p,1}}N_{k,1} \leq 58.09$. Assuming one revolution on each orbit, $0.73 \leq n_{p,1} \leq 58.09$, it is possible to add a minimum of 1 and a maximum of 58 orbits at the first AP. 
\begin{figure}[!htbp]
\begin{subfigure}{.5\textwidth}
  \centering
  \includegraphics[width=1\columnwidth]{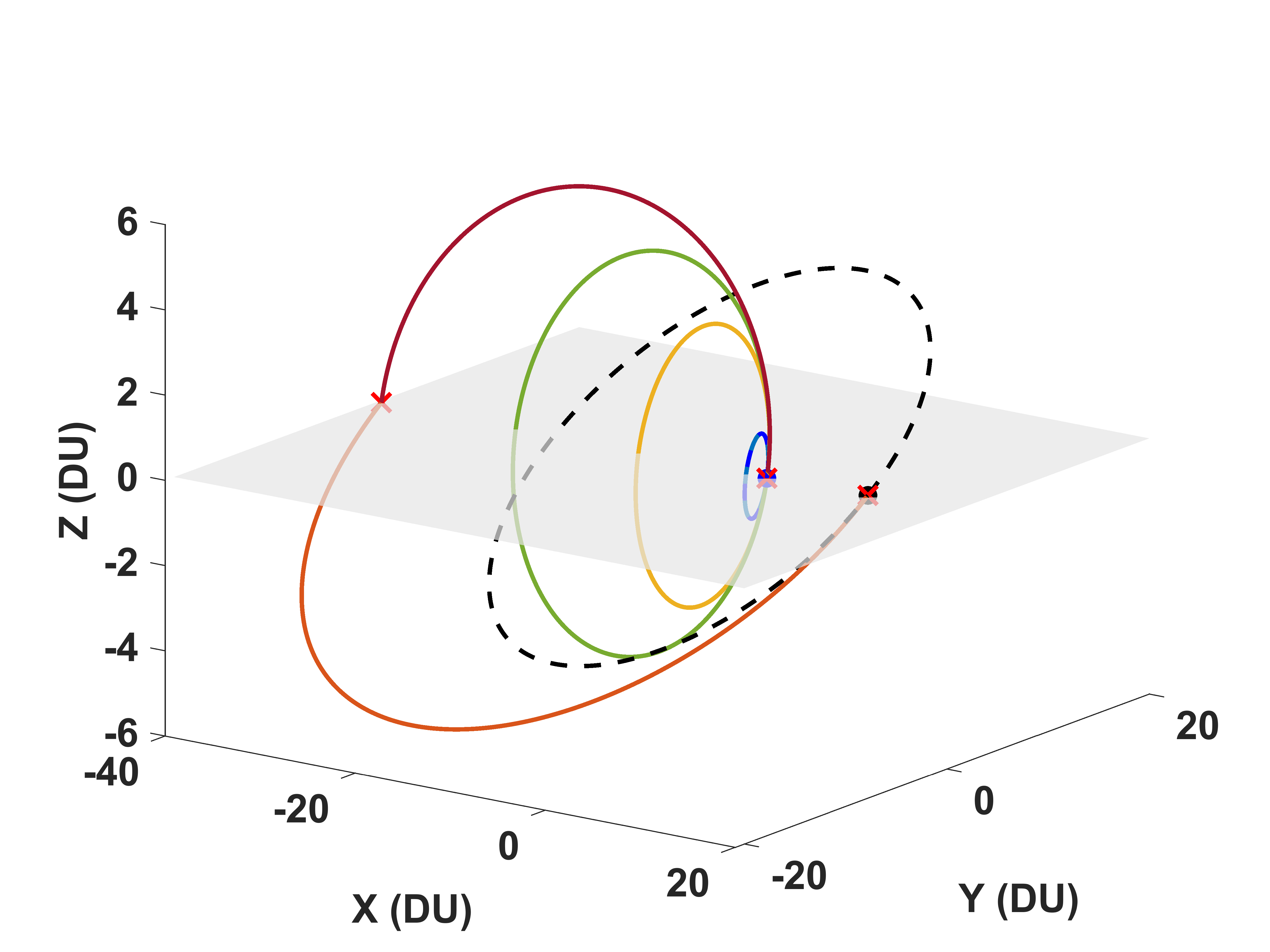}
  \caption{$N_0 = 1, N_{1,1} = 1,N_{1,2}= 1$}
  \label{fig:2ph_a}
\end{subfigure}%
\begin{subfigure}{.5\textwidth}
  \centering
  \includegraphics[width=1\columnwidth]{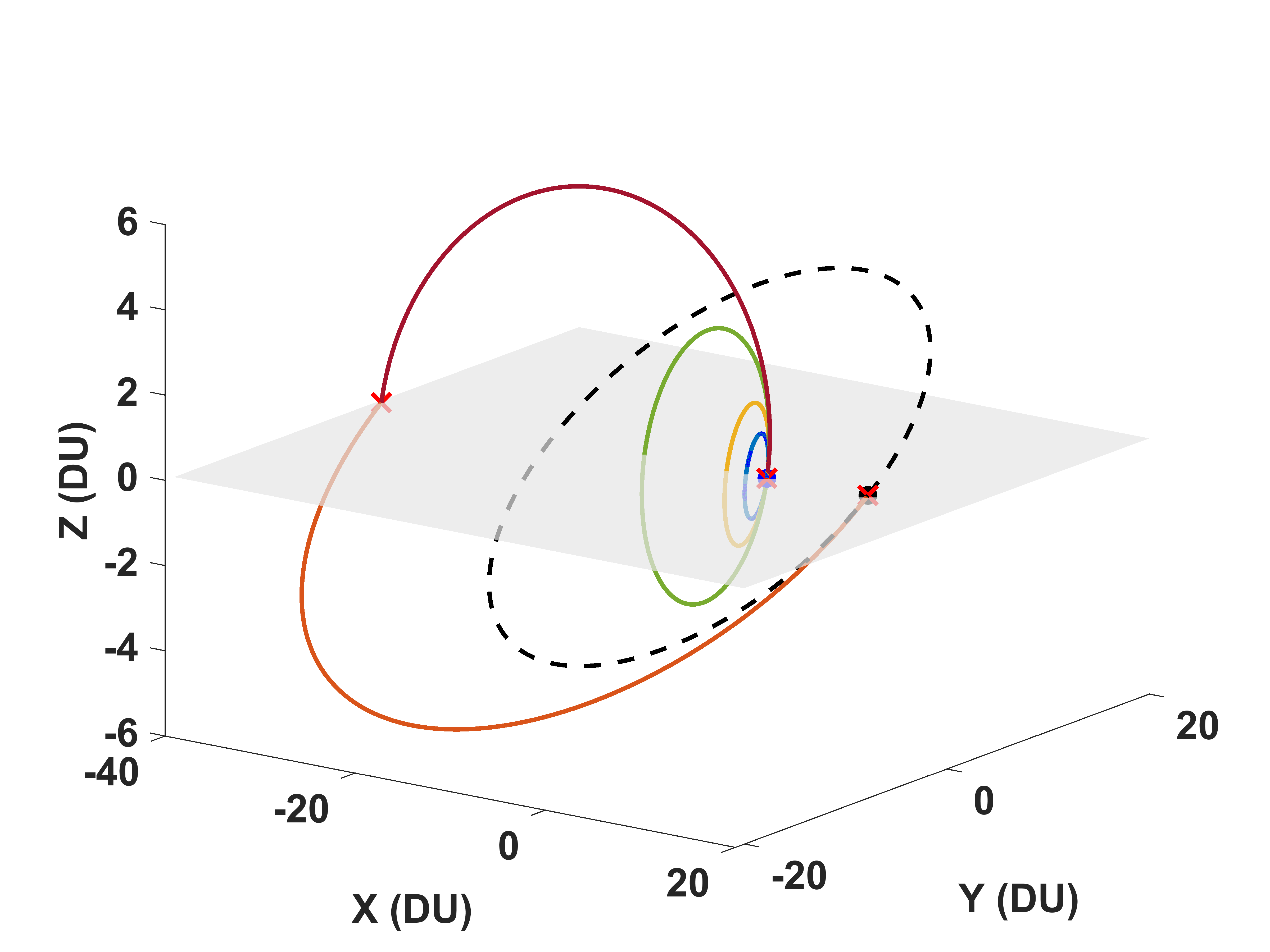}
  \caption{$N_0 = 2, N_{1,1} = 5,N_{2,1}= 3$}
  \label{fig:2ph_b}
\end{subfigure}
\begin{subfigure}{.5\textwidth}
  \centering
  \includegraphics[width=1\columnwidth]{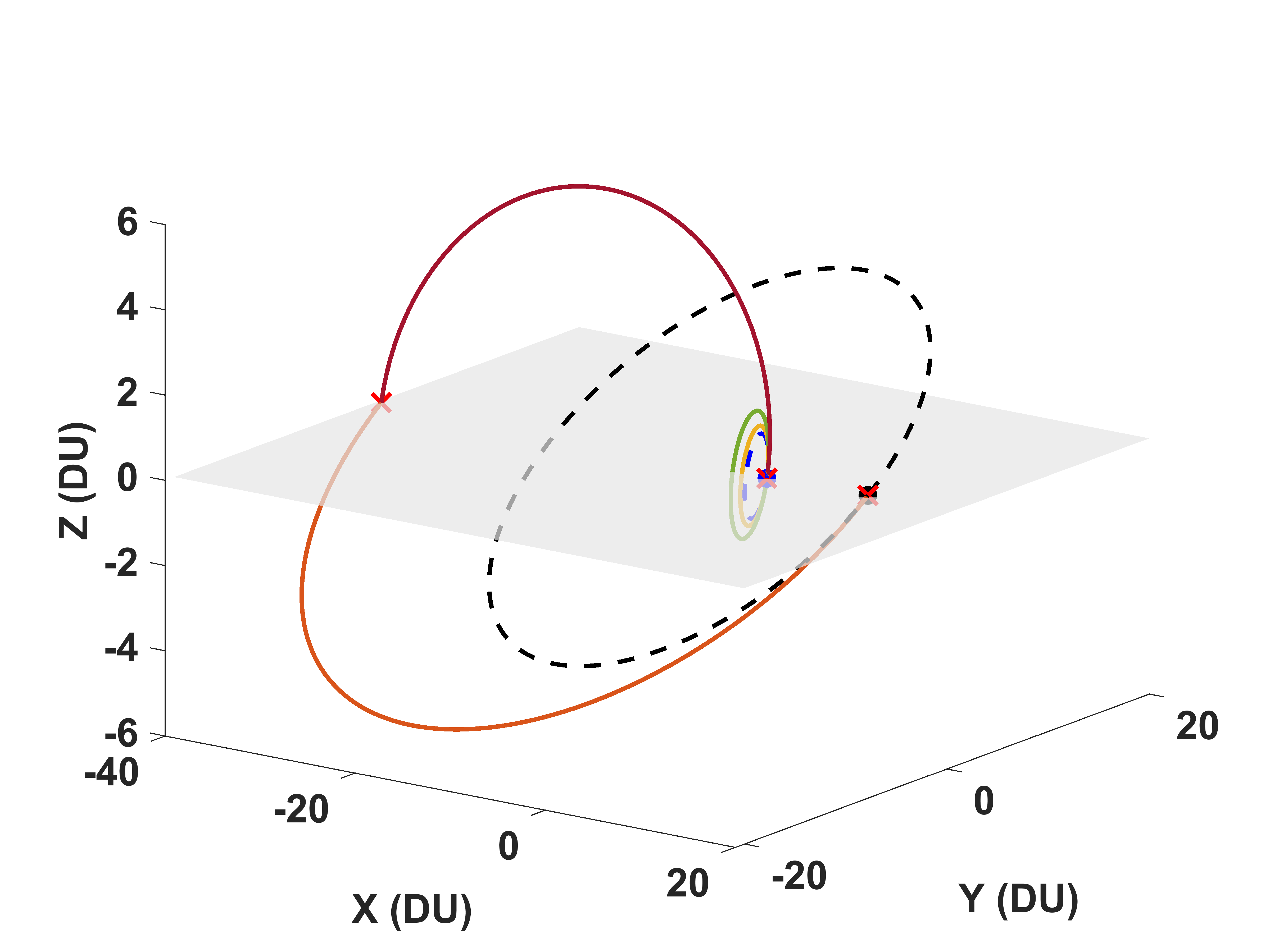}
  \caption{$N_0 = 0, N_{1,1} = 12,N_{2,1} = 20$}
  \label{fig:2ph_c}
\end{subfigure}
\begin{subfigure}{.5\textwidth}
  \centering
  \includegraphics[width=1\columnwidth]{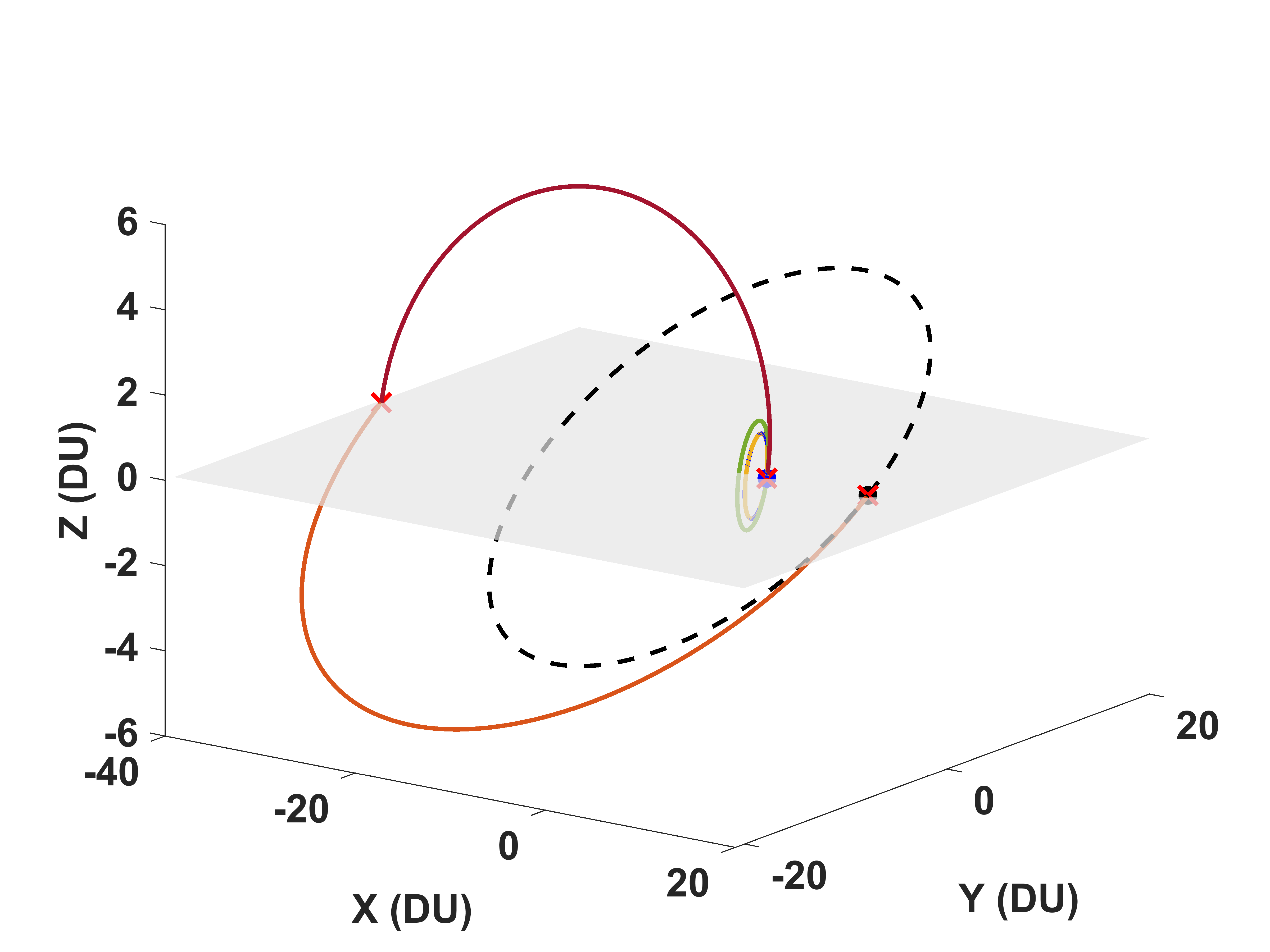}
  \caption{$N_0 = 0, N_{1,1} = 55,N_{2,1} = 2$}
  \label{fig:2ph_d}
\end{subfigure}
\caption{Geocentric problem: five-impulse solutions with various revolutions on phasing and initial orbits.}
\label{fig:2ph}
\end{figure}
As an example, we add two phasing orbits at the first AP. Obtained iso-impulse solutions are shown in Fig.~\ref{fig:2ph}. 
The three-impulse base solution structure is kept in all the iso-impulse solutions. The only change is the number of revolutions on phasing or initial orbits. Each solution family with specific $N_0, N_{1,1}, N_{2,1}$ values has infinitely many iso-impulse solutions since each phasing orbit is a continuous variable. 
Here, we show one solution among the infinitely many options from each family of solutions. The orbital periods of phasing orbits in Fig.~\ref{fig:2ph_a} are the largest among all four subplots since the total number of revolutions is smaller; thus, $TOF_p$ can be spent on phasing orbits with larger orbital period values. In Fig.~\ref{fig:2ph_b}, the solution space is decreased as the number of revolutions on each segment is increased. It means that the largest values of orbital periods for each phasing orbit are much smaller. In Figs.~\ref{fig:2ph_c} and \ref{fig:2ph_d}, the number of revolutions on phasing orbits is even greater. Therefore, the area of the solution envelope is shrunk and orbits are almost identical to the initial orbit. These solutions satisfy the feasibility constraint since $\sum_{k=1}^{n_{p,1}}N_{k,1} =  32$ and $\sum_{k=1}^{n_{p,1}}N_{k,1} =  57$ are both less than 58. Recall that the solution envelopes belong to the fourth layer of the solution hierarchy depicted in Fig.~\ref{fig:vis_abs}.

\subsection{Iso-Impulse Solutions with Multiple Impulse Anchor Positions}
In this section, we consider the $\Delta v$-allocation problem at two or three APs simultaneously. We discuss how the feasibility equations can be leveraged to determine the bounds of the feasible solution space for each case and show the solutions obtained. 

\subsubsection{Solution Families with Two Impulse Anchor Position Solutions (Earth-to-Dionysus Problem)}
\begin{figure} [htbp!]
    \centering
    \includegraphics[width=0.6\linewidth]{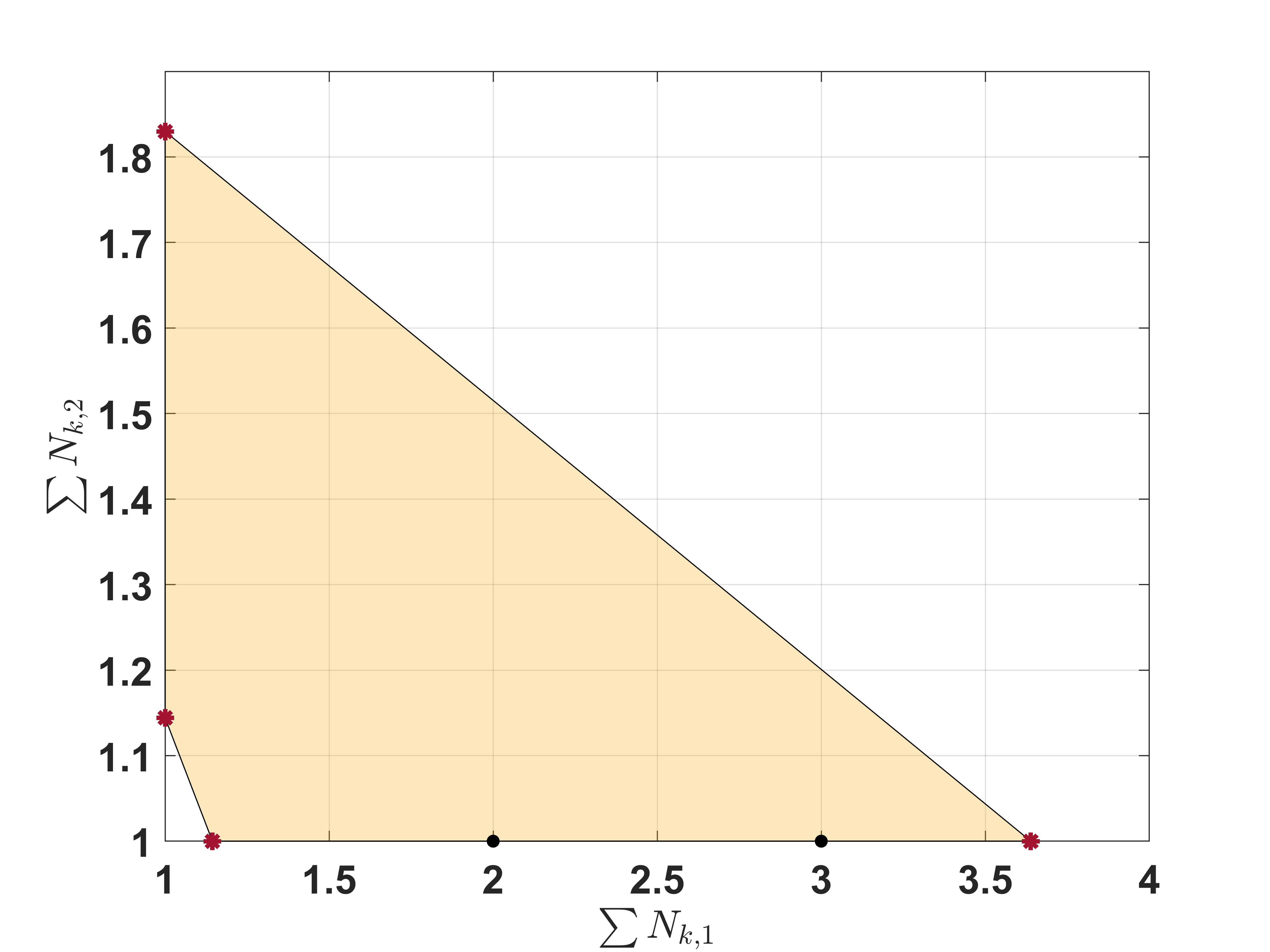}
    \caption{Earth-to-Dionysus problem: feasible solution space with two impulse APs in terms of $\sum N_{k,1}$ and $\sum N_{k,2}$.}
    \label{fig:2AP_solspace_E2D}
\end{figure}

Problem parameters are taken from Ref.~\cite{saloglu_existence_2023}.
Feasibility analysis is performed using Eq.~\eqref{eq:feas2}. We summarized how to obtain the two-impulse base solution in Sec.~\ref{sec:genmultimp}. Here, we use the two-impulse base solution to show the solutions obtained by considering both APs for $\Delta v$ allocation. From Eq.~\eqref{eq:feas1} with $T_0 = 365.25$ days and $T(\alpha_{k,1} = 1) = 1161.47$ days, we calculated $\sum_{k=1}^{n_{p,2}}N_{k,2} \leq 1.83$ as the upper bound. Therefore, we can only add one phasing orbit at the second impulse AP. Then, we use Eq.~\eqref{eq:feas2} to obtain $\sum_{k=1}^{n_{p,1}}N_{k,1} \leq 3.64$ when $\sum_{k=1}^{n_{p,2}}N_{k,2} = 1$. From this analysis, the total number of revolutions that can occur on the first and the second impulse APs ($\sum_{k=1}^{n_{p,1}}N_{k,1}$, $\sum_{k=1}^{n_{p,2}}N_{k,2}$) are obtained as  (2, 1) and (3, 1). 
The feasible solution space in terms of $\sum_{k=1}^{n_{p,1}}N_{k,1}$ and $\sum_{k=1}^{n_{p,2}}N_{k,2}$ is given in Fig.~\ref{fig:2AP_solspace_E2D}. The red star markers represent the upper and lower bounds when the value on the other axis is kept at 1. Even though we provide the upper bounds in the above discussion, we use the lower bounds to fully represent the feasible solution space, as shown in Fig.~\ref{fig:2AP_solspace_E2D}. The lower bounds are calculated as $\sum_{k=1}^{n_{p,2}}N_{k,2} \geq 1.14$ and $\sum_{k=1}^{n_{p,1}}N_{k,1} \geq 1.14$. Every integer combination in this yellow-shaded area represents a feasible solution (cf. layer 2), and each feasible solution can generate multiple solution families (cf. layer 3) with different combinations of numbers of revolutions at intermediate phasing orbits. The only feasible solutions that have infinitely many solutions are (2,1) and (3,1), which are denoted as black dots in Fig.~\ref{fig:2AP_solspace_E2D}. The feasibility we define here is the feasibility of having infinitely many solutions. The bounds of the ``$\sum$'' operator on the labels are dropped.

When considering two impulse APs, we find two solution families for the case with one phasing orbit at each AP. Also, there are three families for two phasing and one phasing orbits, and one family for three phasing and one phasing at the two APs. 
When we consider only one impulse AP, where the largest impulse is applied, we found 134 solution families~\cite{saloglu_existence_2023}. When two impulse APs are considered, additional phasing orbits at the second AP decrease the area of the corresponding solution envelope. Therefore, the number of solution families for two impulse APs is less than the single AP version. 
 
\begin{figure}[!htbp]
\vspace{-4mm}
\begin{subfigure}{.5\textwidth}
  \centering
  \includegraphics[width=1\columnwidth]{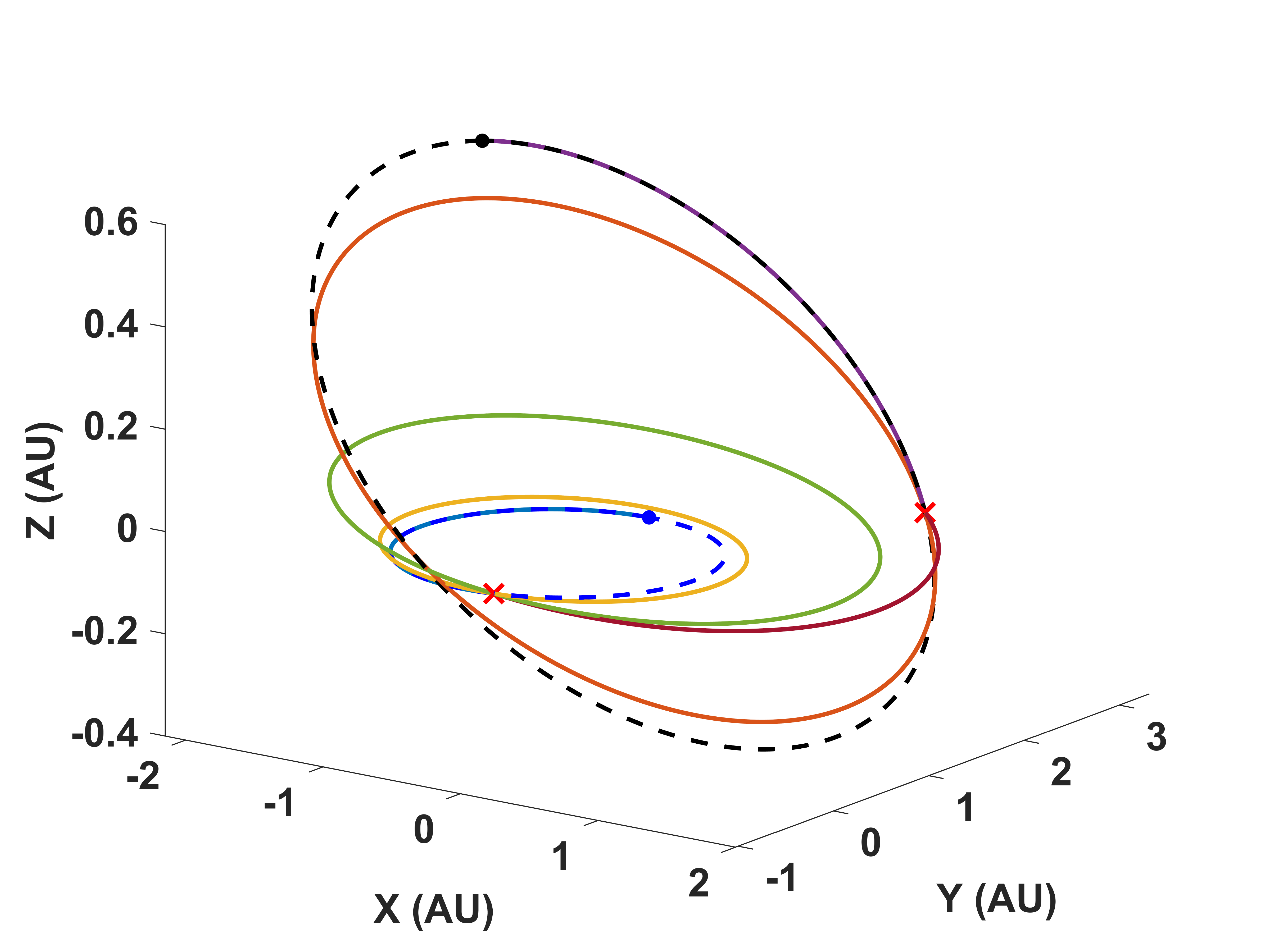}
  \caption{Three-dimensional view of the trajectory}
  \label{fig:traj_2_1_3d}
\end{subfigure}%
\begin{subfigure}{.5\textwidth}
  \centering
  \includegraphics[width=1\columnwidth]{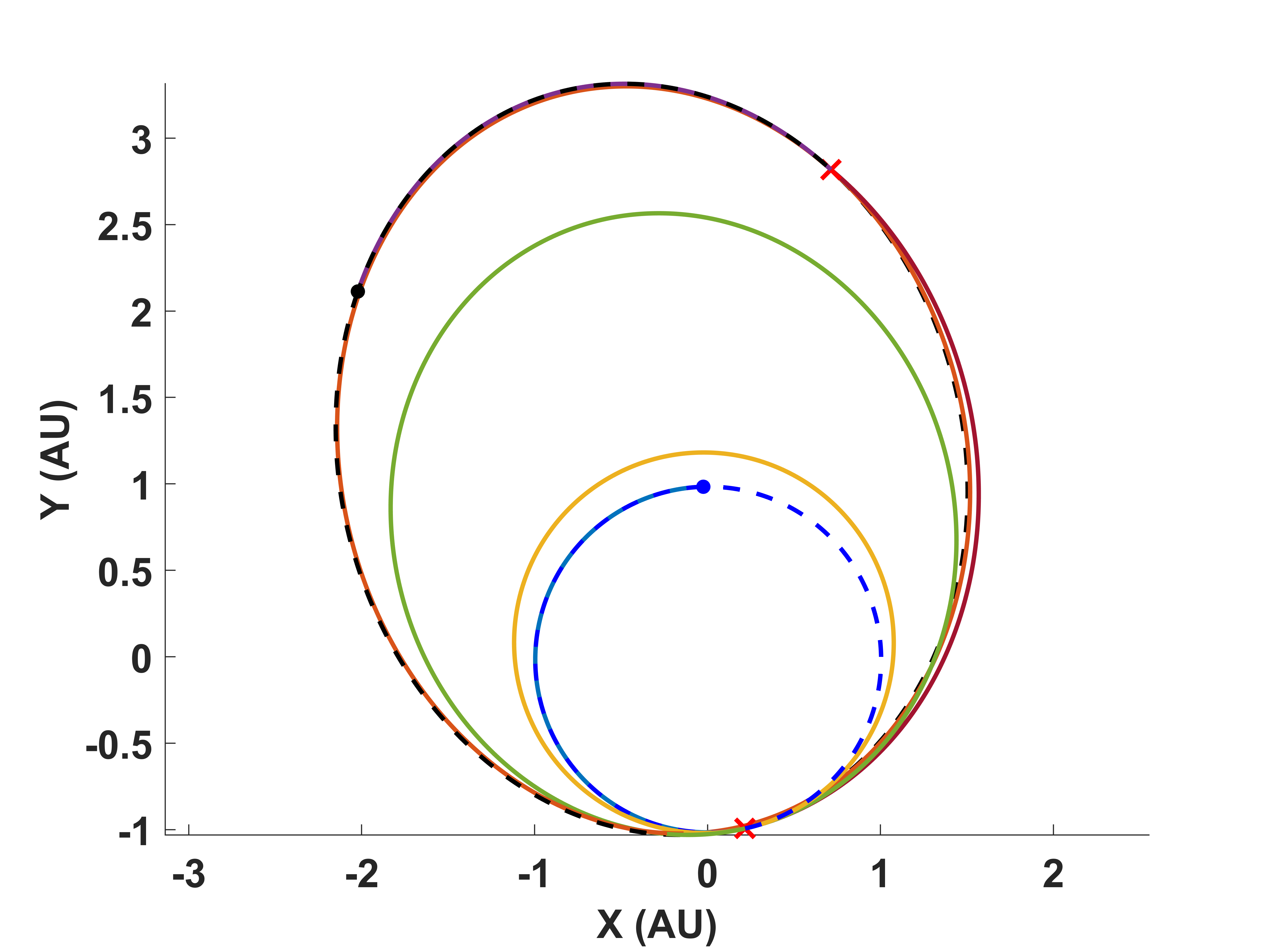}
  \caption{$x-y$ view of the trajectory}
  \label{fig:traj_2_1_2d}
\end{subfigure}
\caption{Earth-to-Dionysus problem: solutions with two impulse APs ($N_{1,1} = N_{2,1} = N_{1,2} = 1$). }
\vspace{-5mm}
\label{fig:traj_2_1}
\end{figure}

\begin{figure}[!htbp]
\begin{subfigure}{.5\textwidth}
  \centering
  \includegraphics[width=1\columnwidth]{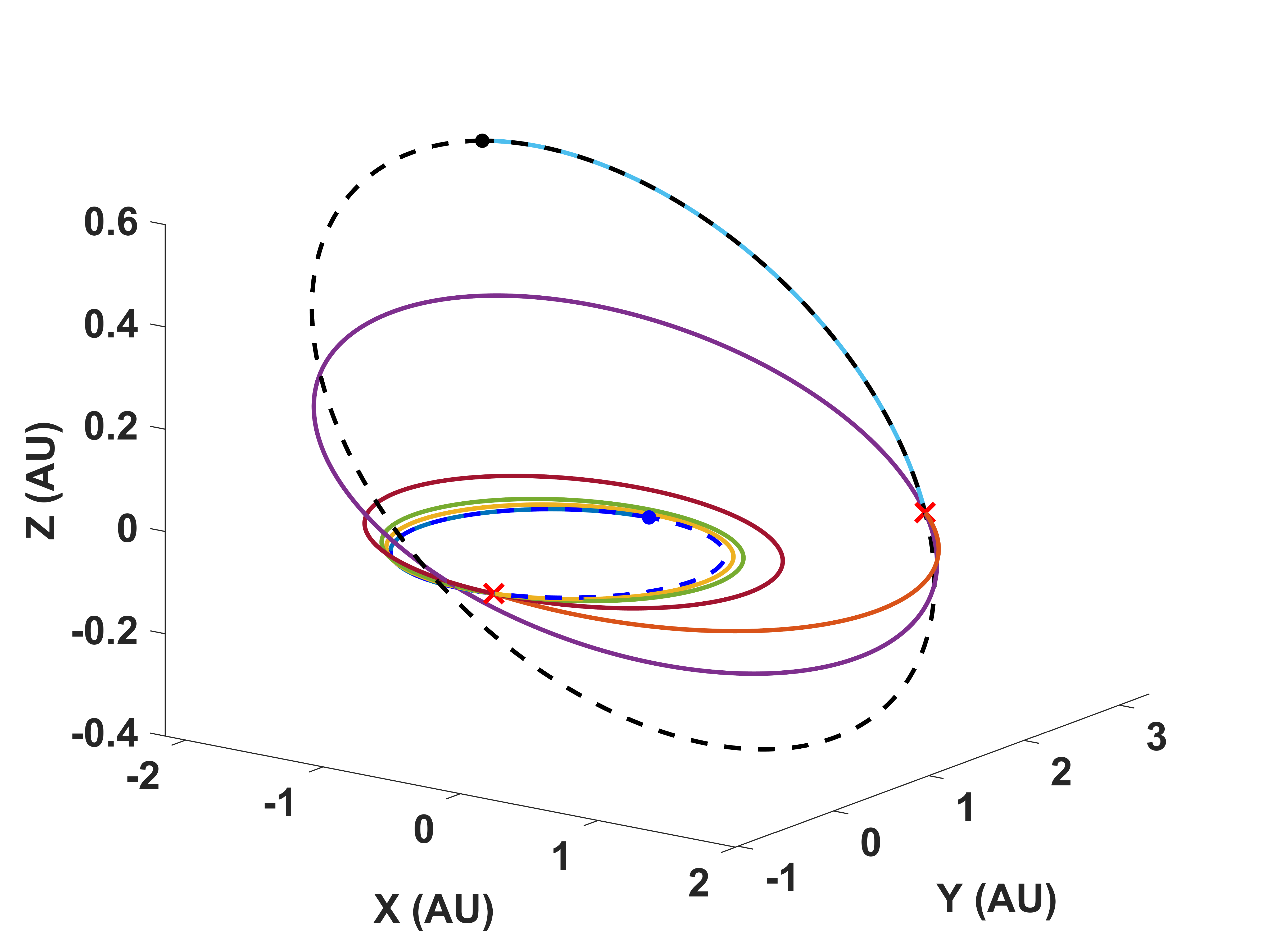}
  \caption{Trajectory in three-dimensions}
  \label{fig:traj_3_1_3d}
\end{subfigure}%
\begin{subfigure}{.5\textwidth}
  \centering
  \includegraphics[width=1\columnwidth]{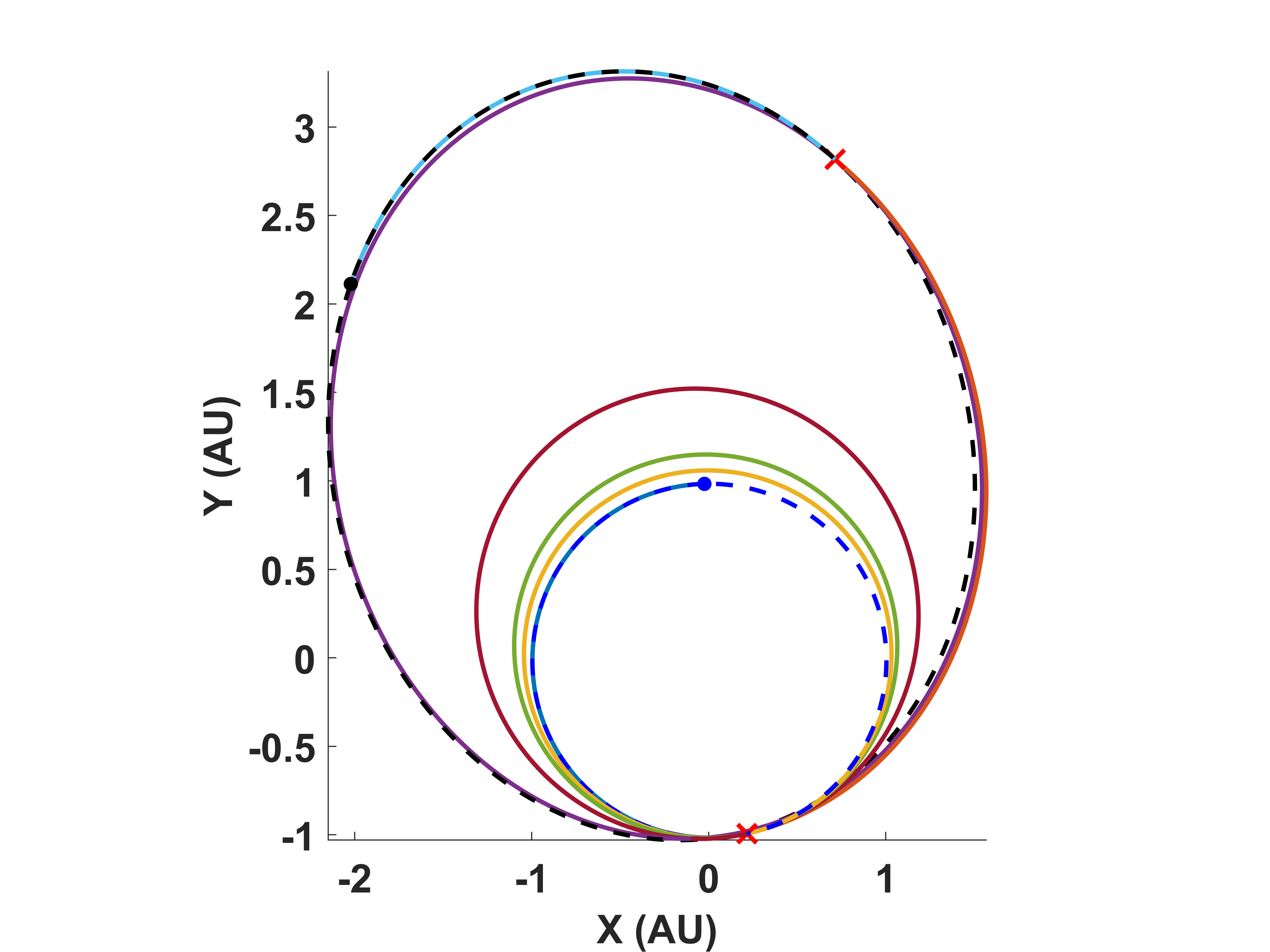}
  \caption{Trajectory in two-dimensions}
  \label{fig:traj_3_1_2d}
\end{subfigure}
\caption{Earth-to-Dionysus problem optimal trajectory with two impulse APs ($N_{1,1} = N_{2,1} = N_{3,1} = N_{1,2} = 1$). }
\label{fig:traj_3_1}
\vspace{-5mm}
\end{figure}

 From the generated solution families, an example solution family with two phasing orbits at the first impulse AP and one phasing orbit at the second impulse AP is shown in Fig.~\ref{fig:traj_2_1}. 
 There are no revolutions on the initial and target orbits, i.e., $N_0 = N_f = 0$. Two phasing orbits at the first impulse AP are shown in colors yellow and green. The last phasing orbit at the second impulse AP is shown in orange. The first and the last impulses are divided into three and two smaller-magnitude impulses, respectively. 
 Since the orbital period of the last phasing orbit is large, the solution space for the phasing orbits at the first impulse AP is reduced. Another solution family with three phasing orbits at the first impulse AP and one phasing orbit at the second impulse AP is shown in Fig.~\ref{fig:traj_3_1}.

 Three phasing orbits are shown in yellow, green, and red at the first impulse AP. Since the solution space is small, orbital periods of the phasing orbits are close to each other. The phasing orbit at the second impulse AP is shown in purple. The total impulse at that position is divided into two smaller-magnitude impulses. 

\subsubsection{Solution Families with Two Impulse Anchor Position (Geocentric Case Study)}
If we have a three-impulse base solution, one can still consider choosing two impulse APs to divide the impulses. We use the base solution presented in Fig.~\ref{fig:3imp_traj} to demonstrate this case. The assumption is that there is no phasing orbit at the last impulse location. First, we perform a feasibility analysis to determine the upper and lower bounds of the total number of revolutions at each AP. The lowest value of the total number of revolutions at the two APs is 1, i.e., $\sum_{k=1}^{n_{p,1}}N_{k,1} = 1$ and $\sum_{k=1}^{n_{p,2}}N_{k,2} = 1$. To determine the upper and lower bounds on the total number of revolutions, we assume $\sum_{k=1}^{n_{p,1}}N_{k,1} = 1$ and solve for the bounds of $6.3096 \leq \sum_{k=1}^{n_{p,2}}N_{k,2} \leq 7.2970$ with $T_0 = 0.0675$ days, $T_{k,2}(\alpha_{k,2} = 0) = 5.3616$, and $T_{k,1} = 5.3616$ using Eq.~\eqref{eq:feas2}. Similarly, when $\sum_{k=1}^{n_{p,2}}N_{k,2} = 1$, $6.3096\leq\sum_{k=1}^{n_{p,1}}N_{k,1}\leq 501.4617$. Therefore, the upper bound on the total number of revolutions on the first and second impulse AP is $7$ and $501$, respectively. 
\begin{figure}[h!]
    \centering
    \includegraphics[width=0.6\linewidth]{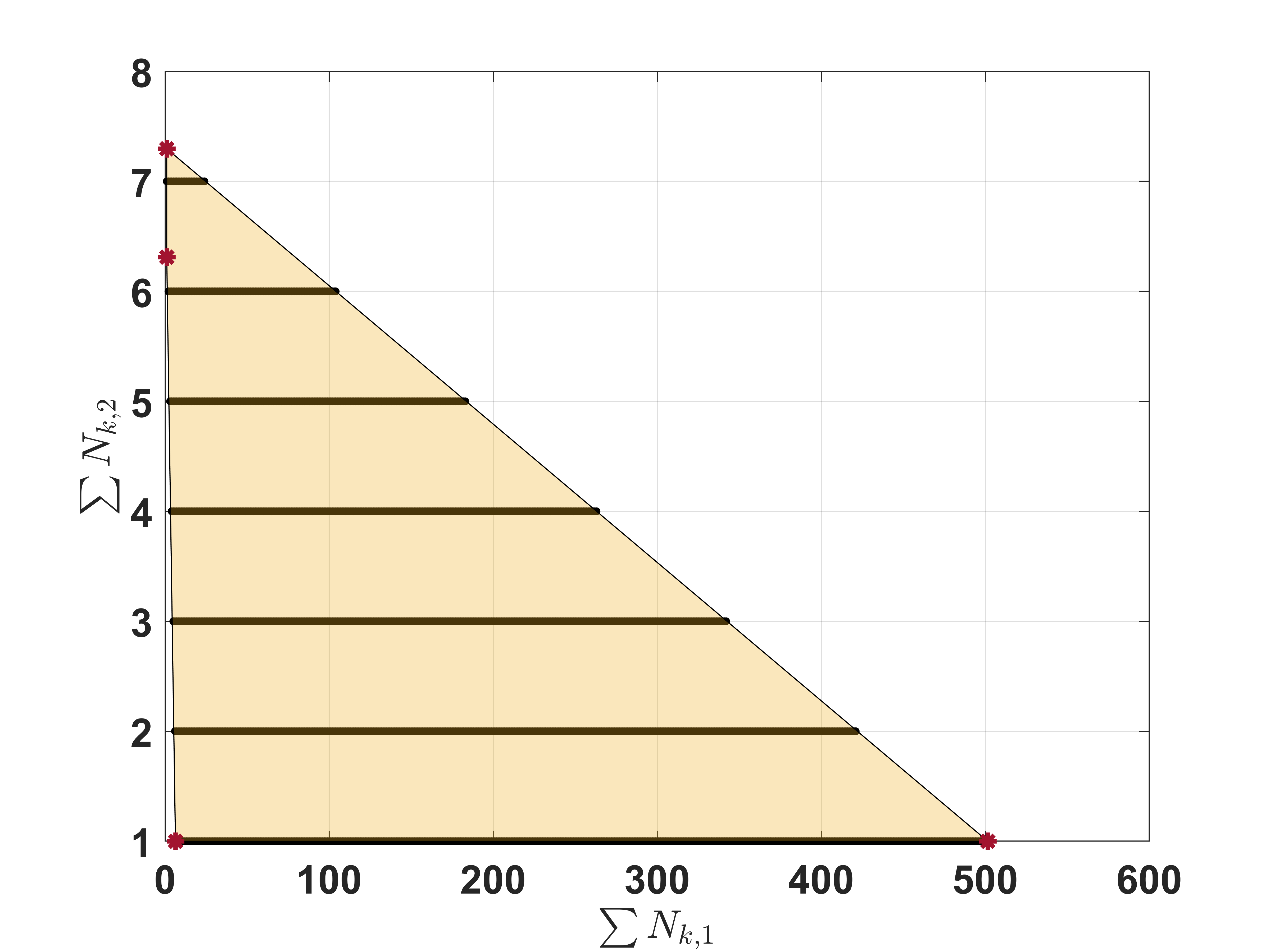}
    \caption{Geocentric problem: feasible solutions of the two impulse APs in terms of $\sum_{k=1}^{n_{p,1}}N_{k,1}$ and $\sum_{k=1}^{n_{p,2}}N_{k,2}$.}
    \label{fig:2AP_solspace}
\end{figure}
If one revolution on each phasing orbit is assumed, one can add 7 and 501 orbits at maximum to each APs independently. The total value of the number of revolutions can be distributed to different numbers of phasing orbits, with different combinations of revolutions on each. The feasible solution space for this example is plotted in Fig.~\ref{fig:2AP_solspace} using the upper and lower bound values of $\sum_{k=1}^{n_{p,1}}N_{k,1}$ and $\sum_{k=1}^{n_{p,2}}N_{k,2}$. There are 1817 feasible solutions, shown with black dot markers. Red star markers show the bounds of the variables on each axis. Every integer combination of $\sum N_{k,1}$ and $\sum N_{k,2}$ within the yellow shaded area corresponds to a solution family with infinitely many solutions. 
\begin{figure}[!htbp]
\begin{subfigure}{.5\textwidth}
  \centering
  \includegraphics[width=1\columnwidth]{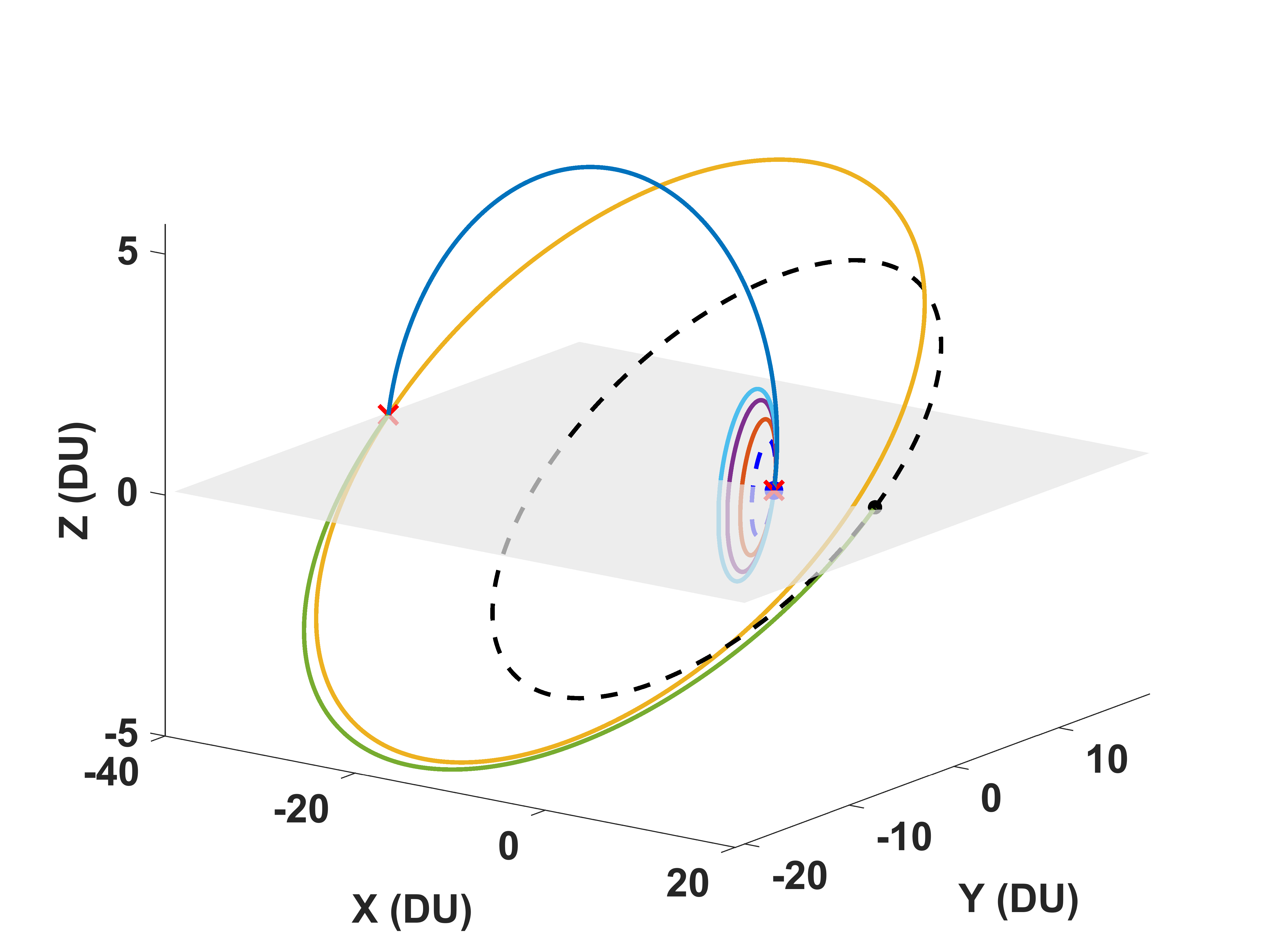}
  \caption{$N_{1,1} = N_{2,1} = N_{3,1} = 12; N_{1,2} = 1$}
  \label{fig:traj_2ap_1orb_a}
\end{subfigure}%
\begin{subfigure}{.5\textwidth}
  \centering
  \includegraphics[width=1\columnwidth]{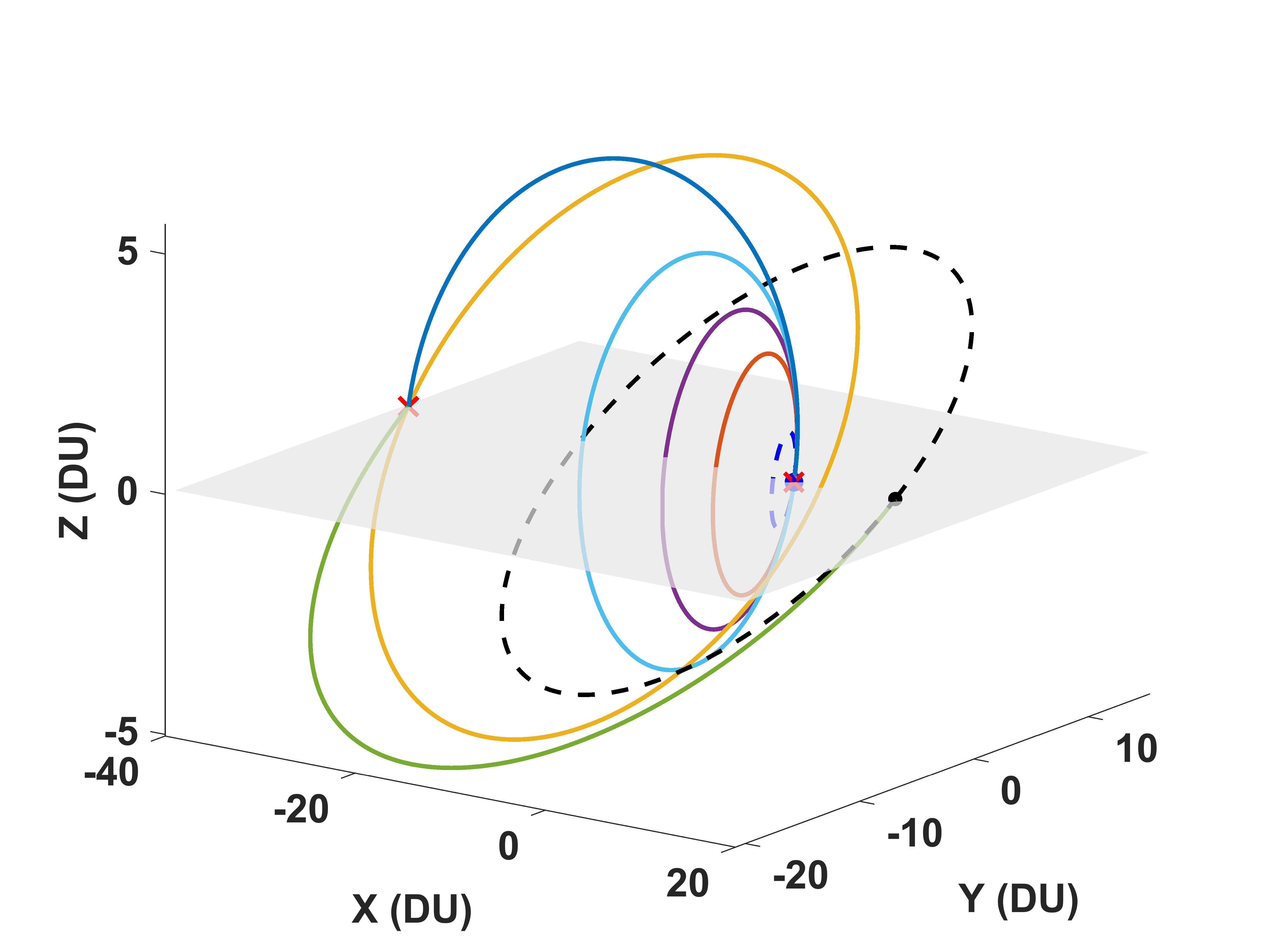}
  \caption{$N_{1,1} =  N_{3,1} = 5, N_{2,1} = 10; N_{1,2} = 1$}
  \label{fig:traj_2ap_1orb_b}
\end{subfigure}
\caption{Geocentric trajectory with two impulse APs, three and one phasing orbits at the first and second APs.  }
\label{fig:traj_2ap_1orb}
\vspace{-5mm}
\end{figure}
One feasible solution is $\sum_{k=1}^{n_{p,1}}N_{k,1} = 36$ and $\sum_{k=1}^{n_{p,2}}N_{k,2} = 1$, which can be realized with a different number of phasing orbits at the first AP. An example solution family is shown in Fig.~\ref{fig:traj_2ap_1orb_a}. Orange-, purple-, and blue-colored orbits are the phasing orbits at the first AP, and yellow is the phasing orbit at the second AP. The second AP phasing orbit period is large, and the possible period range at the first AP is limited. Also, having 12 revolutions on each first AP orbit makes the corresponding solution envelope smaller, which makes these phasing orbit periods small. One of the cases for  $\sum_{k=1}^{n_{p,1}}N_{k,1} = 20$ and $\sum_{k=1}^{n_{p,2}}N_{k,2} = 1$ is plotted in Fig.~\ref{fig:traj_2ap_1orb_b}. The orbit periods at the second AP are smaller and the number of revolutions is smaller at the first AP orbits, thus, the orbits at the first AP  have larger orbital periods compared to the previous case.

\begin{figure}[!htbp]
\begin{subfigure}{.5\textwidth}
  \centering
  \includegraphics[width=1\columnwidth]{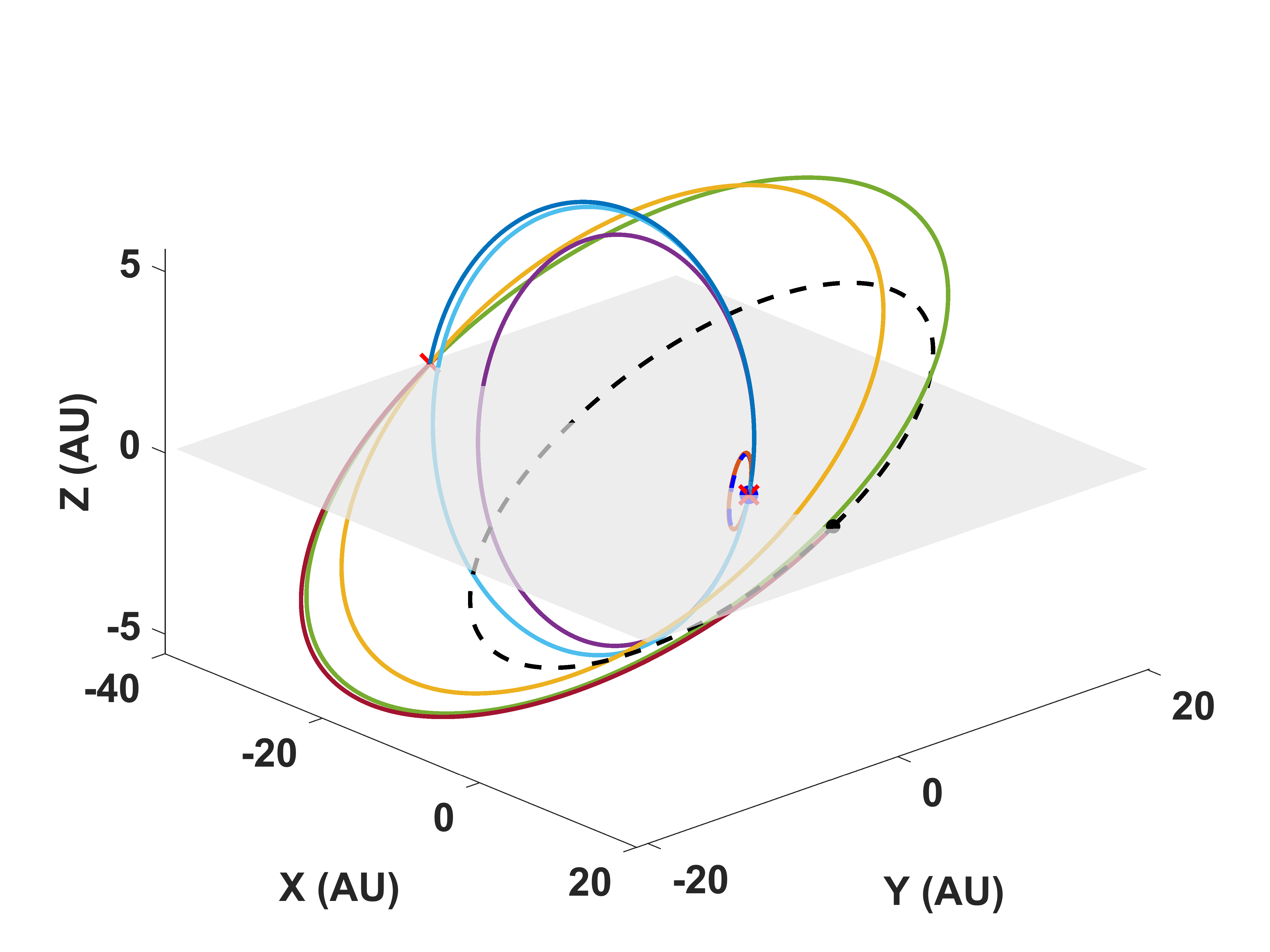}
  \caption{$N_{1,1} = N_{2,1} = N_{3,1} = N_{1,2} = 2, N_{2,2} = 1$}
  \label{fig:traj_2ap_2orb_a}
\end{subfigure}%
\begin{subfigure}{.5\textwidth}
  \centering
  \includegraphics[width=1\columnwidth]{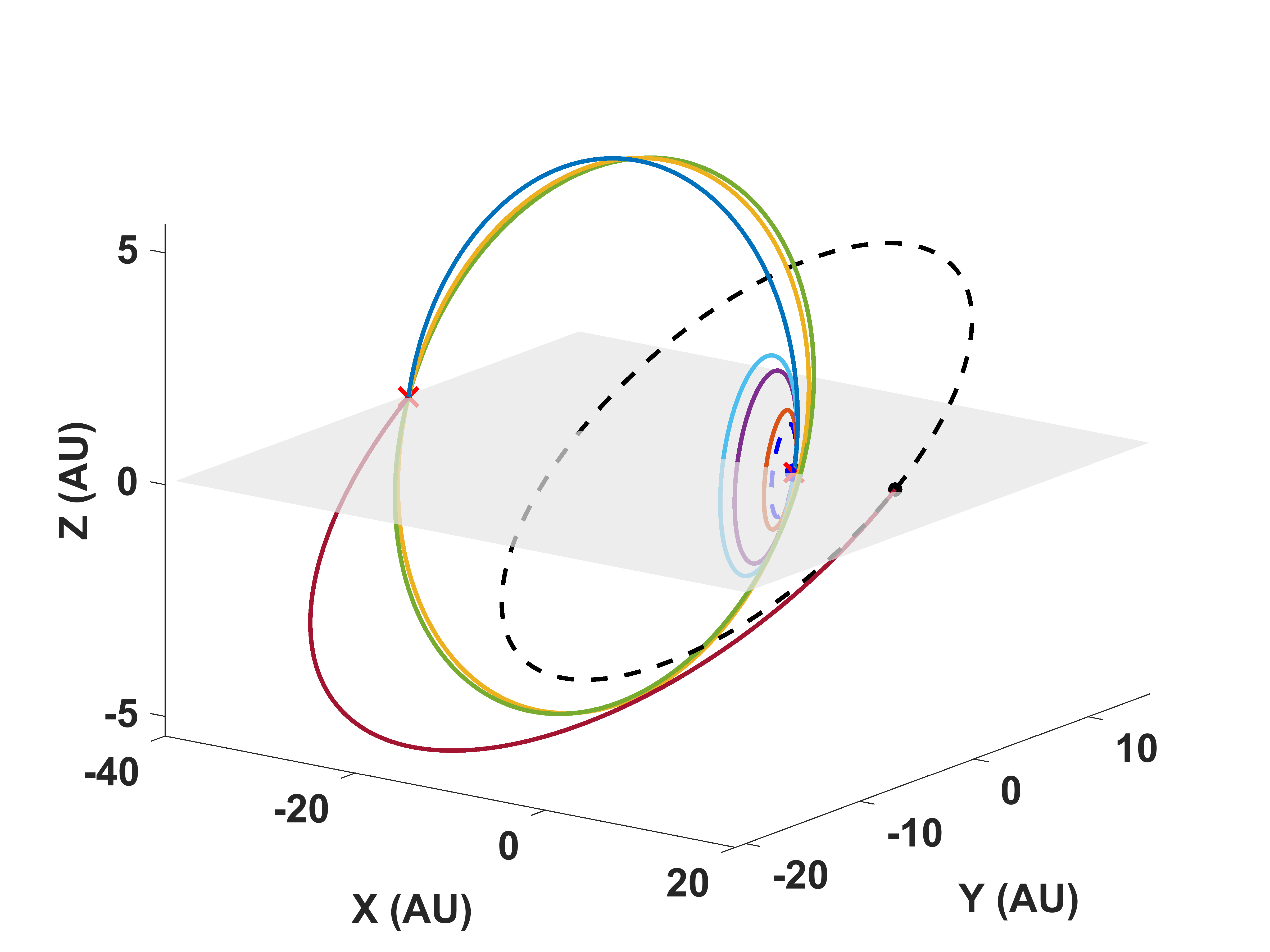}
  \caption{$N_{1,1} = N_{2,1} = N_{3,1} = 1; N_{1,2} = 3, N_{2,2} = 4$}
  \label{fig:traj_2ap_2orb_b}
\end{subfigure}
\caption{Geocentric trajectory with two impulse APs, three and two phasing orbits at the APs. }
\label{fig:traj_2ap_2orb}
\vspace{-5mm}
\end{figure}
If we add two phasing orbits at the second AP, it is still possible to recover minimum-$\Delta v$ solutions as the maximum number of orbits we can add is 7 at the second AP with the assumption of one revolution. Two different solution families are shown in Fig.~\ref{fig:traj_2ap_2orb}. Orange-, purple-, and blue-colored orbits are at the first AP, and green- and yellow-colored orbits are at the second AP. In Fig.~\ref{fig:traj_2ap_1orb_a}, the orange orbit has almost the same period as the initial orbit, which means a small $\Delta v$ is applied. A large semi-major change happened between the orange and purple orbits, since the largest impulse is applied at the first AP for increasing the semi-major axis. The number of revolutions at the first AP orbits is 2 and is relatively small compared to those in Fig.~\ref{fig:traj_2ap_1orb}. Therefore,  
the orbital periods are larger. In Fig.~\ref{fig:traj_2ap_2orb_b}, the number of revolutions at the second AP orbits is higher, resulting in smaller orbital period values compared to the case in Fig.~\ref{fig:traj_2ap_2orb_a}.

\subsubsection{Solution Families with Three Impulse Anchor Positions (Geocentric Case Study)} \label{sec:GEO3APs}
\begin{figure}[htbp!]
    \centering
    \includegraphics[width=0.7\linewidth]{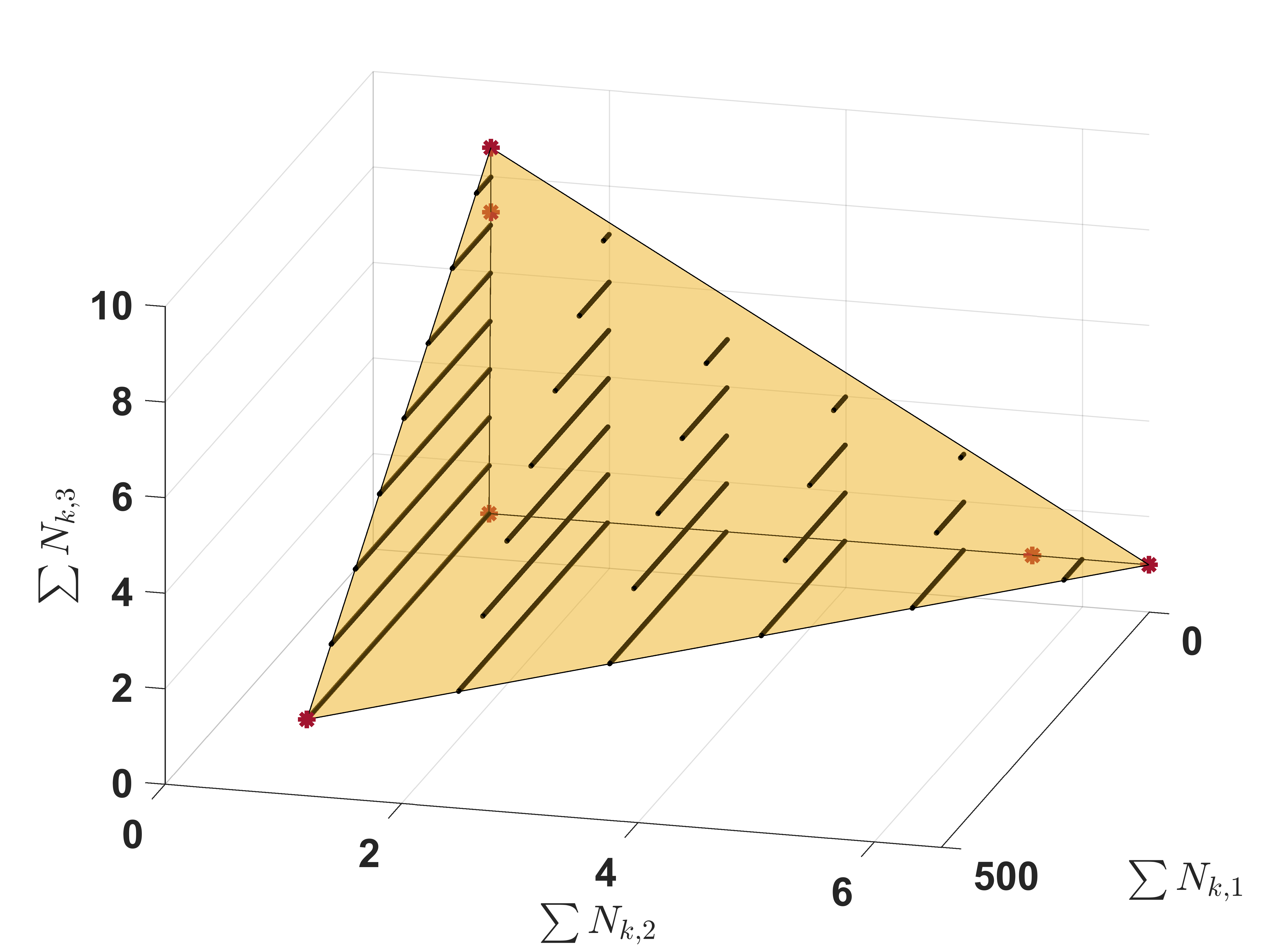}
    \caption{Geocentric problem: feasible solution space with three impulse APs.}
    \label{fig:3AP_solspace}
\end{figure}

Starting from a three-impulse base solution, 
Eq.~\eqref{eq:feas2} is used to determine the solutions' feasibility. Additionally, the feasible solution space (now a polytope) can be bounded to provide upper and lower bounds on total number of revolutions. If there are three impulse APs, we assume there will be at least one phasing orbit at each AP. Keeping two of the $\sum_{k=1}^{n_{p,i}}N_{k,i}$ values at 1, the upper and lower bounds of the remaining term are obtained from Eq.~\eqref{eq:feas1}. For the geocentric example, feasible solutions are denoted as black dot markers within the yellow-shaded polytope in Fig.~\ref{fig:3AP_solspace}. There are 4747 feasible solutions and each solution can be used to obtain different numbers of solution families. 
The upper bound on each total number of revolutions determines the maximum value of the number of phasing orbits that can be added to the trajectory at the respective AP. For instance, $5.5786 \leq \sum_{k=1}^{n_{p,1}}N_{k,1} \leq 443.3749$,  $5.5786 \leq \sum_{k=1}^{n_{p,2}}N_{k,2} \leq 6.5660$, and $7.2639 \leq \sum_{k=1}^{n_{p,3}}N_{k,3} \leq 8.6147$. The maximum number of orbits that can be added is 443, 6, and 8 independently at the respective APs.

\begin{figure}[!htbp]
\begin{subfigure}{.5\textwidth}
  \centering
  \includegraphics[width=1\columnwidth]{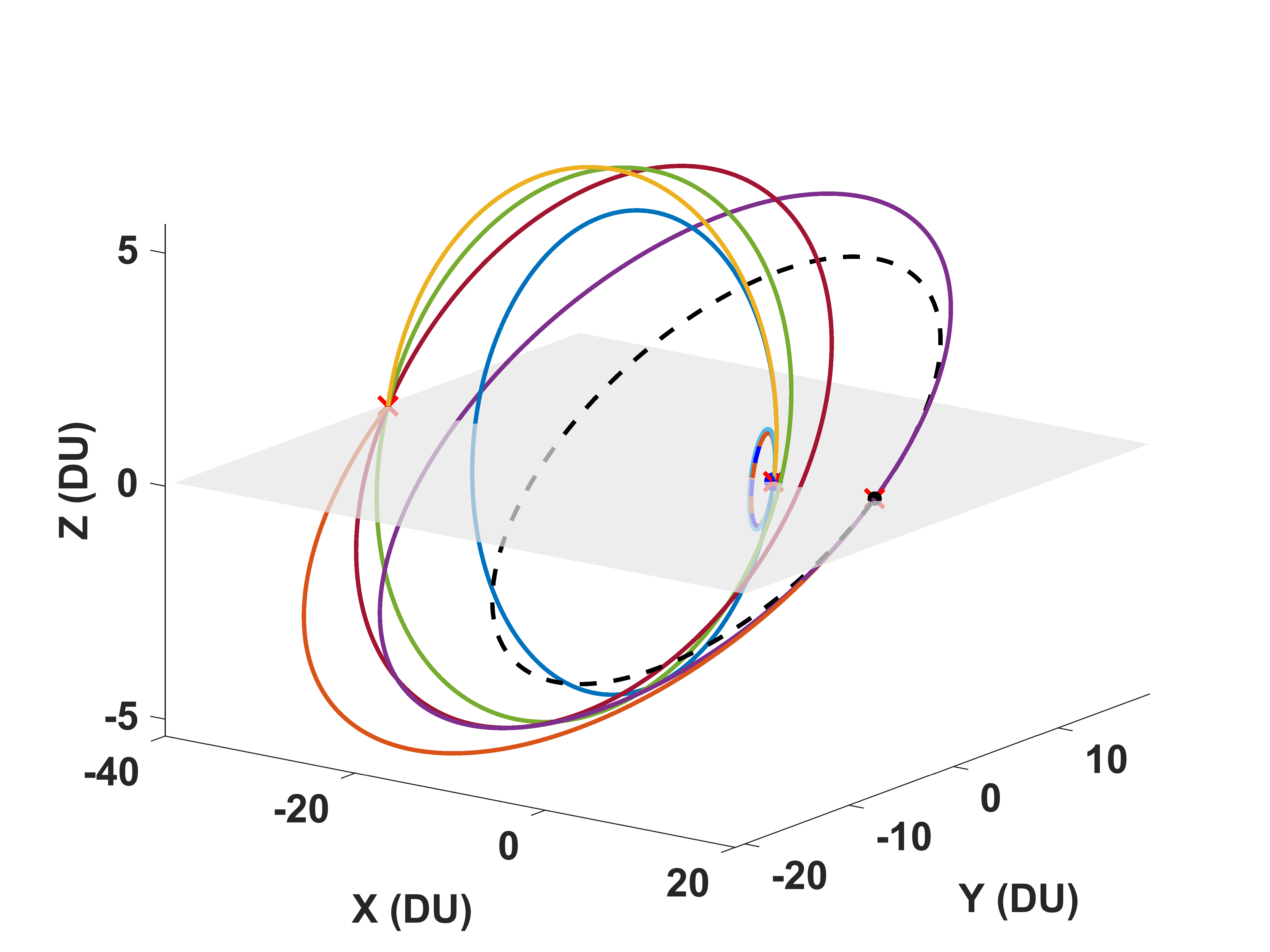}
  \caption{$N_{1,1} = N_{2,1} = N_{3,1} = N_{4,1} = 1$; \\ $N_{1,2} = N_{2,2} = 1; N_{1,3} = 4$}
  \label{fig:traj_3ap_a}
\end{subfigure}%
\begin{subfigure}{.5\textwidth}
  \centering
  \includegraphics[width=1\columnwidth]{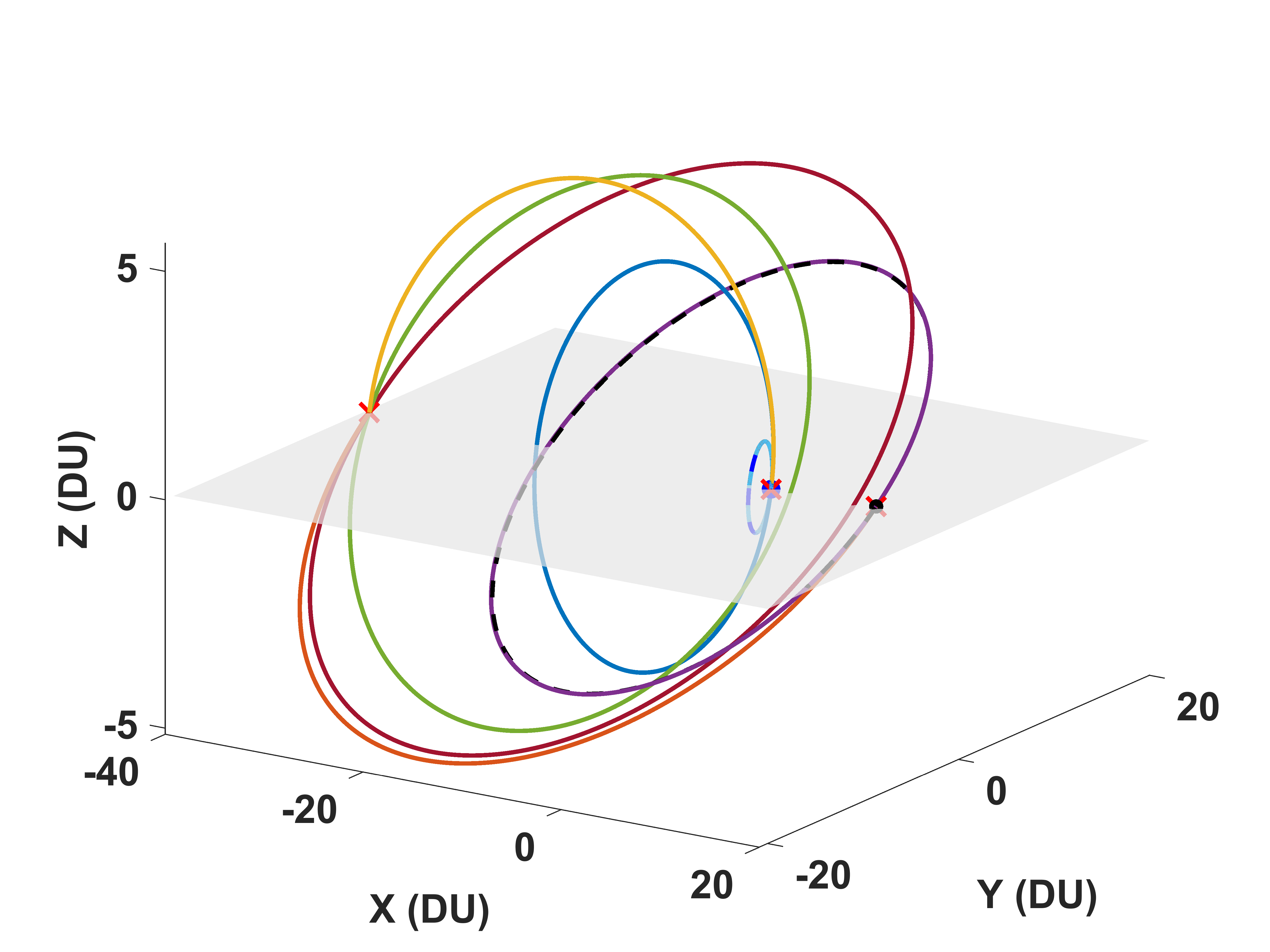}
  \caption{$N_{1,1} = N_{2,1} = 6, N_{3,1} = 8, N_{4,1} = 4; \\N_{1,2} = 3, N_{2,2} = 1; N_{1,3} = 1$}
  \label{fig:traj_3ap_b}
\end{subfigure}
\begin{subfigure}{.5\textwidth}
  \centering
  \includegraphics[width=1\columnwidth]{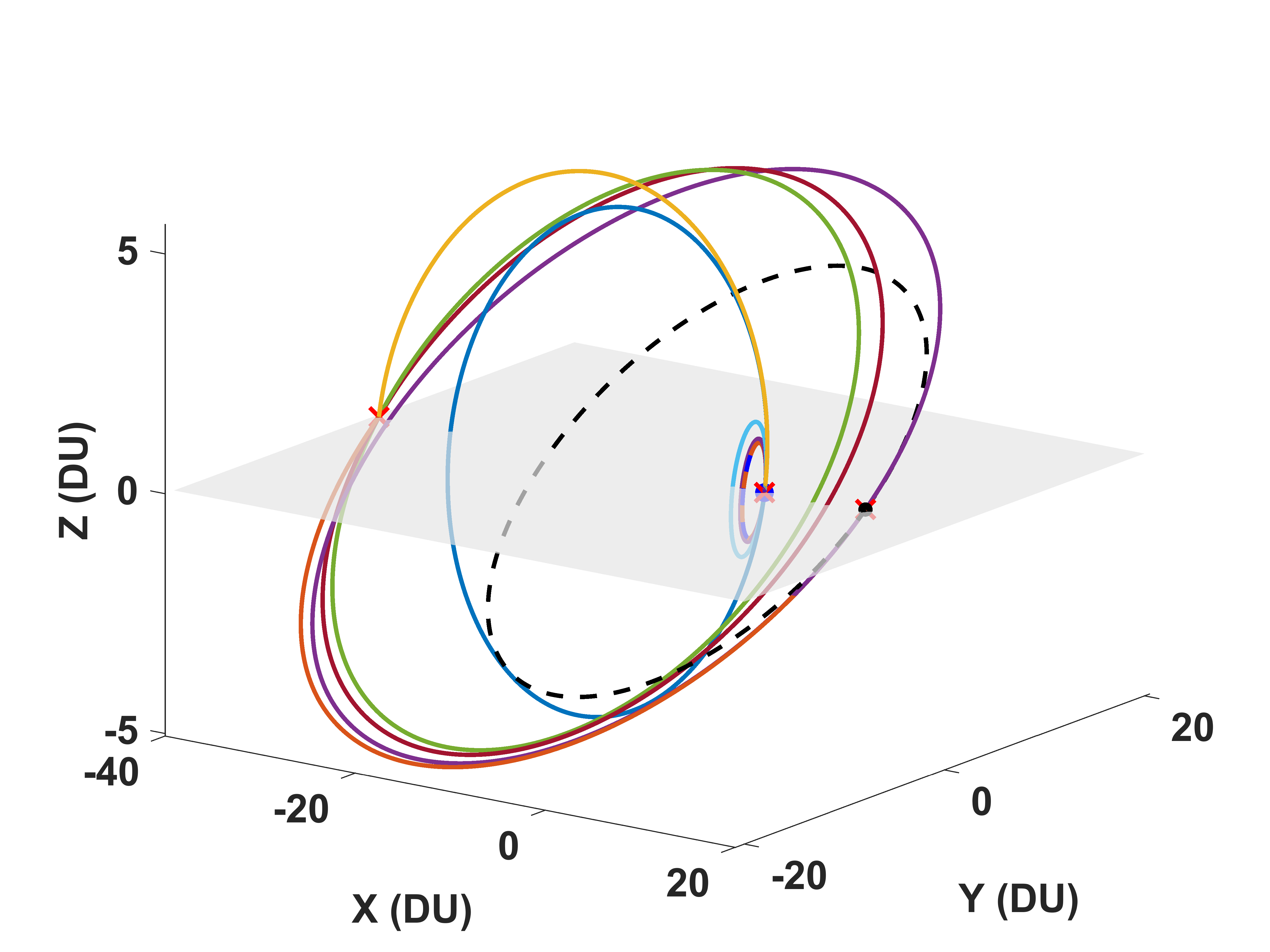}
  \caption{$N_{1,1} = N_{2,1} = 2, N_{3,1} = 4, N_{4,1} = 1; \\N_{1,2} = 1, N_{2,2} = 2; N_{1,3} = 2$}
  \label{fig:traj_3ap_c}
\end{subfigure}
\begin{subfigure}{.5\textwidth}
  \centering
  \includegraphics[width=1\columnwidth]{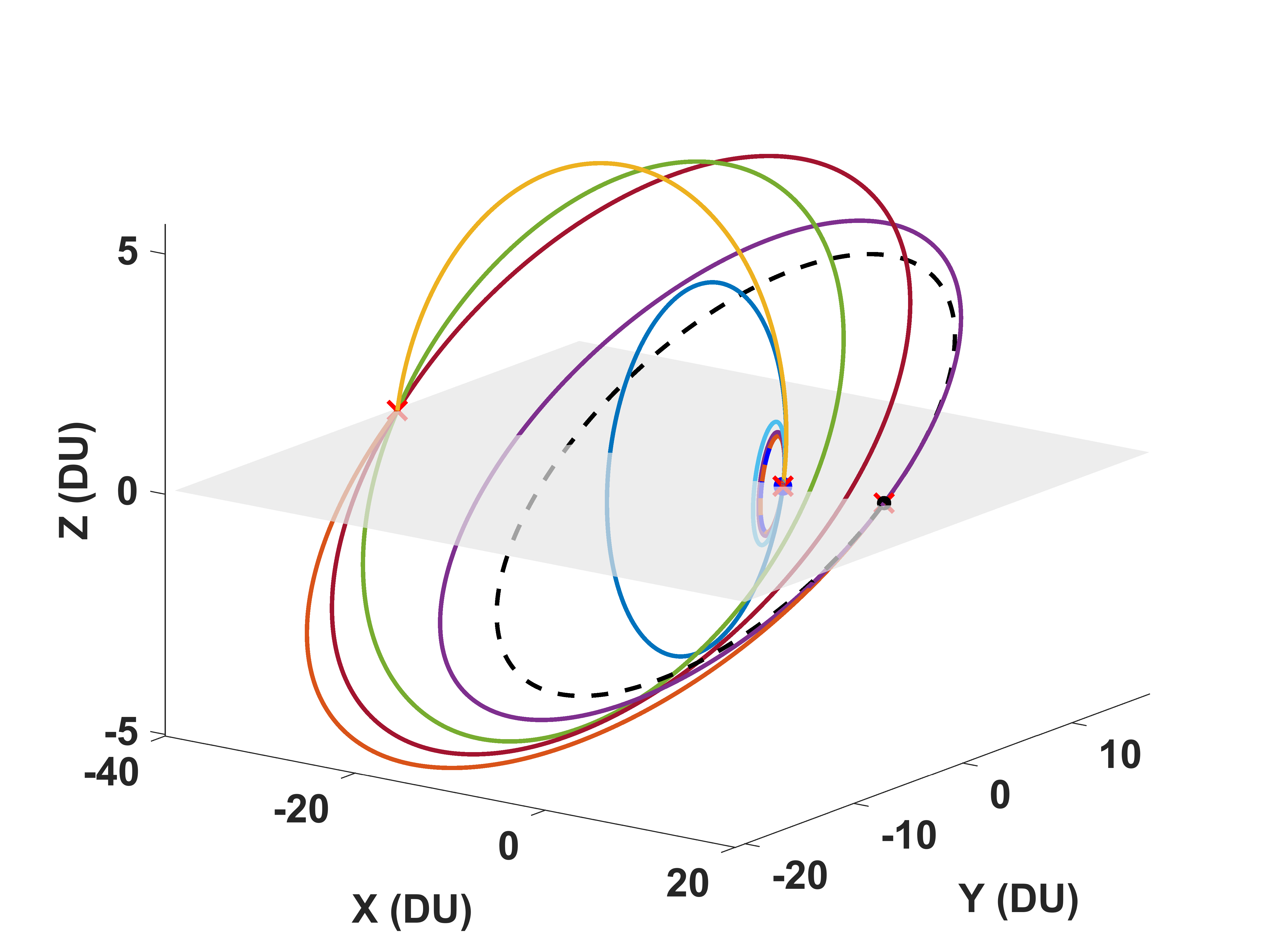}
  \caption{$N_{1,1}=5, N_{2,1} = 3, N_{3,1} =  N_{4,1} = 1; \\N_{1,2} = 2, N_{2,2} = 3; N_{1,3} = 1$}
  \label{fig:traj_3ap_d}
\end{subfigure}
\caption{Geocentric solutions with three impulse APs and four, two, and one phasing orbits at the APs.}
\label{fig:traj_3ap}
\end{figure}

Four different solution families are presented in Fig.~\ref{fig:traj_3ap}. Four, two, and one orbit are added to the first, second, and third APs. Green- and red-colored orbits are at the second AP and the purple orbit is at the third AP. The addition of the orbit at the third AP decreased the area of the solution envelope, which made the first AP orbits' period smaller. Also, the blue-colored orbit period at the first AP is larger than the remaining orbits at the first AP; so, phasing orbits cannot have large orbital periods. All the solutions have the same minimum $\Delta v$ value and each family has infinitely many solutions. 

\subsection{Solution Envelopes for Iso-Impulse Solutions}
\subsubsection{Solution envelopes for the Earth-to-Dionysus Problem}
\begin{figure}[!htbp]
\begin{subfigure}{.5\textwidth}
  \centering
  \includegraphics[width=1\columnwidth]{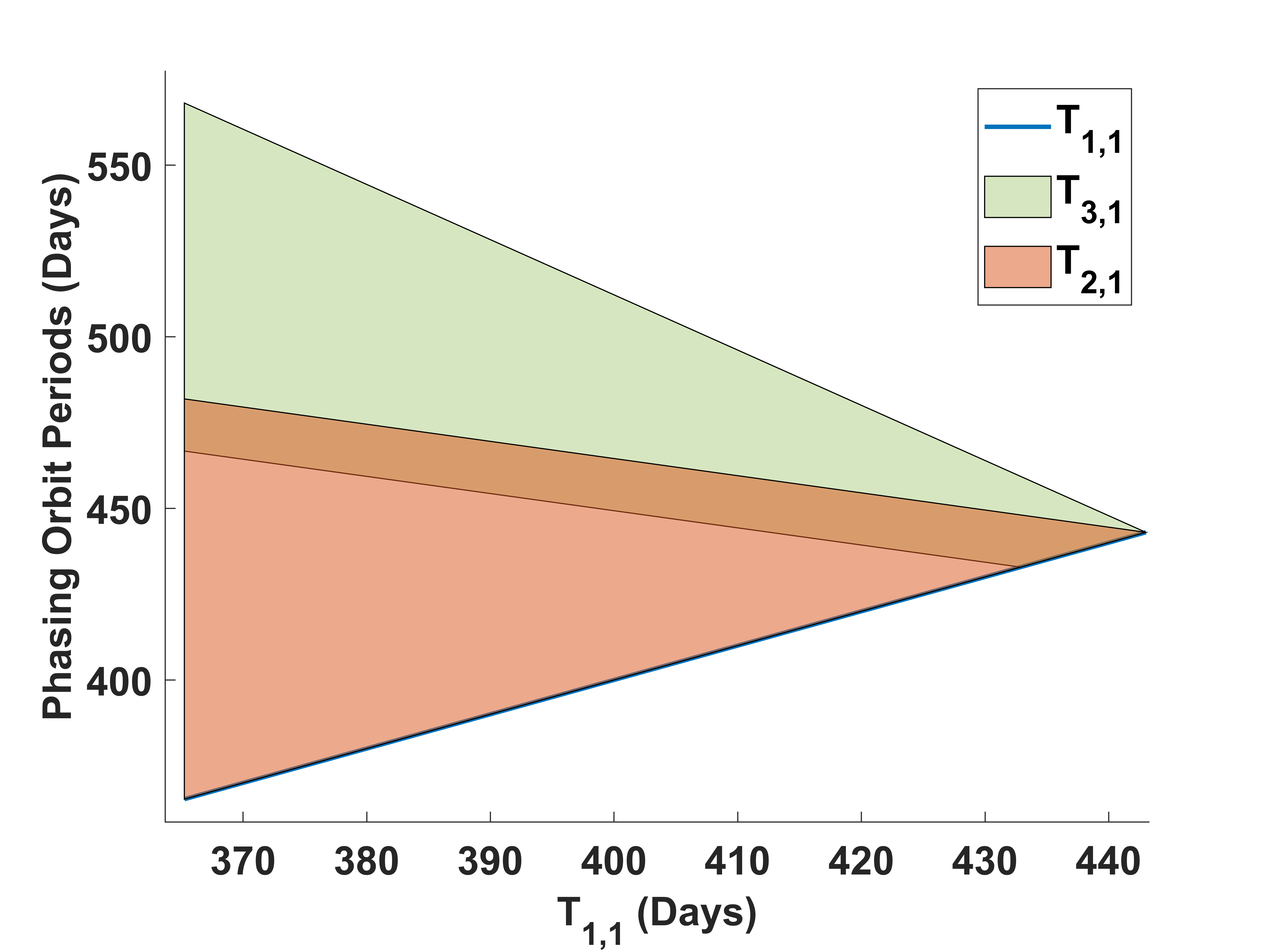}
  \caption{First AP}
  \label{fig:env_2ap_a}
\end{subfigure}%
\begin{subfigure}{.5\textwidth}
  \centering
  \includegraphics[width=1\columnwidth]{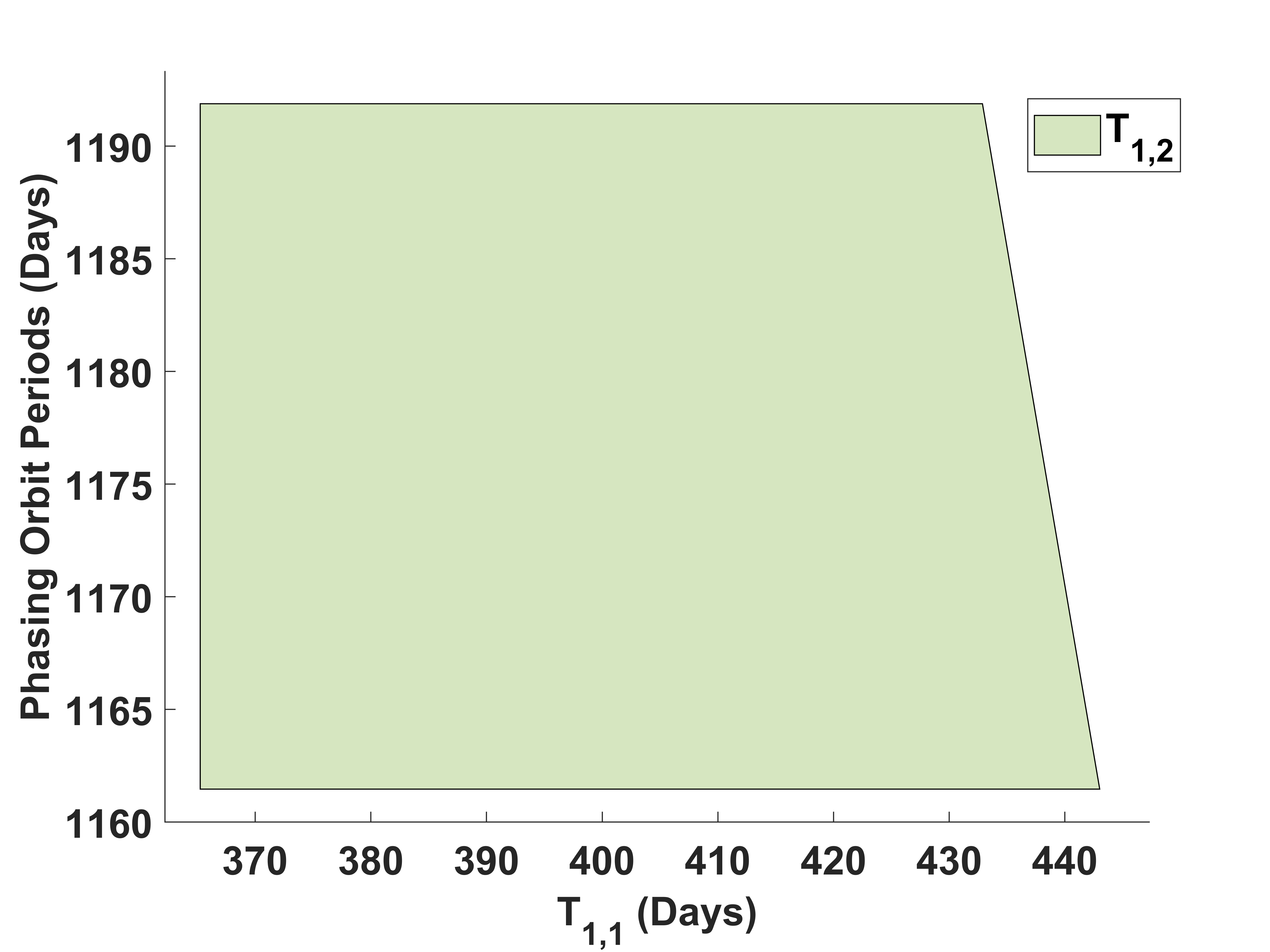}
  \caption{Second AP}
  \label{fig:env_2ap_b}
\end{subfigure}
\caption{Earth-to-Dionysus problem: solution envelopes with $N_{1,1} = N_{2,1} = N_{3,1} = N_{1,2} = 1$. }
\label{fig:env_2ap}
\end{figure}

We generate the solution envelopes for the Earth-to-Dionysus problem with three phasing orbits at the first AP and one phasing orbit at the second AP. For this solution family, one of the solutions is shown in Fig.~\ref{fig:traj_3_1}. The steps to determine the corner points are explained in Section~\ref{sec:sol_env}. The obtained solution envelopes are shown in Fig.~\ref{fig:env_2ap}. Solution envelopes are plotted with respect to the first phasing orbit orbital period, $T_{1,1}$. In Fig.~\ref{fig:env_2ap_a}, $T_{1,1}$ solutions create a line. Orange and green polygons represent the solution envelopes for the orbits at the first AP. When $T_{1,1} = T_0 = 365.25$ days, we first determine the range of values $T_{2,1} \in [365.25, 481.88]$ days and $T_{3,1} \in [466.67, 568.11]$ days, which are found using Eqs.~\eqref{eq:env_min_ap1} and \eqref{eq:env_max_ap1}. All the solutions converge to the rightmost corner when $T_{1,1,\text{max}} = 443$ days using Eq.~\eqref{eq:env_maxT1_ap1}. When $T_{1,1} = 365.25$ days, minimum value of $T_{2,1,\text{min}} \neq 365.25$ days, there is one more corner for the solution envelope of $T_{2,1}$, where $T_{1,1} = T_{3,1} = 432.90$ days found using Eq.~\eqref{eq:env_corner2_ap1}. In Fig.~\ref{fig:env_2ap_b}, the solution envelope of $T_{1,2}$ is shown with a green polygon. When $T_{1,1} = 356.25$ days, $T_{1,2} \in [1161.47, 1191.88]$ days, which are the minimum and maximum bounds of orbital periods at the second AP. This means that when we calculate minimum or maximum values, $T_{1,2}$ is either higher or lower than the given upper and lower bounds, so it becomes equal to the upper and lower bounds. Since $T_{1,2}$ is at the upper bound, there should be one more corner when $T_{1,2}$ starts decreasing and reaches its minimum value when $T_{1,1,\text{max}}$. This corner is at $T_{1,1} = 432.87$ days, when $T_{1,2} = 1191.88$ days calculated using Eq.~\eqref{eq:env_corner1_ap2_ap3}.

\subsubsection{Solution envelopes for the Geocentric Problem}
\begin{figure}[htbp!]
\begin{subfigure}{0.5\textwidth}
  \centering
\includegraphics[width=1\columnwidth]{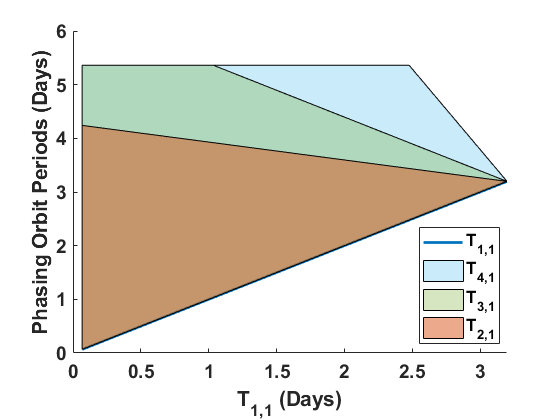}
  \caption{First AP }
  \label{fig:env_3ap_a}
\end{subfigure}
\begin{subfigure}{0.5\textwidth}
  \centering
\includegraphics[width=1\columnwidth]{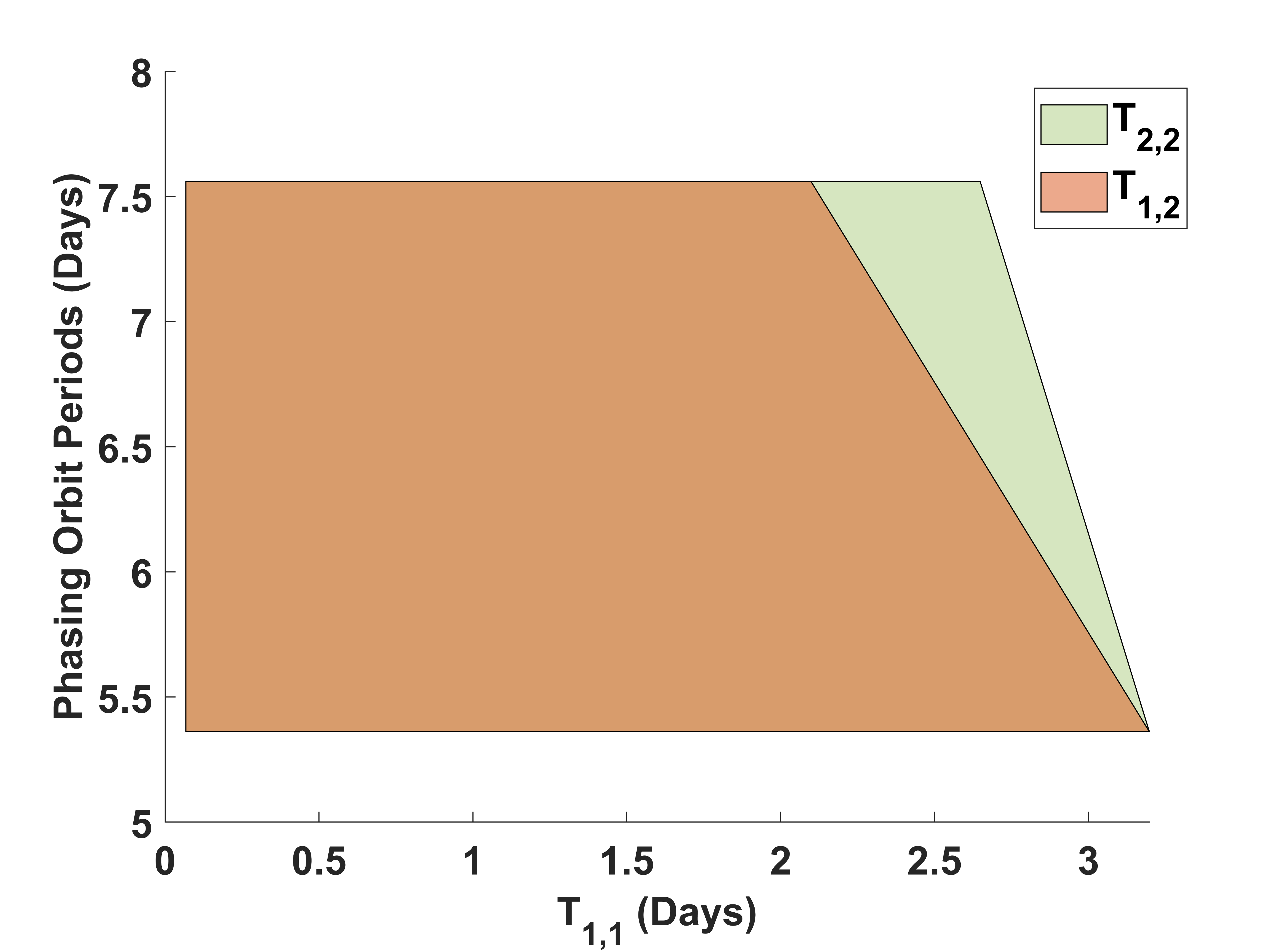}
  \caption{Second AP}
  \label{fig:env_3ap_b}
\end{subfigure}
\begin{subfigure}{0.5\textwidth}
    \centering
\includegraphics[width=1\columnwidth]{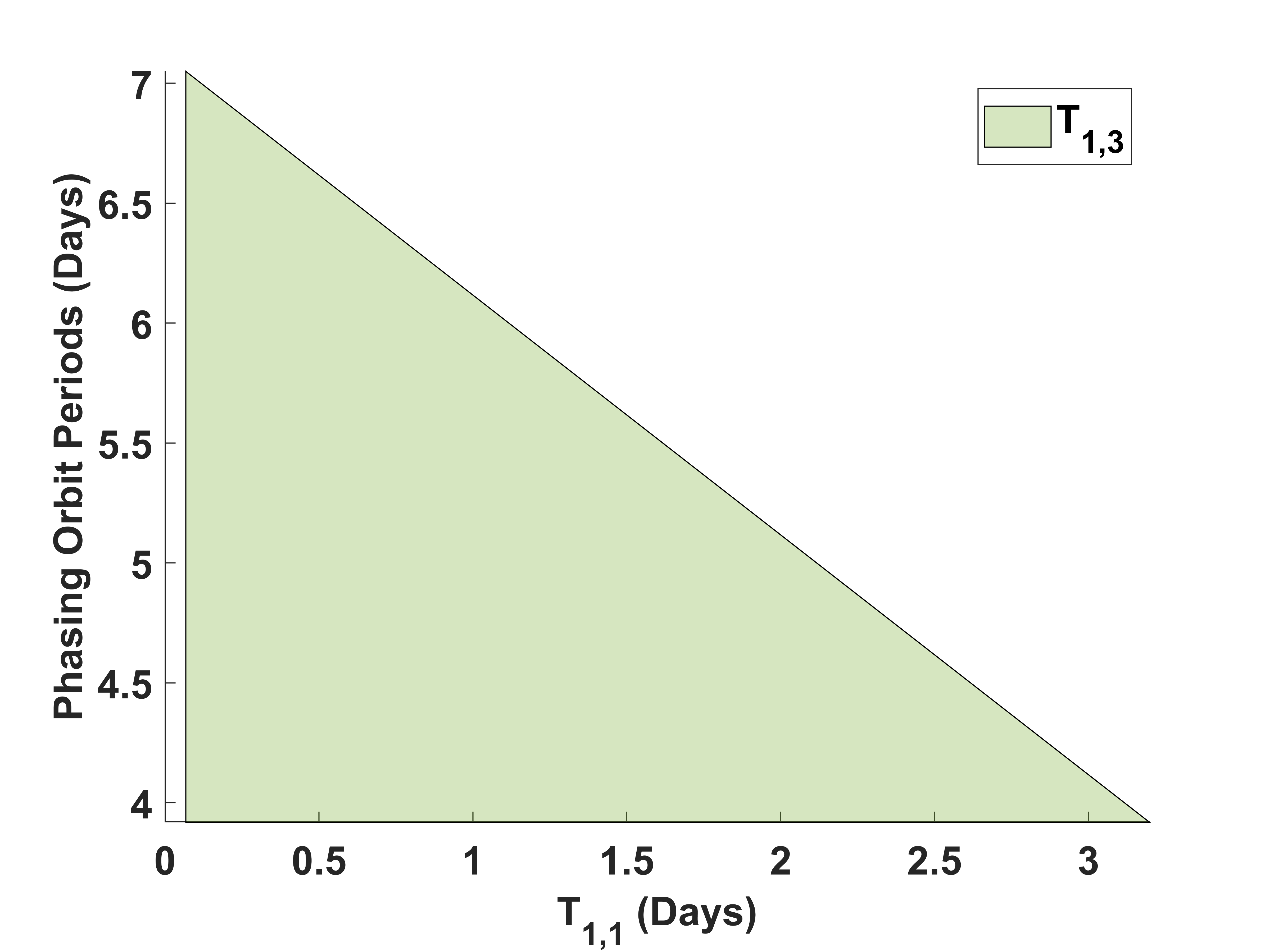}
  \caption{Third AP}
  \label{fig:env_3ap_c}
\end{subfigure}
\caption{Geocentric problem: solution envelopes with $N_{1,1} = N_{2,1} = N_{3,1} = N_{1,2} = N_{2,2}= 1$ and $N_{1,3} = 4$.}
\label{fig:env_3ap}
\end{figure}

Solution envelopes for the geocentric problem with three APs are generated using the provided analytical equations shown in Fig.~\ref{fig:env_3ap}. In Fig.~\ref{fig:env_3ap_a}, solution envelopes for the first AP phasing orbits are shown with orange ($T_{2,1}$), green ($T_{3,1}$), and light-blue ($T_{4,1}$) polygons. The light blue polygon has the largest area and the green one has the second-largest area, which is laid on top of the light blue. Orange is the triangle on top of both light-blue and green polygons. Again, we start determining the ranges of $T_{k,1}$ when $T_{1,1} = T_0 = 0.0675$ days for $k = 2, \cdots, 4$. Since all their minimum values are below $T_0$, they are set to $T_{0}$. For the upper bounds, $T_{3,1}$ and $T_{4,1}$ is at the upper bound with $T(\alpha_1 = 1) = 5.3616$ days. As a summary, $T_{2,1} \in [0.0675, 4.2413]$ days, $T_{3,1} \in [0.0675, 5.3616]$ days, and $T_{4,1} \in [0.0675, 5.3616]$ days, which are determined by using Eq.~\eqref{eq:env_min_ap1} and Eq.~\eqref{eq:env_max_ap1}. 
All the orbital period values converge to the rightmost point, where $T_{1,1,\text{max}} = T_{k,1} = 3.1978$ days determined from Eq.~\eqref{eq:env_maxT1_ap1}. Since $T_{3,1,\text{max}} = T_{4,1,\text{max}} = 5.3616$ days, for these solution envelopes, one more corner exists. These corner points corresponds to the maximum value of $T_{1,1}$, when $T_{3,1}$ or $T_{4,1}$ take their maximum values, which can be calculated by Eq.~\eqref{eq:env_corner1_ap1}. This corner is the point, when $T_{3,1,\text{max}} = 5.3616$ and $T_{1,1} = 1.0341$ days for the green polygon in Fig.~\ref{fig:env_3ap_a}, and $T_{4,1,\text{max}} = 5.3616$ days and $T_{1,1} = 2.4765$ days. In Fig.~\ref{fig:env_3ap_b}, the second AP phasing orbits' solution envelopes vs. $T_{1,1}$ are shown. Again, the corner points are determined for each, starting with the minimum and maximum value ranges, when $T_{1,1} = T_0$ using Eq.~\eqref{eq:solspace_2AP_3AP_eq}. Similar to the Earth-to-Dionysus problem second AP, $T_{1,2}$ and $T_{2,2}$ are between upper and lower bounds for the second AP phasing orbits, $T_{k,1}({\alpha_{k,1} = 1}) = 5.3616$ days and $T_{k,2}({\alpha_{k,2} = 1}) = 7.5606$ days. They converge to their minimum bound at $T_{1,1,\text{max}} = 3.1978$ days. One more corner exists for both orbit solution envelopes, where they start decreasing from their maximum value at $T_{1,1} = 2.0983$ days for $T_{1,2}$ and $T_{1,1} = 2.6481$ days for $T_{2,2}$ found using Eq.~\eqref{eq:env_corner1_ap2_ap3}.  In Fig.~\ref{fig:env_3ap_c}, the solution envelope for the third AP phasing orbit is shown with a green triangular area. Similar to the previous solution envelopes, we start with determining bounds of $T_{1,3} \in [3.9191, 7.0492]$ days when $T_{1,1} = T_0$. Since the upper bound of $T_{1,3} \neq T({\alpha_3 = 0})$, ($T(\alpha_3 = 0) > T(\alpha_3 = 1)$ due to the last $\Delta v$ increases the energy of the orbit to reach to the target orbit), there are only three corners for this solution envelope. The rightmost corner is when we have $T_{1,1\text{max}}$ and $T_{1,3}$ converges to its minimum value.

\section{Discussion on the Selection of the Base Solutions}\label{sec:discussion_base}
\subsection{Earth-to-Mars Problem}
\begin{figure}[!ht]
\begin{subfigure}{.5\textwidth}
  \centering
  \includegraphics[width=1\columnwidth]{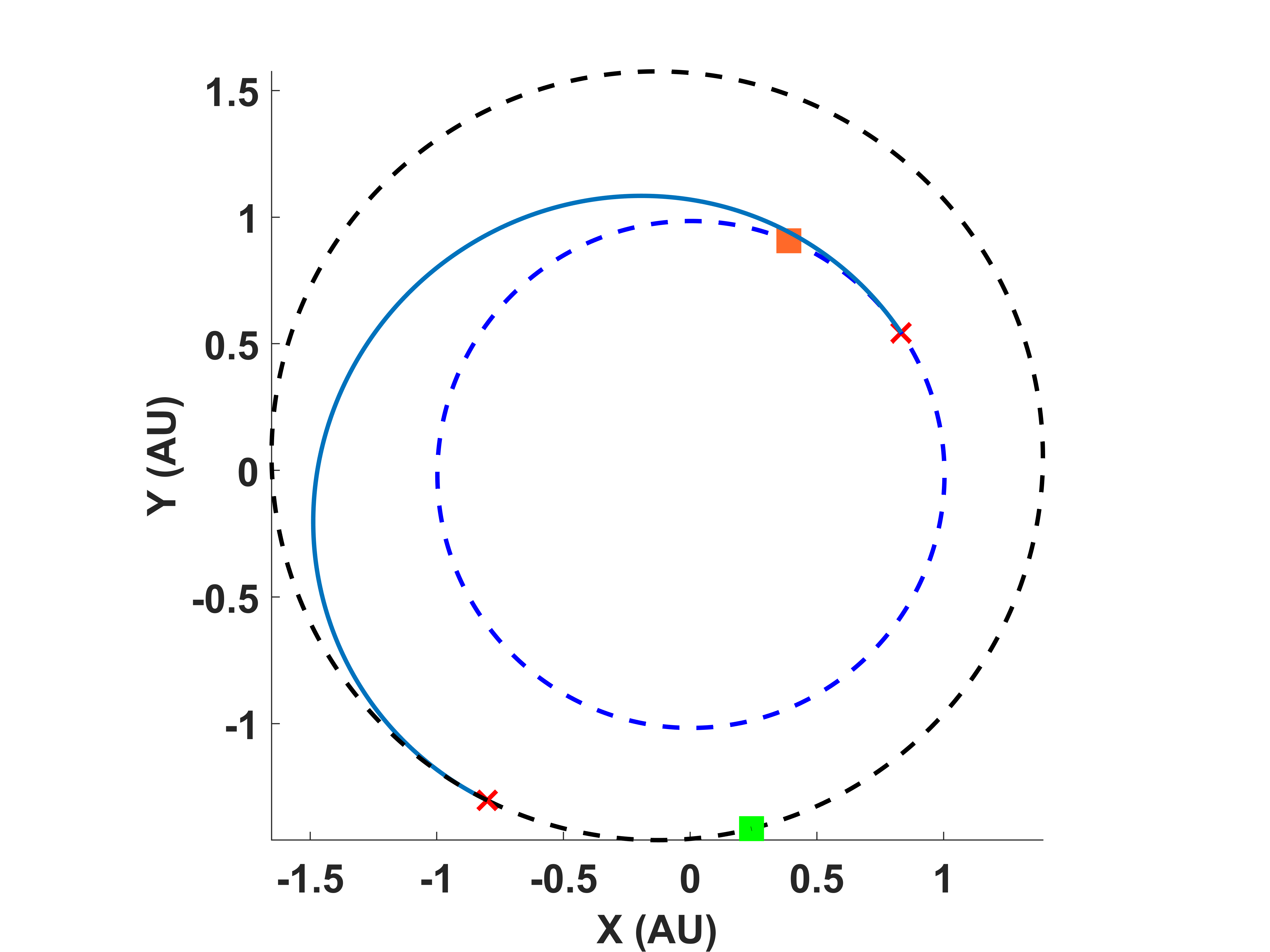}
  \caption{Two-impulse phase free}
  \label{fig:2imp_E2M}
\end{subfigure}%
\begin{subfigure}{.5\textwidth}
  \centering
  \includegraphics[width=1\columnwidth]{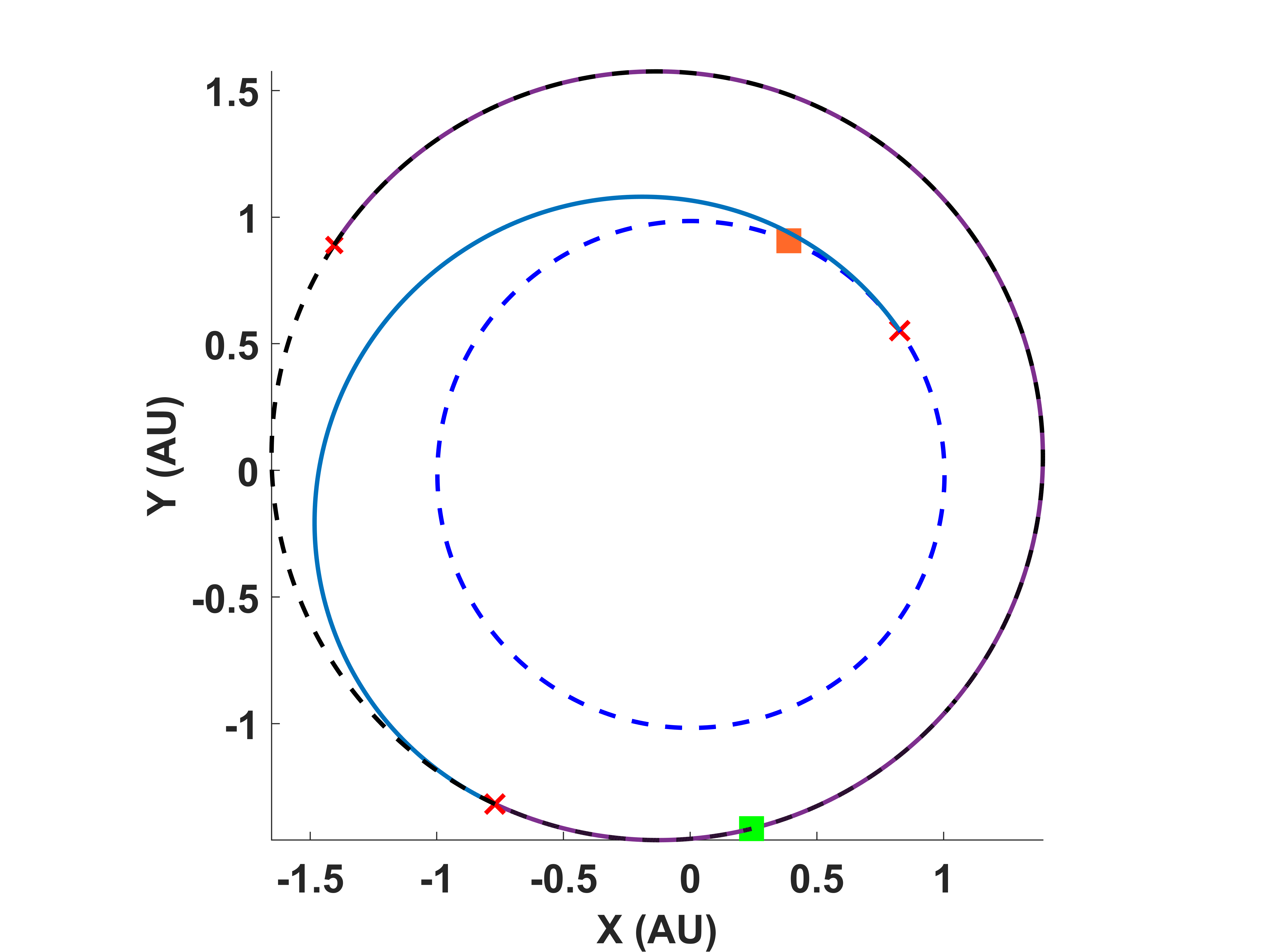}
  \caption{Three-impulse phase free}
  \label{fig:3imp_E2M}
\end{subfigure}
\begin{subfigure}{.5\textwidth}
  \centering
  \includegraphics[width=1\columnwidth]{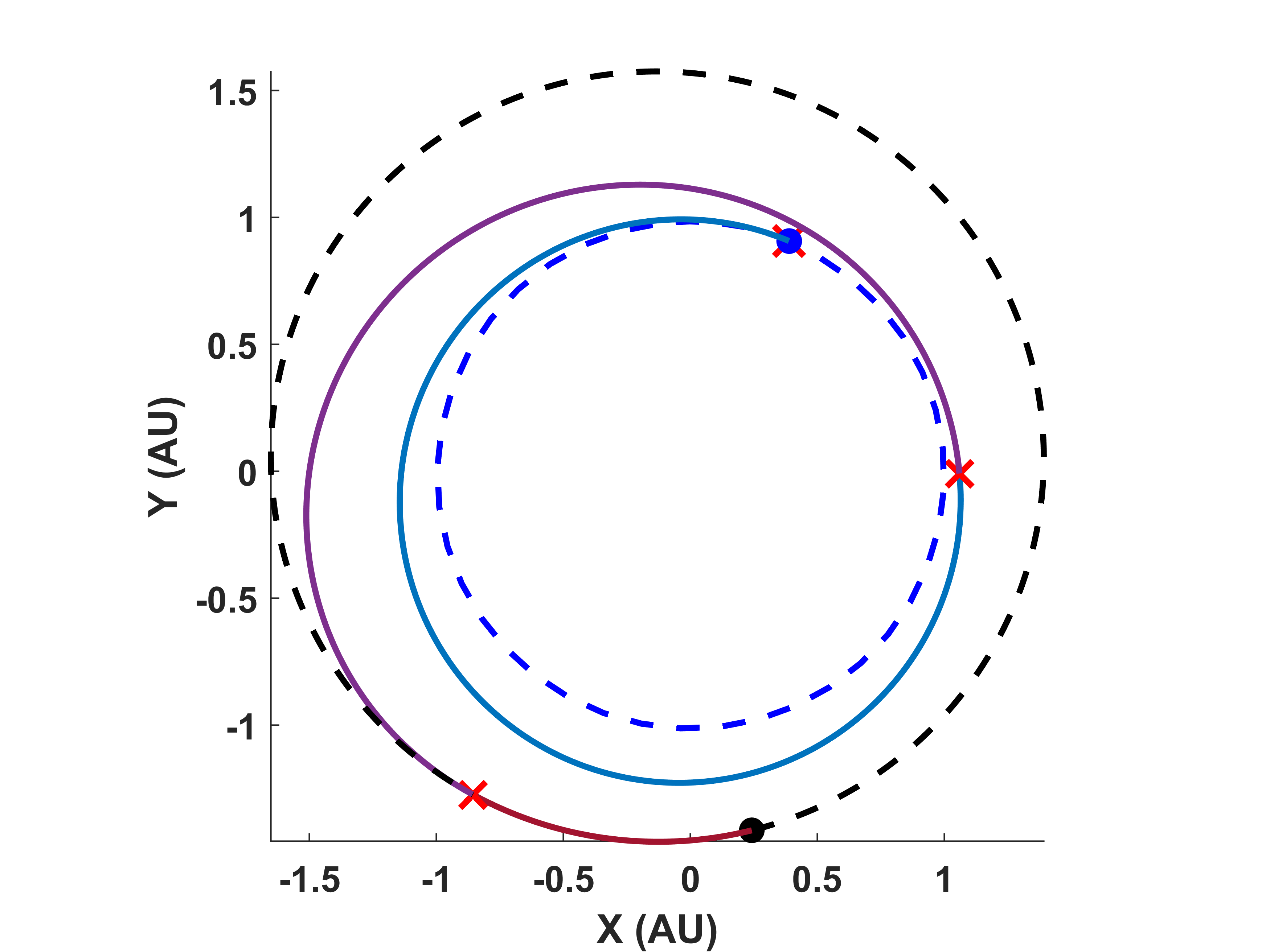}
  \caption{Time-fixed rendezvous}
  \label{fig:3imp_fixedtime_E2M}
\end{subfigure}
\begin{subfigure}{.55\textwidth}
  \centering
  \includegraphics[width=1\columnwidth]{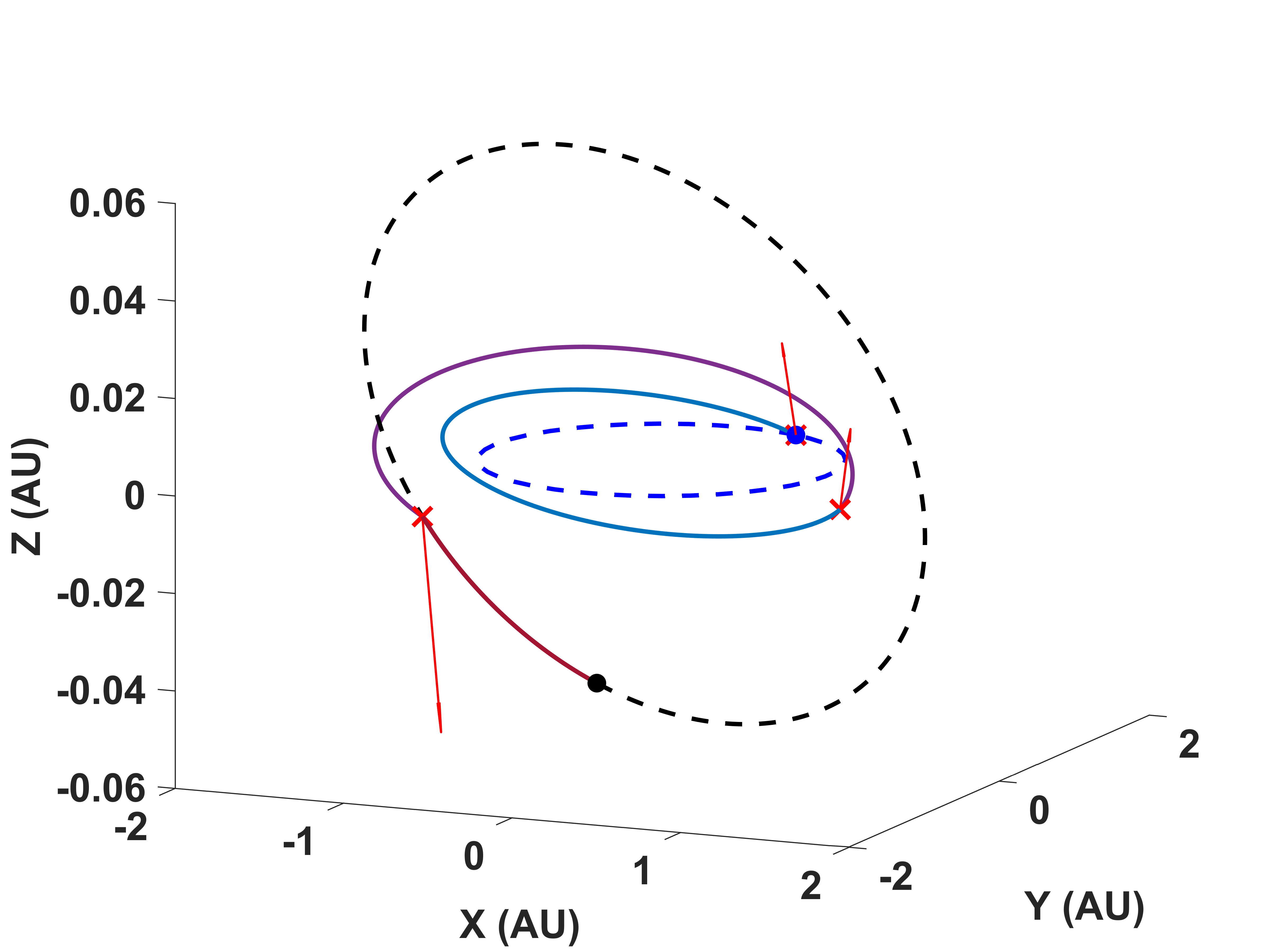}
  \caption{Time-fixed rendezvous three-dimensional trajectory}
  \label{fig:3imp_fixedtime_E2M_3dtraj}
\end{subfigure}
\caption{Solutions for the Earth-to-Mars problem.}
\label{fig:E2M}
\end{figure}

The first example is an Earth-to-Mars transfer from \cite{taheri_how_2020}. 
We proceed by generating the two- and three-impulse base solutions. Obtained solutions are given in Figures~\ref{fig:2imp_E2M} and \ref{fig:3imp_E2M} with $\Delta v_\text{total}$ values of $5.5865$ km/s and $5.5873$ km/s, respectively. The blue and black dot markers represent the departure and arrival points and the dashed blue and black orbits are the departure and arrival orbits, respectively. The $\Delta v_\text{total}$ values of both solutions are close, but the two-impulse solution has a lower value. Regarding flight times, the two-impulse solution occurs over 313.2430 days, whereas the three-impulse occurs over 801.4034 days. For a time-free transfer and time-free rendezvous, the selected base solution is the two-impulse one, with a significantly less mission time and a lower $\Delta v_\text{total}$ value.

Adding phasing orbits or coast arcs to this base solution can provide a rendezvous maneuver if time is not a constraint, i.e., time-free rendezvous. Mission designers can specify departure and arrival points as well as the required mission time using the base solution with the minimum $\Delta v$ value. On the other hand, the base solution provides a certificate of recovering the minimum $\Delta v$ value, if a predetermined phasing and a mission time exists. In Figures~\ref{fig:2imp_E2M} and \ref{fig:3imp_E2M}, the red and green square markers represent the locations of the Earth and Mars for a specific phasing with a mission time of 793 days. The time-fixed rendezvous solution is given in Figure~\ref{fig:3imp_fixedtime_E2M} with $\Delta v_\text{total} = 5.6109$ km/s. The solution is obtained by solving a hybrid PVT-based NLP problem, with an initial guess obtained by thrust acceleration continuation from the low-thrust solution.

When we consider the selection of a base solution for this time-fixed maneuver, the time of flight of the mission and the base solutions are compared. The total mission time is not feasible for both base solutions, since they require a large terminal coast arc before the first impulse. Note that the departure point (red square in Figs.~\ref{fig:2imp_E2M} and \ref{fig:3imp_E2M}) is ahead of the first impulse location. Therefore, the total time of flights becomes 701.2694 days and 1405.8939 days for the two- and three-impulse base solutions. Notice that for the two-impulse case, the total time of flight is shorter than the mission time. However, the constraint $TOF\geq T_0$ is not satisfied since $793-720.2694 = 72.7306 < T_0$ where, $T_0 = 365.25$ which is the Earth orbital period. Therefore, in terms of the time of flight, none of the base solutions are feasible. If the time of flight constraint of the mission is relaxed, the two-impulse base solution with the minimum $\Delta v_\text{total}$ value can be recovered. For example, we can add multiples of Mars's orbital periods to extend the mission time such that the correct phasing is preserved. Adding one orbital period to Mars does not change its location and velocity on its orbit. Then, $793+686.9658 = 1479.9658$ days would be the new mission time, with an available time of flight of $759.6963$ days. The available time can be spent on several phasing orbits. The mission time is sacrificed to recover the minimum $\Delta v$. 

On the other hand, if we keep the mission time at 793 days consistent with \cite{taheri_how_2020}, we cannot choose a two- or three-impulse base solution, since they are not time-feasible. Then, the optimal impulsive solution has to be found by using any of the existing impulsive trajectory optimization methods. However, the base solutions can be used as initial guesses. We note that the $\Delta v_\text{total}$ associated with the resulting time-constrained rendezvous maneuver is guaranteed to be higher than the two-impulse base solution, as the mission time constraint exists. 

Thus, the base solutions represent the lowest $\Delta v$ values for a transfer maneuver and can be used as a certificate of minimum-$\Delta v$ with the given mission time. Since the $\Delta v_\text{total} = 5.6108$ km/s for the time-fixed rendezvous case, the $\Delta v_\text{total}$ is sacrificed to satisfy the mission time. Please also note that the third impulse is located between the orbits of the Earth and Mars in Figure~\ref{fig:3imp_fixedtime_E2M}. For the same problem with the specified flight time of 793 days, the two-impulse Lambert solution results in $6.047$ km/s \cite{taheri_how_2020}. Therefore, the base solution provides a lower bound on the $\Delta v_\text{total}$ where the upper bound is the two-impulse Lambert solution. Even though the base solution $\Delta v_\text{total}$ value is not achievable due to the time constraint, one can determine lower and upper bounds on the $\Delta v_\text{total}$ with the proposed method.  

\subsection{Geocentric Circle-to-Circle Problem} \label{sec:geocircir}

\begin{figure}[!htbp]
\begin{subfigure}{.5\textwidth}
  \centering
  \includegraphics[width=1\columnwidth]{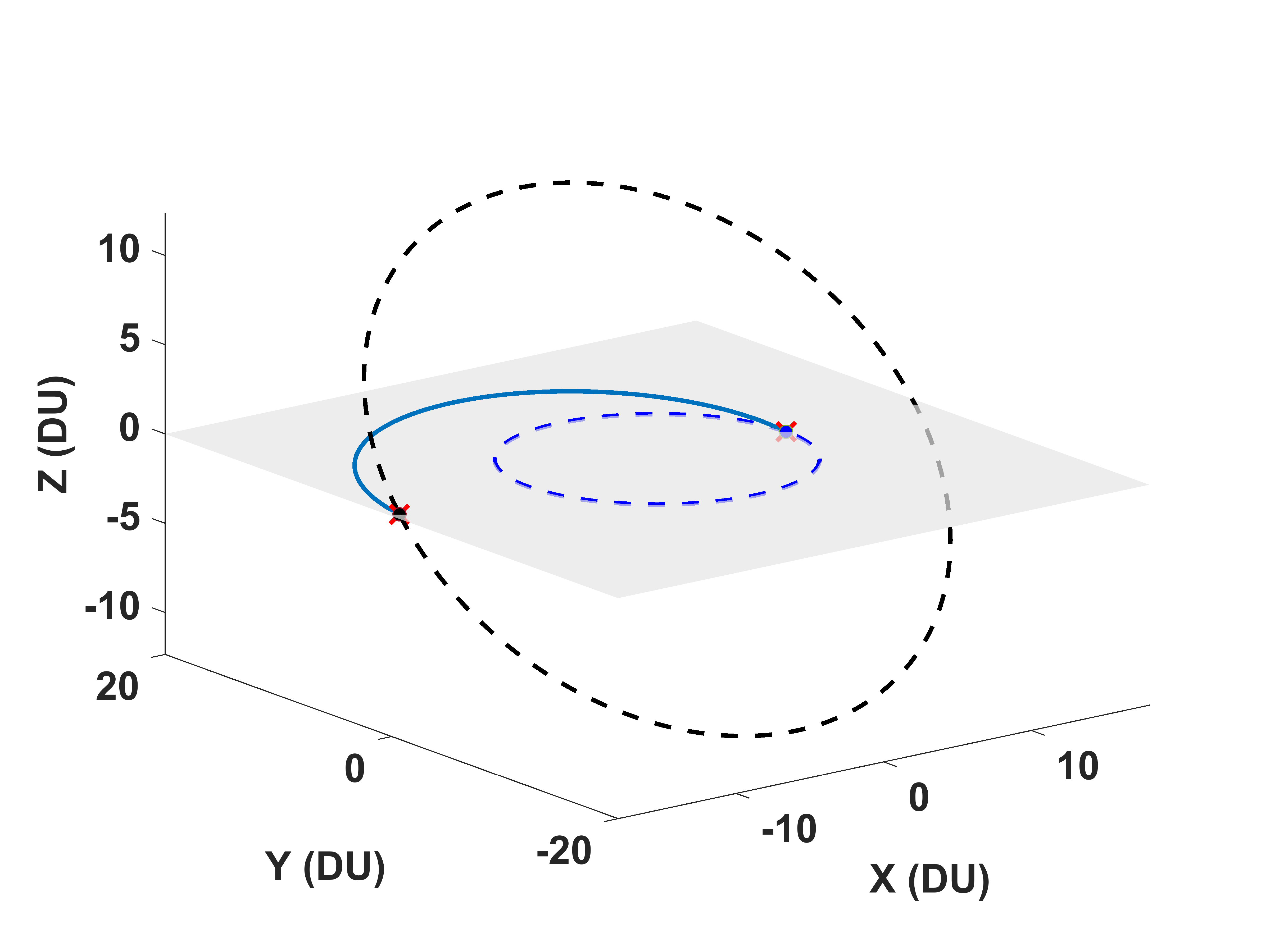}
  \caption{Two-impulse phase-free}
  \label{fig:2imp_c2c}
\end{subfigure}%
\begin{subfigure}{.5\textwidth}
  \centering
  \includegraphics[width=1\columnwidth]{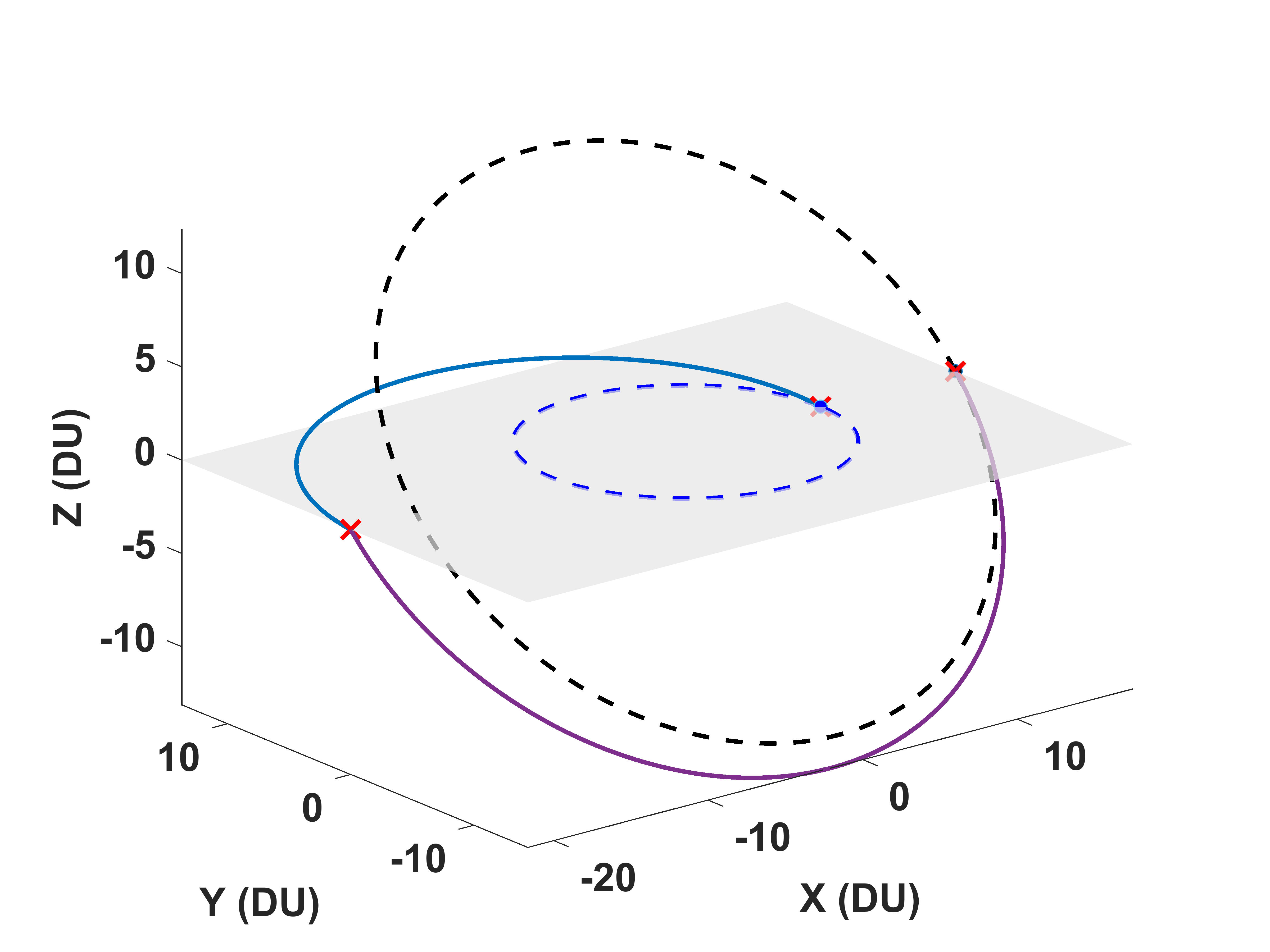}
  \caption{Three-impulse phase-free}
  \label{fig:3imp_c2c}
\end{subfigure}
\caption{Solutions for the circle-to-circle transfer }
\label{fig:circ2circ}
\end{figure}

The next example is a circle-to-circle transfer with a $45^\circ$ inclination change and $\beta = r_0/r_f = 0.5$. The selected base solution is the three-impulse one given in  Figure~\ref{fig:two_three_imp}, as it requires a less value of $\Delta v_\text{total}$. In terms of $\Delta v_\text{total}$, the two- and three-impulse base solutions require $1.7036$ km/s and $1.6853$ km/s, respectively. Regarding the total mission time, the two- and three-impulse solutions occur over 1.3905 days and 4.2767 days, respectively. The solutions are shown in Figure~\ref{fig:circ2circ}. There are two cases for base trajectory selection. The first case is when the time is free concerning the time-free transfer and time-free rendezvous, as shown in the flowchart given in Figure~\ref{fig:flowchart_base}. For those maneuvers, the first criterion is the $\Delta v_\text{total}$. The minimum $\Delta v_\text{total}$ corresponds to the three-impulse base solution. The times of flight for two- and three-impulse solutions are different; so, we proceed with the three-impulse base solution. 

The second case is when the time is fixed. Let $\Delta t$ denote the required mission time. In the time-fixed case, the first concern is the time feasibility. If the mission time, $\Delta t$, is greater than the summation of the three-impulse time of flight ($4.2767$ days) and the initial orbit orbital period ($1.5137$ days), i.e., $\Delta t > 5.7904$ days, the three-impulse solution is feasible in terms of time. Since the three-impulse solution also has the lowest $\Delta v_\text{total}$ value, the selected base solution would be the three-impulse one. In this case, the feasible mission time is connected to the $\Delta v$ optimality in the block diagram given in Figure~\ref{fig:flowchart_base}. If $\Delta t \in (2.9042,5.7904)$ days, then the two-impulse base solution is feasible in terms of time, but it does not have the lowest $\Delta v_\text{total}$ value. Therefore, the $\Delta v$ optimality is sacrificed. This is the branch on the flowchart where the time is feasible, but that branch is not connected to the $\Delta v$ optimality. If $\Delta t < 2.9042$, then none of the base solutions are feasible. Since the problem is fully constrained, and none of the base solutions are feasible, a two-impulse Lambert problem can be solved to obtain the best trajectory. In that case, the required $\Delta v_\text{total}$ is higher than any of the base solutions' $\Delta v_\text{total}$ value. Thus, the proposed method provides information on whether the minimum-$\Delta v$ value can be achieved with the given mission time. In addition, the mission time can be extended with the multiples of the orbital period of the target orbit such that the base solutions become feasible in terms of time. Another constraint for the problem can be a single impulse constraint for the engines. In that case, the mission time can be increased such that the single impulses are below a given threshold by introducing several phasing orbits to the trajectory.

Remark: For both two- and three-impulse base solutions in Fig.~\ref{fig:circ2circ}, there are identical solutions in terms of $\Delta v_\text{total}$. They are the mirror image of the solutions about the $y-z$ plane. Depending on the location of the initial point on the departure orbit and the final point on the arrival orbit, we can have different values of $t_{c,1}$ and $t_{c_2}$ on all the base solutions. Thus, for these equal-$\Delta v_\text{total}$ cases, the lowest terminal coast time trajectory can be chosen for the time-fixed maneuver cases. 
\section{Conclusion} \label{sec:conclusion}
Impulse anchor positions (APs) correspond to the locations of the impulses on time-free, phase-free impulsive base solutions. 
The generation of infinitely many iso-impulse (equal-$\Delta v$) solutions is extended to problems with time-free, phase-free, three-impulse base solutions. Due to the existence of two- and three-impulse base solutions, a selection has to be made among the two options. A selection process for two- or three-impulse base solutions is proposed based on $\Delta v$-optimality and time-feasibility criteria. A circle-to-circle and an Earth-to-Mars maneuver examples are presented to discuss the selection process for the base solution. 
Steps to perform the $\Delta v$-allocation simultaneously at all potential impulse APs are presented. The problem of Earth to asteroid Dionysus and a geocentric problem are considered to demonstrate simultaneous $\Delta v$-allocation at two and three impulse APs and example iso-impulse solution families are presented. Feasible solutions for each problem with two and three impulse APs are analyzed in terms of the total number of revolutions on phasing orbits at each AP to generate feasible solution space plots. The feasibility of a solution family refers to having infinitely many iso-impulse solutions. Solution envelopes are generated and the analytical relations for determining lower and upper bounds on the orbital periods of the phasing orbits are derived. The feasibility of having infinitely many minimum-$\Delta v$ solutions and generating the associated solution space envelopes is explained for solutions with two and three impulse APs. 

Results show that it is possible to divide the impulses and perform $\Delta v$ allocation at multiple impulse APs simultaneously if any of the AP impulse values are beyond the ability of the engine to produce a maximum $\Delta v$ threshold. It is shown that the $\Delta v$-allocation method is applicable to three important classes of orbital maneuvers: 1) time-free transfer, 2) time-free rendezvous, and 3) time-fixed rendezvous. For time-fixed rendezvous problems, the total time associated with the base solution immediately shows whether it is possible or not to recover the minimum-$\Delta v$ solution through $\Delta v$ allocation. This offers a unique method for determining the time-feasibility of minimum-$\Delta v$ maneuvers between two three-dimensional orbits and based on a specified total mission time set by mission designers. The base solution $\Delta v$ value corresponds to the lower bound for time-fixed rendezvous maneuvers. When a base solution cannot be selected for time-fixed rendezvous maneuvers, the generated time-free, phase-free solutions are good initial guesses for impulsive trajectory optimization. 
For transfer maneuvers between co-axial, inclined circular orbits, we provide a plot based on inclination change and orbit ratio to determine the $\Delta v$ optimality based on the number of impulses on the base solutions. We have provided a detailed partitioning of the solution space to facilitate the process of determining the base solutions. For transfers between non-coaxial, elliptical, inclined transfers, it is concluded that the selection process of base solutions has two main features: time of flight and the $\Delta v$ optimality. 
We have shown that the provided time-feasibility equation can be used to determine the total possible number of feasible solutions analytically and the maximum number of phasing orbits that can be added at any AP. Each feasible solution is expanded with different combinations of numbers of revolutions on each phasing orbit and each solution family has infinitely many solutions. Solution envelopes provide information on these solution families in terms of ranges of orbital periods, which can be related to the required $\Delta v$. A remarkable result is that all solutions can be categorized into four layers: 1) the base solutions, 2) feasible solution spaces, 3) solution families, and 4) solution envelopes.

\appendix
\setcounter{equation}{0}\renewcommand\theequation{A\arabic{equation}}
\section*{Appendix: Solution Envelope Equations}\label{sec:sol_env_app}
\subsection{Corners at the First AP}
 We determine the leftmost corners ``1'' and ``2'' as shown in Fig.~\ref{fig:solenv1}, where $T_{1,1} = T_0$ for the first AP. The minimum values of orbital periods (corner ``1'' in Fig.~\ref{fig:solenv1}) are
\begin{equation}
    \begin{aligned}
        T_{n_{p,1},1}^{\text{min}} &= \frac{TOF_p - N_{1,1}T_{1,1}^{\text{min}} - C_1}{\sum_{k=2}^{n_{p,1}}N_{k,1}}, & & 
        T_{s,1}^{\text{min}} = \frac{TOF_p - N_{1,1}T_{1,1}^{\text{min}} - \sum_{k=s+1}^{n_{p,1}}N_{k,1}T_{k,1}^{\text{max}} - C_1}{\sum_{k=2}^{s}N_{k,1}},
    \end{aligned}
    \label{eq:env_min_ap1}
\end{equation}
where $s=2,\cdots,n_{p,1}-1$, and $C_1 = \sum_{k=1}^{n_{p,2}}N_{k,2}T_{k,2}^{\text{max}} + \sum_{k=1}^{n_{p,3}}N_{k,3}T_{k,3}^{\text{max}}$. If any of the $T_{s,1,\text{min}}$ values is less than $T_0$, then we set them to $T_0$. Then, we obtain the maximum orbital period values (corner ``2'' in Fig.~\ref{fig:solenv1}) by
\begin{equation}
    \begin{aligned}
        T_{m,1}^{\text{max}} = \frac{TOF_p - \sum_{k=1}^{m-1}N_{k,1}T_{k,1}^{\text{min}} - C_2}{\sum_{k=m}^{n_{p,1}}N_{k,1}},
    \end{aligned}
    \label{eq:env_max_ap1}
\end{equation}
where $m=2,\cdots,n_{p,1}$ and $C_2 = \sum_{k=1}^{n_{p,2}}N_{k,2}T_{k,2}^{\text{min}} + \sum_{k=1}^{n_{p,3}}N_{k,3}T_{k,3}^{\text{min}}$. If any of these maximum values are larger than the upper bound, $T_{k,1}(\alpha_{k,1} = 1)$, we set them equal to the upper bound. The maximum value of $T_{1,1}$ is when all the orbital periods of phasing orbits are equal (corner ``3'' in Fig.~\ref{fig:solenv1}), which can be written as, 
\begin{equation}
    T_{1, 1}^{ \text{max}}=\frac{TOF_p-C_2}{\sum_{k=1}^{n_{p,1}} N_{k,1}}.
    \label{eq:env_maxT1_ap1}
\end{equation}
\subsection{Corners at the Second and Third APs}
The maximum and minimum values can be determined when $T_{1,1} = T_0$, which are the corners ``1'' and ``2'' in Fig.~\ref{fig:solenv1},
\begin{equation}
    \begin{aligned}
        T_{m,i}^\text{max} &= \frac{TOF_p - \sum_{k=1}^{m-1}N_{k,i}T_{k,i}^{\text{min}} - C_{3,q}}{\sum_{k=m}^{n_{p,i}}N_{k,i}}, \\
        T_{n_{p,i},i}^{\text{min}} &= \frac{TOF_p-N_{1,1}T_{1,1}^{\text{min}}- C_{4,q}}{\sum_{k=1}^{n_{p,i}}N_{k,i}},\\
        T_{s,i}^{\text{min}} &= \frac{TOF_p -N_{1,1}T_{1,1}^{\text{min}}- \sum_{k=s+1}^{n_{p,2}}N_{n_{p,i},i}T_{n_{p,i},i}^{\text{max}} - C_{4,q}}{\sum_{k=1}^{s}N_{k,i}},   
    \end{aligned} \quad \text{for} \quad i = 2,3,
    \label{eq:solspace_2AP_3AP_eq}
\end{equation}
where $m = 2,\cdots,n_{p,i}$, and $s=1,\cdots,n_{p,i}-1$ and 
\begin{align*}
    C_{3,q} & = \sum_{k=1}^{n_{p,1}}~N_{k,1}~T_{k,1}^{\text{min}} + \sum_{k=1}^{n_{p,q}}~N_{k,q}~T_{k,q}^{\text{min}}, \\
    C_{4,q} & = \sum_{k=2}^{n_{p,1}}~N_{k,1}~T_{k,1}^{\text{max}} + \sum_{k=1}^{n_{p,q}}~N_{k,q}~T_{k,q}^{\text{max}},
\end{align*}
with $q = 3$ if $i = 2$ or vice versa. Note that the minimum and maximum values we use in $C$ constants are the upper and lower bounds those orbital periods possibly can take.
To calculate the maximum value of the first orbit, shown as corner ``3'' in Fig.~\ref{fig:solenv1}, when $i \in \{2,3\}$, we have
\begin{equation}
    T_{1, i}^{\text{max}}=\frac{TOF_p-C_{3,q}}{\sum_{k=1}^{n_{p,i}} N_{k,i}}.
    \label{eq:env_maxT1_ap23}
\end{equation}

\subsection{Additional Corners}
If any of the resultant values are not between the bounds, then the upper or lower bound is set for the corresponding orbital period value. In the case of the upper bound, there is a possibility of having one additional corner in the solution envelope. That corner is at the point when $T_{1,1}$ is maximized when the other orbital period is at its maximum bound. The corner is shown in Fig.~\ref{fig:solenv2} as corner ``4''. The other case is when the minimum value of the orbital periods is higher than its lower bound, when $T_{1,1} = T_0$. Therefore, there exists a corner where the orbital period value decreases and takes its minimum value at a point when $T_{1,i} \neq T_0$ which is the corner ``5'' in Fig.~\ref{fig:solenv2}.

When there is a possibility of having a corner point due to an orbital period variable being on the upper bound, corner ``4'' in Fig.~\ref{fig:solenv2}, one can determine this point by maximizing $T_{1,1}$ when $T_{j,1}$ is at the maximum value for any $j = 2,\cdots,n_{p,1}$, 
\begin{equation}
\begin{aligned}
    T_{1,1}^c &= \frac{TOF_p -\sum_{k=j}^{n_{p,1}}N_{k,1}~T_{k,1}^{\text{max}} -C_2}{\sum_{k=1}^{j-1}N_{k,1}}.
    \end{aligned}
    \label{eq:env_corner1_ap1}
\end{equation}

If an additional corner is possible due to $T_{j,1,\text{min}}$ being greater than the lower bound, we have 
\begin{equation}
\begin{aligned}
     T_{1,1} = T_{j,1}^{\text{min}} &= \frac{TOF_p -C_1}{\sum_{k=1}^{j}N_{k,1}}, \quad \text{for}\quad j = n_{p,1}, \\
     T_{1,1} = T_{j,1}^{\text{min}} &= \frac{TOF_p -\sum_{k=j+1}^{n_{p,1}}N_{k,1}~T_{k,1}^{\text{max}} -C_1}{\sum_{k=1}^{j}N_{k,1}}, \quad \text{for}\quad j \neq n_{p,1}.
    \end{aligned}
    \label{eq:env_corner2_ap1}
\end{equation}

At the second or third AP, similar relations can be derived for corner points for the case when any of the variables are at the upper bound: 

\begin{equation}
\begin{aligned}
    T_{1,1}^{\text{max}} &= \frac{TOF_p -\sum_{k=j}^{n_{p,2}}N_{k,2}~T_{k,2}^{\text{max}} - \sum_{k=1}^{j-1}N_{k,2}T_{k,2}^{\text{min}}-\sum_{k=1}^{n_{p,3}}N_{k,3}~T_{k,3}^{\text{min}}}{\sum_{k=1}^{n_{p,1}}N_{k,1}},
    \end{aligned}
    \label{eq:env_corner1_ap2_ap3}
\end{equation}
where the subscripts $(., 2)$ and $(.,3)$ can be swapped to determine the corner points for the third AP. 

There can be a corner point due to  $T_{j,2}^{\text{min}}$ or $T_{j,3}^{\text{min}}$ taking values greater than the lower bound for any $j$ which shown as corner ``5'' in Fig.~\ref{fig:solenv2}. If such a case is possible, the additional corner is at the point of minimum value of $T_{1,1}$ when $T_{j,2}^{\text{min}}$ or $T_{j,3}^{\text{min}}$, which are the lowest values they can get, regardless of the value of $T_{1,1}$. Unlike the first AP case, $T_{1,1} \neq T_{j,2}^{\text{min}}$ or  $T_{1,1} \neq T_{j,3}^{\text{min}}$ because $T_{1,1}$ is not related to other AP phasing orbits through total $\Delta v$. At this corner point, the minimum value of $T_{1,1}$ is determined by 
 \begin{equation}
 \begin{aligned}
      T_{1,1} &= \frac{TOF_p - \sum_{k=1}^{j}N_{k,2}~T_{k,2}^{\text{min}} - \sum_{k=j+1}^{n_{p,2}}N_{k,2}~T_{k,2}^{\text{max}}-T_{n_{p,1},1}^{\text{max}} - \sum_{k=1}^{n_{p,3}}N_{k,3}~T_{k,3}^{\text{max}}}{\sum_{k=1}^{n_{p,1}-1}N_{k,1}}, \quad \text{for} \quad j \neq n_{p,i}, \\
       T_{1,1} &= \frac{TOF_p - \sum_{k=1}^{j}N_{k,2}~T_{k,2}^{\text{min}} -T_{n_{p,1},1}^{\text{max}} - \sum_{k=1}^{n_{p,3}}N_{k,3}~T_{k,3}^{\text{max}}}{\sum_{k=1}^{n_{p,1}-1}N_{k,1}}, \quad \text{for} \quad j = n_{p,i}.
     \end{aligned}
     \label{eq:env_corner2_ap2_ap3}
 \end{equation}
\section*{Acknowledgments}
Authors acknowledge the NASA Alabama EPSCoR Research Seed Grant for partially supporting this research.

\bibliography{ref_export}

\end{document}